\documentclass[a4paper,12pt,reqno]{amsart}

\usepackage{amssymb}

\input epsf.tex
\newdimen\xsize
\newdimen\oldbaselineskip
\newdimen\oldlineskiplimit
\def\restorelineskip{\baselineskip=\oldbaselineskip%
\lineskiplimit=\oldlineskiplimit}
\def\putm[#1][#2]#3{
\hbox{\vbox to 0pt{\parindent=0pt%
\vskip#2\xsize\hbox to0pt{\hskip#1\xsize $#3$\hss}\vss}}}%
\long\def\Line#1{\hbox to \hsize{#1}}
\def\putt[#1][#2]#3{
\vbox to 0pt{\noindent\hskip#1\xsize\lower#2\xsize%
\vtop{\restorelineskip#3}\vss}}

\makeatletter
\def\xbig[#1]#2{{\hbox{$\m@th\left#2\vbox to#1\xsize{}%
\right.\n@space$}}}
\makeatother
\def\xlar[#1]#2{%
\smash{\mathop{ \hbox to #1\xsize{\leftarrowfill}}\limits^{#2}}}

\def\xline[#1]{\hbox to #1\xsize{\leaders\hrule\hfill}}

\DeclareFontFamily{U}{rsf}{\skewchar\font'177}%
\DeclareFontShape{U}{rsf}{m}{n}{<-6>rsfs5<6-8>rsfs7<8->rsfs10}{}%
\DeclareFontShape{U}{rsf}{b}{n}{<-6>rsfs5<6-8>rsfs7<8->rsfs10}{}%
\DeclareMathAlphabet\RSFS{U}{rsf}{m}{n}
\SetMathAlphabet\RSFS{bold}{U}{rsf}{b}{n}
\let\scr=\rfs
\def\mib#1{\boldsymbol{#1}}
\def\sf#1{{\mathsf{#1}}}

\def\slsf{\slshape \sffamily }

\def\msmall#1{\mathchoice{\hbox{\small$\displaystyle {#1}$}}{#1}{#1}{#1}}

\def\bfu{{\mathbf{U}}}

\let\3=\ss
\def\sss{{\mathbb S}}
\def\sph{{\mathbb S}}
\let\ss=\sss

\def\cc{{\mathbb C}}

\def\jj{{\mathbb J}}
\def\rr{{\mathbb R}}
\def\nn{{\mathbb N}}

\def\pp{{\mathbb P}}
  
\def\zz{{\mathbb Z}}

\def\adyn{\sf{1}}
\def\arg{{\sf{Arg}}}

\def\loc{{\sf{loc}}}
\def\st{_{\mathsf{st}}}

\def\deg{\sf{deg}\,}

\def\dim{\sf{dim}\,}

\def\endo{\sf{End}} \def\End{\sf{End}}

\def\ext{\sf{ext}}

\def\ind{\sf{ind}}
\def\bindo{\sf{ind}^b_0}

\def\exp{\sf{exp}}

\def\sfh{\sf{H}}

\def\id{\sf{Id}}
\def\im{\sf{Im}\,}

\def\re{\sf{Re}\,}

\def\lim{\mathop{\sf{lim}}}

\def\min{\sf{min}}

\def\pr{\sf{pr}}

\def\bfc{{\mib{c}}}

\def\bfu{{\mib{U}}}

\def\eps{\varepsilon}

\def\<{\langle}\let\la=\<
\def\>{\rangle}\let\ra=\>
 \let\bs=\bss 

\def\db{\partial} 
\def\dbar{{\barr\partial}}

\def\ddef{\mathrel{{=}\raise0.3pt\hbox{:}}}
\def\deff{\mathrel{\raise0.3pt\hbox{\rm:}{=}}}
\def\inv{^{-1}}

\def\fraction#1/#2{\mathchoice{{\msmall{ \frac{#1}{#2}}}}%
{{ \frac{#1}{#2 }}}{{#1/#2}}{{#1/#2}}}
\def\norm#1{\left\Vert{#1}\right\Vert}
\def\half{\frac{1}{2}}
\def\le{\leqslant}
\def\vph{^{\mathstrut}}

\def\scirc{\mathop{\mathchoice{\hbox{\small$\circ$}}{\hbox{\small$\circ$}}%
{{\scriptscriptstyle\circ}}{{\scriptscriptstyle\circ}}}}
\def\longpoints{\leaders\hbox to 0.5em{\hss.\hss}\hfill \hskip0pt}
\def\stateskip{\smallskip}
\def\state#1. {\stateskip\noindent{\bf#1. }} 
\def\statep#1. {\stateskip\noindent{\bf#1 }} 
\def\proof{\state Proof. }

\def\Chi{\raise 2pt\hbox{$\chi$}}

\def\ie{\hskip1pt plus1pt{\sl i.e.\/,\ \hskip1pt plus1pt}}

\def\sli{{\sl i)} } 
\def\slii{{\sl i$\!$i)} } 
\def\sliii{{\sl i$\!$i$\!$i)} }
\def\sliv{{\sl i$\!$v)} }

\def\barr#1{\mskip1mu\overline{\mskip-1mu{#1}\mskip-1mu}\mskip1mu}
\def\vect{{\mathrm{v}}}
\def\wect{{\mathrm{w}}}

\def\sft{{\sf t}}

\def\Chi{\raise 2pt\hbox{$\chi$}}

\let\phI=\phi\let\phi=\varphi\let\varphi=\phI
\def\bfmu{{\boldsymbol\mu}}%
\def\bfchi{{\boldsymbol\chi}}%
\def\scrw{{\scr W}}

\def\calc{{\cal C}}

\def\cald{{\cal D}}

\def\calk{{\cal K}}

\def\calo{{\cal O}}

\def\eps{\varepsilon}

\def\bs{\backslash}

\def\d{\partial}
\def\dbar{{\barr\partial}}
\def\1{{1\mkern-5mu{\rom l}}}

\def\ge{\geqslant}

\def\inv{^{-1}}

\let\wt=\widetilde

\def\fraction#1/#2{\mathchoice{{\msmall{ \frac{#1}{#2}}}}%
{{ \frac{#1}{#2 }}}{{#1/#2}}{{#1/#2}}}
\def\half{{\fraction1/2}}
\def\le{\leqslant}
\def\leq{\leqslant}
\def\vph{^{\mathstrut}}

\def\ti#1{{\tilde{#1}}}

\def\qed{\ \ \hfill\hbox to .1pt{}\hfill\hbox to .1pt{}\hfill $\square$\par}

\def\comment#1\endcomment{}

 \belowdisplayskip=\abovedisplayskip

\def\lineeqqno(#1){\hfill\llap{\vbox to 10pt%
{\vss\begin{align} \eqqno(#1)\end{align}\vss}}\vskip1pt}

\advance\headheight 1.2pt

\def\ShowwLLabel#1{}

\def\thechpt{\Roman{chpt}}
\def\newchapt[#1]#2{\newpage%
\refstepcounter{chpt}\setcounter{subsection}{0}%
\setcounter{thm}{0}\setcounter{defi}{0}%
\setcounter{rema}{0}\setcounter{exrc}{0}%
\renewcommand{\thesubsection}{\thechpt.\arabic{subsection}}%
\section*{\begin{center}\huge \bf Chapter \thechpt\\
#2 \end{center}}\label{#1}%
\ \smallskip%
\addcontentsline{toc}{part}{Chapter \thechpt. #2}%
\markboth{Chapter \thechpt}{#2}%
}

\def\newsect[#1]#2{\refstepcounter{section}\setcounter{equation}{0}%
\renewcommand{\thesubsection}{\arabic{section}.\arabic{subsection}}%
\section*{\arabic{section}.
#2}\vspace{-20pt}\label{#1}\vspace{20pt}%
\markboth{Section \arabic{section}}{#2}}

\def\newlect[#1]#2{\refstepcounter{section}%
\renewcommand{\thesubsection}{\arabic{section}.\arabic{subsection}}%
\section*{Lecture \arabic{section}\\
#2}\label{#1}%
\markboth{Lecture \arabic{section}}{#2}}

\def\newprg[#1]#2{\refstepcounter{subsection}%
\subsection*{{\thesubsection.\ #2}} \label{#1}%
}

\def\newappx[#1]#2{%
\refstepcounter{appx}\setcounter{section}{0}%
\renewcommand{\thesubsection}{A\arabic{appx}.\arabic{subsection}}%
\section*{Appendix \arabic{appx}\\ #2}
\label{#1}%
\markboth{Appendix A\arabic{appx}}{#2}
}

\newtheorem{thm}{Theorem}
   \def\newthm#1{\begin{thm}\label{#1}}
   
\newtheorem{lem}{Lemma}[section]
   \def\newlemma#1{\begin{lem} \label{#1}}

\newtheorem{prop}[lem]{Proposition}
   \def\newprop#1{\begin{prop}\label{#1}}

\newtheorem{corol}[lem]{Corollary}
   \def\newcorol#1{\begin{corol} \label{#1}}

\newtheorem{defi}{Definition}[section]
   \def\newdefi#1{\begin{defi} \label{#1}\rm }

\newtheorem{exmp}{Example}[section]
   \def\newexmp#1{\begin{exmp} \label{#1}\rm }

\newtheorem{exrc}{Exercise}[subsection]
   \def\newexrc#1{\begin{exrc} \label{#1}\rm }

\newtheorem{quest}{Question}
   \def\newquest#1{\begin{quest} \label{#1}\rm }

\newtheorem{rema}{Remark}[section]
   \def\newrema#1{\begin{rema} \label{#1}\rm }

\newtheorem{nnrema}{Remark}
   \def\newrema#1{\begin{rema} \label{#1}\rm }

\def\eqqno(#1){\label{(#1)}}
\def\eqqref(#1){(\ref{(#1)})}

\let\3=\ss
\def\sss{{\mathbb S}}
\let\ss=\sss

\numberwithin{equation}{section}
\usepackage[dvips]{graphicx}
\usepackage[dvipsnames]{xcolor}  
\usepackage{colortbl}

\def\db{\partial}
\def\re{\sf{Re}\,}

\def\cal#1{{\mathcal{#1}}}

\title{Positivity of Boundary Intersections of Complex \\ 
Curves  and Adjunction Formula}
\author{Ivashkovych Serge}
\address{Universit\'e de Lille, UFR de Math\'ematiques, 59655 Villeneuve
d'Ascq, France.}
\email{serge.ivashkovych@univ-lille.fr}

\thanks{* Partially supported by the  R-CDP-24-004-C2EMPI project.}
\date{\today}

\begin{document}
\begin{abstract}
One of the goals of this paper is to prove that the index of intersection of two complex curves 
in a two-dimensional complex manifold tangent to each other at a common boundary point 
is positive. This is achieved via the construction of a totally real surface such that the 
curves in question are attached to it by some parts of their boundaries and then 
defining a certain ``boundary intersection index'' of two complex ``half-disks'' with their 
edges on a totally real surface. We prove that this index is always positive. This second 
result holds true, more generally, for complex curves in a two dimensional almost
complex manifold with boundaries on totally real submanifold, but to our best knowledge 
is new even for integrable structures, unless the totally real surface in question is supposed 
to be real analytic. We also formulate and prove the Adjunction Formula for complex
curves with boundaries on totally real submanifolds in an almost complex manifold of dimension
four.
\end{abstract}

\maketitle

\setcounter{tocdepth}{1}
\tableofcontents

\newsect[INT]{Introduction and Statement of Main Results}

\newprg[INT.b-cusp]{Boundary intersection index}

Let $(X,J)$ be an almost complex manifold, $p_0\in X$ a point  and $W\ni p_0$ a germ of a $J$-totally 
real submanifold of $X$ passing through $p_0$. We assume that $J$ is smooth of class $\calc^{1,\alpha}$, 
$W$ is smooth of class $\calc^{2,\alpha}$, both for some $0<\alpha <1$, and that $\dim_{\rr}W = 
\frac{1}{2}\dim_{\rr}X$. Denote by $\Delta^+\deff \{\zeta\in \Delta:\im \zeta \ge 0\}$ the upper 
half-disk 
and by $\db_0\Delta^+ \deff (-1,1)$ its edge. Let a $J$-holomorphic map $u: \Delta^+ \to X$ be given, 
assume it is continuous up to $\db_0\Delta^+$, $u(\db_0\Delta^+)\subset W$ and $u(0)=p_0$.  We shall call $C = 
u(\Delta^+)$ a $J$-complex half-disk attached to $W$. Under our assumptions $u$ is of class $\calc^{2,
\alpha}$ up to the edge. We shall prove in Lemma \ref{van-order} below that 
there exists $\mu\in\nn$ such that $u(\zeta) - u(0)$ vanishes to order $\mu$ at zero, \ie in local 
coordinates near $p_0$ one has $u(\zeta) - u(0) = \zeta^{\mu}v(\zeta)$ for $\zeta$ in a neighbourhood 
of zero in $\Delta^+$. Here $v\in L^{1,p}_{loc}$ and $\zeta v\in L^{2,p}_{loc}$ for every $p>2$ 
and $v(0)\not=0$.



\smallskip There exists (locally) a $\calc^{2,\alpha}$-diffeomorphism $\Psi$  of $(X,W,p_0)$ to 
$(\rr^{2n},\rr^n,0)$ such that $\Psi_*J=J\st$ - the standard structure of $\cc^n$. 
Such coordinate change we shall call a {\slsf rectification map}. Notice that the structure 
$\Psi_*J$ obtained this way is still of class $\calc^{k,\alpha}$. Therefore modulo a 
rectification map  we may assume in the 
sequel that  $X=\rr^{2n}$, $W=\rr^n$ and $J|_{\rr^n} = J\st$. 

\smallskip Let $u:(\Delta^+,\db_0\Delta^+)\to (\rr^{2n}, \rr^n)$ be a $J$-holomorphic map. Set
\begin{equation}
\eqqno(ext-v1)
\tilde u(\zeta) =
\begin{cases}
u(\zeta) \text{ if } \im \zeta\ge 0\\[5pt]
\overline{u(\bar \zeta)} \text{ if } \im \zeta < 0,
\end{cases}
\end{equation}
and call $\tilde u$ the extension of $u$ by reflection with respect to $W = \rr^n$. 
$\tilde u$ is defined on $\Delta$ and we shall explain later that it is as a 
$\calc^{1,\delta}$-regular map for all $0<\delta <1$. But attention, we do not 
claim that $\tilde u$ is holomorphic for some structure and that it is of the same 
regularity $\calc^{2,\alpha}$ as was $u$! Now we can give the following.

\smallskip Now let $u_1, u_2: (\Delta^+,\db_0\Delta^+)\to (X,W)$ be two $J$-holomorphic 
mappings. We shall consider the case $\dim_{\rr}X=4$, $u_k(0)=p_0$ for $k=1,2$, and 
$\mu =1$ for both of them. We say that $u_1$ is a reparametrization of $u_2$ \emph{on the 
boundary} if there exists a biholomorphism  $\psi :  \Delta_r \to \Delta_r$,
$\psi (\zeta) = \pm \zeta + O(\zeta^2)$, such that one has
\begin{equation}
\eqqno(repar)
u_1(\xi) = u_2(\psi(\xi)) \qquad\text{ for } \qquad \xi \in (-r,r).
\end{equation}

Relation \eqqref(repar) has the clear geometric meaning. If $\psi (\zeta) = \zeta +
O(\zeta^2)$ then the half-discs $C_k\deff u_k(\Delta_r)$ coincide. In the case when 
$\psi(\zeta) = -\zeta + O(\zeta^2)$ one has that $C\deff C_1\cup C_2$ is a $J$-complex 
disc. Fix some rectifying map $\Psi : (X,W,p_0)\to (\rr^4,\rr^2,0)$, and notice that
$\Psi\circ u_k$ are $\Psi_*J$-holomorphic.

\begin{defi}
\label{index1}
We define the boundary intersection index $\ind_{\sf p_0}^{\sf b} (u_1,u_2)$ of $u_1$ and $u_2$
at $p_0$ as the intersection index at $0$ of the extensions by reflection 
$\widetilde{\Psi\circ  u_1}$ and $\widetilde{\Psi\circ u_2}$ of their rectifications 
$\Psi\circ u_1$ and $\Psi\circ u_2$. 
\end{defi}


We shall prove that $\bindo (u_1,u_2)$ is correctly defined, meaning that it
doesn't depend on the rectification map $\Psi$ as above and, moreover, the extensions 
of $\Psi\circ u_1$ and $\Psi\circ u_2$ by reflection locally intersect by a point provided one is not
a reparametrization of another on the boundary, see Remark \ref{pnt-2} and Proposition 
\ref{point-int}.


\medskip Our first goal in this paper is to prove the following statement.

\begin{thm}
\label{b-pos-thm}
Let $J$ and $W$ be as above and let $u_k:(\Delta^+,\db_0\Delta^+,0)\to (X, W, p_0)$ be 
two $J$-holomorphic maps continuous up to $\db_0\Delta^+$ such that the order of vanishing 
of $u_k(\zeta)-u_k(0)$ is one for both of them, and  such that one is not a reparametrization
of another on the boundary. Then:

\smallskip\sli the boundary intersection index $\ind_{\sf p_0}^{\sf b}(u_1,u_2)$ is positive;

\smallskip \slii $\ind_{\sf p_0}^{\sf b}(u_1,u_2)=1$ if and only if $u_k(\Delta^+)$ 
intersect at $p_0$ transversely;

\smallskip\sliii there exist an arbitrarily small perturbations $u_k^{prt}
:(\Delta^+_r,\db_0\Delta^+_r,0) \to  (X, W, p_0)$ 

\quad of $u_k$, defined on a smaller half-disk $\Delta^+_r$, which 
intersect only at  points of $W$, 

\quad and all their intersections are transverse;

\smallskip\sliv for any $u_k^{prt}$ as in (\sliii one has 
\begin{equation}
\eqqno(prt-int-f)
\ind_{{\sf p_0}}^{\sf b} (u_1,u_2) = \texttt{\#}\{u_1^{prt}(\db_0\Delta_r^+)\cap 
u_2^{prt}(\db_0\Delta_r^+)\}.
\end{equation} 
\end{thm}
\smallskip Notice that since such $u_k$-s are of class $\calc^{2,\alpha}$ up to the edge 
the notion of transversality (or tangency) at their boundary point has sense.






\begin{nnrema} \rm 
{\bf a)} We use the notation $\ind_{\sf p_0}^{\sf b} (u_1,u_2)$, but can write also 
$\ind_{\sf p_0}^{\sf b}(C_1,C_2)$, where $C_k=u_k(\Delta^+)$, since our $u_k$-s 
are embeddings near zero.

\smallskip\noindent {\bf b)} To our best knowledge this result is new even in the case of an  
integrable $J$, unless $W$ is assumed to be real analytic. For real analytic $W\subset \cc^2$ 
positivity of intersections can be obtained after reflecting both of $u_k(\Delta^+)$ with respect 
to $W$. This result is due to Alexander, \cite{A}. If $W$ and $J$ are both real analytic, but 
$J$ is not necessarily integrable, positivity of intersections follows from the possibility
to extend $u_k$-s to a neighbourhood of $\db_0\Delta^+$ as $J$-holomorphic maps, this was 
proved in \cite{IS3}, and then from the positivity of intersections of $J$-complex disks at 
their interior points.  The latter statement was proved in \cite{MW}.
\end{nnrema}

\newprg[INT.touch]{Index of intersection of two tangent half-disks}

\smallskip As the first application consider two {\sl complex} half-disks $C_1$, $C_2$ in 
$\cc^2$, \ie $C_k$ is the image  of a smooth up to the edge embedding $u_k:\Delta^+\to \cc^2$
which is holomorphic in the interior of $\Delta^+$. Assume that $u_1(0)=u_2(0)=0$ and
$du_1(0)[e_1] = \pm du_2(0)[e_1]$, where $e_1=(1,0)$. The 
order of tangency of such half-discs is well defined, see subsection \ref{TOUCH.touch-r}. 
We prove the following statement. 

\begin{thm}
\label{touch-thm}
Assume that $C_1$, $C_2$ are {\sl analytic}, \ie $C_k$ is the image  of a holomorphic in 
the standard sense  embedding $u_k:\Delta^+\to \cc^2$ smooth up to $\db_0\Delta^+$. Suppose 
that the order of tangency of $C_1$ with $C_2$ at zero is finite and equal to $d$. 
Then:

\smallskip\sli the index of intersection of $C_1$ with $C_2$ at zero is equal to $d$, 
in particular it is positive;

\smallskip\slii is equal to $\adyn$ if and only if $C_1$ intersects $C_2$ 
at zero transversely.
\end{thm}

\smallskip The proof of this theorem goes as follows. First we define the index of 
intersection of two \emph{real} half-discs at their common boundary point having the 
finite order of tangency at this point. This index can be any number, positive or negative. 
This definition goes through a construction of a real surface $W$ such that our hapf-discs 
are attached to it. Then in the case when these discs are complex we observe that 
$W$, as it is constructed, is totally real. Therefore the index in question 
is nothing but the boundary intersection index $\ind_{\sf p_0}^{\sf b}(u_1,u_2)$ as it
was defined above. By Theorem \ref{b-pos-thm} it occurs to be positive.

\newprg[INT.adj-f]{Adjunction Formula}
As the second application we prove the following formula. Let $M$ be an irreducible compact $J$-complex curve in an almost complex $4$-manifold $X$ with boundary on a totally real surface $W$.
Denote by  $\delta^{(b)}$ the sum of the boundary self-intersection indices of $M$ and 
by $\delta^{(i)}$) the sum  inner self-intersection indices of $M$. Further denote by
$\varkappa^{(i)}$ is the sum of inner cusps indices of $M$.

\begin{thm} {\slsf (Adjunction Formula).} 
\label{adj-1} 
Assume that $W$ is of class $\calc^{2,\alpha}$ and $J$ of class $\calc^{1,\alpha}$ for 
some $0<\alpha <1$, and assume that $M$ has no cusps on $W$.  Then 
\begin{equation}
\eqqno(adj-frm1)
2g + \sigma = \frac{[M^d]^2- \mu_M (TX,TW)}{2} + 2 - \delta^{(b)}
-2\delta^{(i)} - 2\varkappa^{(i)}.
\end{equation}
Here $g$ is the genus of $M$, $\sigma$ is the number of boundary circles, $[M^d]$ 
is the homological self-intersection of the Schottky double of $M$ and 
$\mu_M (TX,TW)$ is the Maslov index of the pair $(TX,TW)$ over $M$.
\end{thm}

The main new moment here (apart from the positivity of $\delta^{(i)}$) is the notion 
of the homological self-intersection of the Schottky double of $M$. We define it 
via the construction of a double of $(X,W)$. 

\newprg[INT.cusp]{Smoothing of boundary cusps}

We cannot say much about boundary intersections of complex curves with cusps on the 
boundary. But we prove the following partial statement in this direction. 

\smallskip
\begin{thm}
\label{b-pert-thm}
Let $(X,J,W)$ be as above and let a $J$-holomorphic $u:(\Delta^+,\db_0\Delta^+,0)\to 
(X, W,p_0)$ map be given. Then there exists an arbitrarily small perturbation 
$u^{prt}$ of $u$   defined on a smaller half-disc $\Delta_r^+ = \{\zeta \in \Delta^+:
|\zeta|<r\}$, such that $u^{prt}$ has no cusps.
\end{thm}

Here $u^{prt}: (\Delta_r^+, \db_0\Delta^+_r,$ $0) \to (X, W, p_0)$, \ie stays to 
be attached to $W$, to be more precise.

\newprg[INT.struc]{The structure of the paper.} Let us list the main ingredients 
of the proofs.

\smallskip\noindent{\slsf 1.} First we recall the reflection principle from \cite{IS2},
which makes possible to transfer  certain boundary value problems for $J$-complex 
disks to the inner ones. Using this tool we prove that in our settings that a 
non-constant $J$-holomorphic map satisfying a totally real boundary condition 
vanishes up to the finite order. This makes possible to treat the simple item 
(\slii of Theorem \ref{b-pos-thm}.

\smallskip\noindent{\slsf 2.} To treat the item (\sli  we prove the Comparison Lemma 
\ref{comp-lem-r}. 

\smallskip\noindent{\slsf 3.} In order to prove item (\sliii of Theorem \ref{b-pos-thm}
as well as Theorem \ref{b-pert-thm} we need perform certain perturbation of $J$-complex 
curves keeping them to be attached to the totally real submanifold. This should be
dome with estimates in $L^{1,p}$-norm. This is done in Lemma \ref{disk-pert-r} .

\smallskip\noindent{\slsf 4.} Theorem \ref{touch-thm} is proved afterwards and should be 
considered as an illustration of the results and methods of the first part of this paper. 

\smallskip\noindent{\slsf 5.} As for the Adjunction Formula the main point is to define 
the homological invariant which we call the self-intersection of the Schottky double
of a complex curve with boundary on totally real submanifold. This is done by 
doubling the almost complex manifold itself.

\smallskip\noindent{\slsf 6.} To prove Theorem \ref{b-pert-thm}, \ie to remove cusps, 
we perturbations with $L^{2,p}$-estimates. 

\smallskip\noindent{\slsf Acknowledgment.}  It should be underlined that the present 
exposition heavily depends on our
previous works  \cite{IS1,IS2,IS3,IS4}. We take care to recall the needed statements 
and explain how they should be modified to serve our ``boundary case'' here. I would 
like to use this opportunity to thank my co-authors V. Shevchishin and A. Sukhov, especially 
Vsevolod, for the fruitful collaboration along many years as well as for useful discussions 
on the subject of the present paper.

\newsect[REFL]{Extension by Reflection and the Order of Vanishing}

\newprg[REFL.redr]{Rectification a totally real submanifold}

Let $(X,J)$ be an almost complex manifold, $W$ a real submanifold of $X$.
$W$ is called {\slsf $J$-totally real} if $T_wW\cap J(T_wW) = \{0\}$ for every point 
$w\in W$. We shall consider only the case when $\dim_{\rr}W = \dim_{\cc}X$ and in 
this case the previous condition is equivalent to $T_wW\oplus J(T_wW) = T_{w}X$. 
We shall also always assume that $W$ is closed as a subset of $X$.

\begin{prop}
\label{redres-1}
Let $W$ be a totally real submanifold in an almost complex manifold $(X,J)$ of real 
dimension $n=\dim_{\cc }X$ of class $\calc^{k+1,\alpha}$, where the structure $J$ is 
supposed to be  of class $\calc^{k,\alpha}$, $k\ge 0, 0<\alpha <1$. Then for any point 
$x_0\in W$ there exists a neighbourhood 
$x_0\in W_0\subset W$, an open subset $X_0$ of $X$ with $X_0\cap W = W_0$ and a 
$\calc^{k+1,\alpha}$ - diffeomorphism $\Psi :(X_0,W_0)\to (\rr^{2n}, \rr^n)$ such that 
$\Psi_*J\mid_{\rr^n} = J\st $. 
\end{prop}
\proof  An appropriate $\calc^{k+1,\alpha}$-diffeomorphism of an appropriate
open $X_0$ with $X_0\cap W=:W_0$ will map the pair $(W_0,X_0)$ to $(\rr^n,\rr^{2n})$. 
Therefore, without loss of generality, we can suppose that $W_0$ is a relatively compact 
open subset of $\rr^n$, the same for $X_0$ and a $\calc^{k,\alpha}$-continuous almost 
complex structure $J$ is defined in the neighbourhood of the closure of $X_0$. Finally 
that $W_0\subset \rr^n$ is $J$-totally real. Coordinates in $\rr^{2n}$ we denote by
$x_1,...,x_n, y_1,..., y_n$. Let $\{ \frac{\db}{\db x_1},....,
\frac{\db}{\db x_n}\} $ be the standard basis of $T\rr^n$. Set
$v_j(x):=J(x) \frac{\db}{\db x_j}$. 

\smallskip\noindent {\slsf Step 1.} {\it We need to find $\calc^{k+1,\alpha}$-
functions $\phi_1,...,\phi_{2n}$ in the neighbourhood of $W_0$ such that}

\begin{equation}
\eqqno(alpha3)
\begin{cases}
\phi_j(x,0)=x_j \text{ for } j=1,...,n; \\[3pt]

\phi_j(x,0)=0 \text{ for } j=n+1,...,2n; \\[3pt]

\frac{\db\phi_j}{\db y_i}(x,0) = v_i^j(x) \text{ for } j=1,...,2n;

\end{cases}
\end{equation}

\smallskip\noindent where $v_i(x)= \sum_{j=1}^n v_i^j(x)\frac{\db}{\db x_j}
+ \sum_{j=n+1}^{2n} v_i^j\frac{\db}{\db y_{j-n}}$. Indeed, these relations for the mapping 
$\phi = (\phi_1,...,\phi_{2n})$ mean that $\phi$ preserves $\rr^n$ and 
\begin{equation}
\eqqno(alpha2)
\begin{cases}
\frac{\db\phi (x,0)}{\db y_i} = v_i(x)\\[3pt]
\frac{\db\phi (x,0)}{\db x_i} = e_i
\end{cases}
\text{ and that in its turn means that \ }
\begin{cases}
d\phi (x,0)\left[e_{n+i}\right] = v_i \\[3pt]
d\phi (x,0)\left[e_i\right] =e_i.
\end{cases}
\end{equation}
This implies that $(\phi^{-1}_*J|_{\rr^n})[e_i] = \big(d\phi^{-1}(x,0)
\circ J(x,0)\circ d\phi (x,0)\big)[e_i] = e_{n+i}$, 
\ie $\phi^{-1}_*J|_{\rr^n}$ is standard. Our problem is local. Therefore we can assume
that $W_0$ is the cube $Q$
\[
Q := \{ x=(x_1,\ldots,x_n)\in \rr^n :-\pi \leq x_j \leq \pi  \}
\]
and all involved functions have compact support in $Q$. We extend such
functions from $Q$ to $\rr^n$ being  $2\pi$-periodic
in each variable $x_j$. This means that those extensions
are functions on the torus $T^n := (\rr / 2\pi \zz)^n$.

\smallskip\noindent{\slsf Step 2.} {\it Let $v(x)$ be an $\calc^{k,\alpha}$-continuous 
function and $w(x)$ a $\calc^{k+1,\alpha}$-continuous function on the torus $T^n$. 
Then there exists a $\calc^{k+1,\alpha}$-continuous function $\phi(x,y)$ on
$T^n\times [-1,1]$ such that $\phi (x,0) = w(x)$ and $\frac{\db \phi}{\db y}(x,0)
= v(x)$.} Existence of such $\phi$ is the subject of the Trace theorem for H\"older 
classes, see \cite{Tr}.

\smallskip\noindent
{\slsf Step 3.} {\it There exist $\calc^{k+1,\alpha}$-continuous functions 
$\phi_1,..., \phi_{2n}$ on the unit ball in $\rr^{n+1}$ such that 
$\frac{\db \phi^{(1)}}{\db y_1}(x,0)= v_1(x)$ for $\phi^{(1)} := 
(\phi_1,...,\phi_{2n})$}. Indeed,  write $v_1=(v_1^1,...,v_1^{2n})$ and
apply Step 2 to each coordinate $v_1^j$ of $v_1$. More precisely for each 
$v_1^j$ find a $\calc^{k+1,\alpha}$-continuous function $\phi_j$ in the unit 
ball $W_0$ of $\rr^{n+1}$ such that $\phi_j(x,0)=x_j$ and $\phi_{n+j}(x,0)=0$ 
for $j=1,...,n$, and such that $\frac{\db \phi}{\db y_1}(x,0)=v_1^j$.

\smallskip\noindent
{\slsf Step 4.} {\it End of the proof.} Extend $v_2$ from $W_0$
to the unit ball $W_2$ in $\rr^{n+1}$ not depending on $y_1$. Apply
Step 3 to find $\phi^{(2)}(x,y_1,y_2)$ which is of class $\calc^{k+1,\alpha}$ 
in the unit ball $W_3$ of $\rr^{n+2}$ and satisfies the initial conditions
$\phi^{(2)}(x,y_1,0)=\phi^{(1)}(x,y_1)$ and $\frac{\db \phi^{(2)}}{\db
y_2}(x,y_1,0)=v_2(x,y_1)=v_2(x)$. In particular $\frac{\db \phi^{(2)}}{
\db y_2}(x,0,0) =v_2(x)$. Repeat this $n-2$ times to get $\phi^{(n)}$
which is $\calc^{k+1,\alpha}$ in the unit ball $W_{n+1}$ of $\rr^{2n}$ such
that $\phi^{(n)}(x,0,...,0)=...=\phi^{(1)}(x,0)=x$ and $\frac{\db
\phi^{(n)}}{\db y_j}(x;0) =v_j(x)$. Therefore $\phi^{(n)}$ is the
needed solution  of \eqqref(alpha3).

\smallskip\qed

\begin{rema}\rm
\label{red-rem}
We shall need later the following observation: after a coordinate
change sending $W$ to $\rr^n$ the further ``correction'' of the structure $J$ is 
performed by a diffeomorphism which is identity on $\rr^n$.

\end{rema}

\newprg[REFL.refl]{Reflection Principle}

\smallskip 
It turns out that {\sl one can reflect} a $J$-holomorphic mapping with
respect to $W$, preserving its holomorphicity (!), but not in $X$. In order to do 
this one should ``change the range'' of the mapping and ``reflect'' not only $u$ but 
also the structure $J$. This trick was introduced in \cite{IS2}, let us briefly recall
it here. 

\smallskip We say that a continuous map $u:\Delta^+ \to X$ satisfies the {\slsf 
totally real boundary condition} $W$ along $\db_0\Delta^+$ if $u|_{\db_0\Delta^+}$ 
takes its values in $W$. In the sequel for a mapping $u:(\Delta^+,\db_0\Delta^+)\to 
(X,W)$ to be $J$-holomorphic means that $u$ is of the class $(\calc^0\cap L^{1,2})
(\Delta^+)$, \ie  the norm $\norm{du}_{L^2(\Delta^+)}$ is bounded and that
\begin{equation}
\eqqno(cr1)
\bar\db_Ju:=du+J(u)\circ du\circ J\st =0 
\end{equation}
a.e. in  $\Delta^+$. Under these assumptions $u$ is of class $\calc^{k+1,\alpha}$ up
to $\db_0 \Delta^+$ provided $J\in \calc^{k,\alpha}$ and $W\in\calc^{k+1,\alpha}$, 
see Theorem 1.3 in \cite{IS3}. 

\smallskip Let such $J$-holomorphic $u$ be given. We assume, if the opposite is not
explicitly mentioned, that $(X,W) = (\rr^{2n}, \rr^n)$ and $J|_{\rr^n}=J\st$. Set 
$E^+\deff u^*TX$ and $F\deff u^*TW$. Then $E^+$ is a trivial vector bundle over 
$\Delta^+$, \ie  $E^+=\Delta^+ \times \rr^{2n}$. Likewise $F$ is a trivial vector 
bundle over $\db_0\Delta^+$, \ie  $F=\db_0\Delta^+ \times \rr^{n}$. Triviality of 
bundles is understood here as triviality of real bundles of regularity $\calc^{k+1,\alpha}$. Bundle $E^+$ is complex in the sense that it possesses a natural complex linear structure $J_u$, here $J_u(\zeta)\deff J(u(\zeta))$ acts on the fibre $E^+_\zeta \deff \{\zeta\}\times \rr^{2n}$. This structure $J_u$ is of class $\calc^{k,\alpha}$ if $J$ was such because $u$ is of class
$\calc^{k+1,\alpha}$ up to $\db_0\Delta^+$ as explained above. Notice that $F$ is 
the real part of $E^+|_{\db_0\Delta^+}$, \ie $E^+|_{\db_0\Delta^+}= F\oplus J_uF$.
Moreover, due to our identifications we have that $F= \db_0 \Delta^+\times \rr^n$ 
and the fact that $J|_{\rr^n}=J\st$ implies that $J_u(\zeta) = J\st$ for $\zeta 
\in \db_0\Delta^+$. Our $J$-holomorphic map $u:\Delta^+\to \rr^{2n}$ can be regarded
as a section of $E^+$ over $\Delta^+$, and this section (which we denote still as $u$) 
satisfies the Cauchy-Riemann equation 
\begin{equation}
\eqqno(cr2)
\dbar_{J_u}u\deff du + J_u\circ du\circ J\st = 0,
\end{equation}
\ie is $J_u$-holomorphic. Equation \eqqref(cr2) is just the rewritten \eqqref(cr1). 
In addition $u|_{\db_0\Delta^+}$ is a section of $\db_0\Delta^+\times \rr^n$. We 
shall use the same notation $\dbar_Ju$ for the operator $du + J_u\circ du\circ J\st$ 
as in \eqqref(cr1).

\begin{rema} \rm
{\bf a)} And one more observation: given some section $v$ of $E^+$. It can be 
$J_u$-holomorphic or not. It can be naturally considered as a mapping to 
$\rr^{2n}$, we denote this mapping as $v$ as well. And vice versa a mapping 
$v:\Delta^+\to \rr^{2n}$ we can consider as a section of $E^+$ denoting it
with the same letter. I.e., we do not distinguish between sections of $E^+$
and mappings to $\rr^{2n}$ hoping that this will not lead to a confusion. 
And in the same spirit, interpreting mapping $v:\Delta^+\to \cc$ as a section of $E^+$,
after proving some properties of its extension by reflection $\tilde v$ as 
a section of $E$ (this provided $v$ is $\rr^n$-valued on $\db_0\Delta^+$), 
we shall interpret backwards this properties as properties of $\tilde v$ as a
mapping from $\Delta$ to $\cc^n$.

\smallskip\noindent{\bf b)} If $E^+$ is a complex vector bundle over $\Delta^+$ 
and $F$ a real subbundle of the restriction $E^+|_{\db_0\Delta^+}$ we denote by  
$L^{1,p}(\Delta^+, E^+, F)$ and respectively by $\calc^{k,\alpha}(\Delta^+, E^+, F)$ 
the spaces of sections of $E^+$ over $\Delta^+$ of the corresponding smoothness 
whose restrictions to $\db_0\Delta^+$ take values in $F$.  
\end{rema}

Denote by $\tau:\cc^n \to \cc^n$ the standard complex conjugation in $\cc^n$ and
by $\tau_{\Delta}$ the standard complex conjugation in the unit disk, set $\Delta^-
\deff \tau_{\Delta}(\Delta^+)$. Let $E\deff \Delta\times \cc^n$ be the trivial 
vector bundle over the disk. Extend $J_u$ to $E$ as follows 
\begin{equation}
\eqqno(ext-ref-j)
\tilde J_u(\zeta)[\vect ] \deff - \tau\left(J_u(\tau_{\Delta}(\zeta))[\tau (\vect)]\right)
= - \overline{J_u(\bar \zeta)[\bar\vect]} \quad\text{ for }  \zeta\in \Delta^- 
\text{ and } \vect \in E_\zeta\deff\{\zeta\}\times \cc^n.
\end{equation}

Notice that the extended $\tilde J_u$ is a complex structure on $E$ and the 
following natural global involution $\tau_E$ of the total space of $E$
\[
\tau_E (\zeta,v)\deff (\tau_{\Delta }(\zeta),\tau (v)).
\]
is $\tilde J_u$-antiholomorphic, \ie 
\[
\tau_E\circ\tilde J_u = -\tilde J_u\circ \tau_E.
\]
To check this write 
\[
\tilde J_u^2(\zeta)[\vect] = \tilde J_u(\zeta)\left(-\overline{J_u(\bar \zeta)[\bar\vect]}\right)
= \overline{J_u^2(\bar \zeta)[\bar\vect]} = -\vect ,
\]
and also
\[
(\ti J_u \scirc \tau_E)(\zeta)[\vect] = \ti J_u(\tau_\Delta \zeta)[\tau \vect] =
-\tau J_{u(\tau^2_\Delta \zeta)} [\tau^2 \vect] = - \tau_E J_{u(\zeta)} [\vect] =
-(\tau_E \scirc J_u)(\zeta)[\vect ].
\]

\begin{rema} \rm 
{\bf a)} Notice that the structure $\tilde J_u$ is  Lipschitz-continuous only 
whatever was the assumed smoothness of $J$! Let us give an example. Consider 
the following complex structure $J(\zeta) = J(\xi + i\eta)$ on $\Delta^+\times
\rr^4$:
\begin{equation}
J(\xi + i\eta) = 
\left(
\begin{matrix} 
 0 & -1 & 0 & 0 \cr
 1 & 0 & 0 & 0 \cr
 0 &\eta & 0 &-1\cr
 \eta & 0 & 1 & 0
\end{matrix}
\right) \quad\text{ for }\quad \eta \ge 0, \text{ \ie } J \text{ is real analytic }. 
\end{equation}

It is easy to see that extended by reflection structure has the form  
\[
\tilde J(\xi + i\eta) = 
 \left(
\begin{matrix}
 0 & -1 & 0 & 0 \cr
 1 & 0 & 0 & 0 \cr
 0 &-\eta & 0 &-1\cr
 -\eta & 0 & 1 & 0
\end{matrix}
\right) \quad \text{ for } \quad \eta <0.
\]

Therefore the extended to the whole of $\Delta\times \rr^4$ structure has the form 
\[
\tilde J (\xi + i\eta ) = 
\left(
\begin{matrix}
 0 & -1 & 0 & 0 \cr
 1 & 0 & 0 & 0 \cr
 0 & |\eta| & 0 &-1\cr
 |\eta | & 0 & 1 & 0
\end{matrix}
\right) \quad \text{ for } \quad \zeta = \xi + i\eta\in \Delta . 
\]

I.e., a smooth structure extends by reflection only to a Lipschitz-continuous one.

\medskip\noindent{\bf b)} At the same time it will be important in the sequel to
notice that if original $J$ was of class $\calc^{k,\alpha}, k\ge 1$, then 
$\tilde J_u$ is not only Lipschitz-continuous but also "piecewise $\calc^{k,\alpha}$"
in the sense that it is $\calc^{k,\alpha}$-regular on $\Delta^{\pm}$ {\slsf up to}
$\beta_0=(-1,1)$. 
\end{rema}

\smallskip For $p>1$ we define the following continuous extension operator 
\[
\ext : L^{1,p} (\Delta^+,E^+,F)\to  L^{1,p}(\Delta, E)
\]
\begin{equation}
\eqqno(ext-v2)
\ext (v)(\zeta) = \tau v (\tau_{\Delta }\zeta) =\overline{v(\bar\zeta)}.
\end{equation}

\smallskip
Set $\tilde v = \ext (v)$ for short. In \cite{IS2} it was proved that for a 
$J_u$-holomorphic section $v$ of $E^+$ with boundary condition $v|_{\db_0\Delta^+}
\subset F$ its extension $\tilde v$ is a $\tilde J_u$-holomorphic section of $E$. 
In particular, since the section $u$ is $J_u$-holomorphic and of class 
$\calc^{k+1,\alpha}$ up to $\db_0\Delta^+$ its extension 
by reflection $\tilde u$ will be $\tilde J_u$-holomorphic. Moreover, since the
structure $\tilde J_u$ is Lipschitz-continuous, \ie belongs to $L^{1,\infty}_{loc}
\subset \bigcap_{p>2}L^{1,p}_{loc}$ we have that $\tilde u$ is of class 
$L^{2,p}_{loc}$ for all $p>2$. Finally by Sobolev  embedding  
$L^{2,p}_{loc}\subset \calc^{1,\delta}$ with $\delta = 1-\frac{2}{p}$ our
$\tilde u$ is of class $\calc^{1,\delta}$ for all $0<\delta<1$. 
In  the sequel we use the abbreviation ``$L^{2,p}$-regular'' instead 
of ``$L^{2,p}$-regular for any $2<p<\infty$'' and a similar abbreviation for 
``$\calc^{1,\delta}$-regular'', \ie $\calc^{1,\delta}$ for all $0<\delta <1$.

\begin{rema} \rm
Notice that via identification of sections of $E^+$ (resp. of $E$) and mappings from 
$\Delta^+$ to $\cc^n$ (resp. from $\Delta$ to $\cc^n$) our extension operator is 
exactly the extension by reflection as in \eqqref(ext-v1). 
\end{rema}

\newprg[REFL.dbar]{Uniqueness theorems for $\dbar$-inequalities}
We shall need a generalization of Lemma 1.4.1 from \cite{IS1} as well as one 
corollary from it. In that lemma we considered a $\dbar$-inequality for a 
$\cc^n$-valued $L^{1,2}$-function $u$ in the unit disk of the form 
\begin{equation}
\eqqno(dbar-in1)
|\dbar u| \leq h\cdot |u|, 
\end{equation}
where $h\in L^p_{loc}(\Delta)$ for some $p>2$. The statement is that a solution $u$ 
of \eqqref(dbar-in1) obeys the unique continuation theorem, the same as for holomorphic 
functions. Let us state this more  explicitly .
\begin{lem}
\label{thm3.1.1}
Suppose that a function $u\in L^2_\loc(\Delta, \cc^n)$ with $\dbar 
u\in L^1_\loc(\Delta,\cc^n)$ satisfies {\sl a.e.} the inequality \eqqref(dbar-in1) 
for some non-negative $h\in L^p_\loc(\Delta)$ with $p>2$. Then:

\smallskip
\sli $u\in L^{1, p}_\loc(\Delta)$, in particular $u\in
\calc^{0,\alpha}_\loc (\Delta)$ with $\alpha\deff1-\frac{2}{p}$;

\slii if $u$ is not identically zero then for any $\zeta_0\in\Delta$ such that
$u(\zeta_0)=0$ there exists $\mu\in\nn$,  

- the~multiplicity of zero of $u$ at $\zeta_0$ - such that for $\zeta$ in a
neighbourhood of $\zeta_0$ one has
\begin{equation}
\eqqno(dbar-in2)
u(\zeta) = (\zeta - \zeta_0)^{\mu}v(\zeta),
\end{equation}
\quad for some $v\in L^{1,p}_{loc}(\Delta)$ with $v(\zeta_0) \not= 0$.
\end{lem}

Now let us state a generalization of this lemma. Let $J$ be a continuous 
complex structure in the trivial $\rr^{2n}$-bundle over the unit disk
$\Delta$. By this we mean that $J=\{J_{\zeta}\}_{\zeta\in \Delta}$ is 
a continuous family of linear operators $J_{\zeta}:\rr^{2n} \to \rr^{2n}$ 
such that $J_{\zeta}^2 = -\id$ for every $\zeta$. With $J$ we can associate
a $\dbar$-operator $\dbar-J:L^{1,p}(\Delta , E)\to L^p(\Delta ,E)$ as 
follows
\begin{equation}
\eqqno(dbar_j1)
\dbar_Ju(\zeta) = \db_{\xi} u(\zeta) + J_{\zeta}\db_{\eta} u(\zeta).
\end{equation}

\begin{lem} 
\label{dbar-in1}
Let $J$ be an $L^{1,p}$-regular for some $p>2$ complex structure in the trivial 
$\rr^{2n}$-bundle over the disk $\Delta$ such that $J(0) =J\st$. 
Suppose that a function $u\in L^{1,2}_\loc(\Delta, \cc^n)$ is not identically zero
and satisfies for some non-negative $h\in L^p_\loc(\Delta)$ {\sl a.e.} the inequality
\begin{equation}
\eqqno(3.1)
|\dbar_J u| \leq h\cdot | u | 
\end{equation}
Then:

\smallskip
\sli $u\in L^{1, p}_\loc(\Delta)$, in particular $u\in
\calc^{\alpha}_\loc (\Delta)$ with $\alpha\deff 1-\frac{2}{p}$;

\slii for any $\zeta_0\in\Delta$ such that $u(\zeta_0)=0$ there exists
$\mu\in\nn$ - the~multiplicity of zero of $u$ 

in $\zeta_0$ - such that $u(\zeta)=(\zeta - \zeta_0)^\mu \cdot v(z)$ for some 
$v\in L^{1, p}_\loc (\Delta)$ with $v(\zeta_0)\not=0$.
\end{lem}
\proof Fix a $(J\st, J)$-complex bundle isomorphism $\Phi: \Delta\times\cc^n \to
\Delta\times\cc^n$ of regularity $L^{1,p}$, \ie $\Phi$ is such that $\Phi\inv \circ 
J \circ \Phi = J\st$. Then any section $u(\zeta)$ of $\Delta\times\cc^n$ has the 
form $u(\zeta) = \Phi (w(\zeta))$ and $u(\zeta)$ is $L^{1,p}$-regular if and only 
if such is $w(\zeta)$. Moreover,
\[
\dbar_J u(\zeta) = (\partial_{\xi} + J(\zeta)\partial_{\eta})\Phi (w(\zeta)) =
\Phi\big( \partial_{\xi} + \Phi\inv \cdot J(\zeta) \cdot \Phi \partial_{\eta}
\big)\,w(\zeta) +
\]
\[
 \;+\; (\partial_{\xi} \Phi + J(\zeta)\partial_{\eta}\Phi)\,w(\zeta)
= \Phi \dbar_{J\st}w + (\partial_{\xi} \Phi + J(\zeta)\partial_{\eta}\Phi)\,w(\zeta).
\]
Consequently, \eqqref(3.1) is equivalent to the differential inequality
\[
| \dbar\st w | \leq | \Phi\inv( \dbar_J u(\zeta) )| 
+  |\Phi\inv (\partial_{\xi} \Phi+ J(\zeta)\partial_{\eta}\Phi)\,w(\zeta) | \leq
\]
\begin{equation}
\leq  h\cdot | \Phi\inv ||u| +  |\Phi\inv (\partial_{\xi}\Phi + J(\zeta)\partial_{\eta}
\Phi)|\,w(\zeta)| \leq h_1\cdot |w| \eqqno(3.2)
\end{equation}
with a new $h_1\in L^p(\Delta)$. The statement of the lemma  is reduced now 
to Lemma \ref{thm3.1.1}.

\smallskip\qed

The following corollary from Lemma \ref{thm3.1.1} was given in \cite{IS1}, but 
without proof. Since it will be used here we shall state and prove it below. First, recall the following theorem of Harvey-Polking, see \cite{HP}.
\begin{thm}
\label{thm-Ha-Po}
Let $f:\Delta \to \cc^n$ be locally $L^2$-integrable.  Assume that for some $g \in L^1_\loc
(\Delta, \cc^n)$ the equation $\dbar f =g$ holds (in the weak sense)
in the punctured disc $\check\Delta$. Then $\dbar f =g$ holds in the
whole disc $\Delta$.
\end{thm}

Now we have the following

\begin{corol}
\label{thm3.1.2}
Under the~hypothesis of Lemma \ref{thm3.1.1} suppose additionally that
$u$ satisfies {\sl a.e.} the inequality
\begin{equation}
\eqqno(3.1.2)
|\dbar u(\zeta)| \leq |\zeta - \zeta_0 |^\nu h(\zeta)\cdot |u(\zeta)|,
\end{equation}
with some $\zeta_0\in \Delta$, $\nu\in\nn$, and $h\in L^p_\loc(\Delta)$ for some 
$2<p<\infty$. Then
\begin{equation}
\eqqno(3.1.3)
u(\zeta) = (\zeta - \zeta_0)^\mu\left[ \bigl(P^{(\nu)}(\zeta) + (\zeta - \zeta_0)^\nu
v(\zeta)\bigr)\right], 
\end{equation}
where $\mu\in\nn$ is the multiplicity of zero of $u$ at $\zeta_0$,
defined as above in Lemma \ref{thm3.1.1}, $P^{(\nu)}(\zeta)$ is a holomorphic 
polynomial in $\zeta$ of degree 
$\leq \nu$ with $P^{(\nu)}(\zeta_0) \not= 0$, and $v\in L^{1,p}_{loc}$ with  
$v(\zeta_0)=0$. In  particular $v\in \calc^{0,\alpha}$, $\alpha=1-\frac{2}{p}$,
and $v(\zeta)=O(|\zeta - \zeta_0|^\alpha)$.
\end{corol}
\proof Set $u_0(\zeta)= \frac{u(\zeta)}{(\zeta-\zeta_0)^\mu}$ and $h_1(\zeta) \deff 
h(\zeta)\cdot | u_0(\zeta) |$. By Lemma \ref{thm3.1.1}, $u_0\in \calc^{0,\alpha}$,
$u_0(\zeta_0)\not=0$, $h_1\in L^p_\loc$, and $u_0$ satisfies {\sl a.e.}
the~inequality 
\begin{equation}
\eqqno(dbar-u0)
|\dbar u_0(\zeta)| \leq | \zeta-\zeta_0 |^\nu h_1(\zeta).
\end{equation}
Indeed, for $\zeta\not=\zeta_0$ one has $|\dbar u_0(\zeta)|=\left|\frac{\dbar u(\zeta)}
{(\zeta-\zeta_0)^{\mu}}\right|$$\leq \frac{|\zeta-\zeta_0|^{\nu}h(\zeta)|u(\zeta)|}
{|\zeta-\zeta_0|^{\mu}}=|\zeta-\zeta_0|^{\nu}h(\zeta)|u_0(\zeta)|=
|\zeta-\zeta_0|^{\nu}h_1(\zeta)$, and then the claim follows from the Harvey-Polking 
theorem. Set $a_0=u_0(\zeta_0) = \lim_{\zeta\to \zeta_0}
\frac{u(\zeta)}{(\zeta-\zeta_0)^{\mu}}$. Since $u_0(\zeta)-a_0 = O(| \zeta-\zeta_0 |^\alpha)$ 
we have $u_1\deff\msmall{\frac{u_0(\zeta) - a_0}{\zeta-\zeta_0}}\in L^2_\loc$. 
Observe that again for $\zeta\not=\zeta_0$ we have $|\dbar u_1(\zeta)| = \left|\frac{\dbar 
u_0(\zeta)}{\zeta - \zeta_0}\right|\le |\zeta - \zeta_0|h_1(\zeta)$.
Applying the~theorem of Harvey-Polking once more, we obtain  $|\dbar u_1 | \leq 
|\zeta-\zeta_0|^{\nu-1} h_1$, and consequently $u_1\in\calc^{0,\alpha}$, $u_1(\zeta) - 
u_1(\zeta_0)=O(| \zeta-\zeta_0 |^\alpha)$. Repeating this procedure $\nu$ times, we obtain 
the~polynomial
\[
P^{(\nu)}(\zeta)=a_0 + (\zeta-\zeta_0)a_1 + \cdots + (\zeta-\zeta_0)^\nu a_\nu
\]
with
\begin{equation}
\eqqno(a-k)
a_k\deff\lim_{\zeta\to \zeta_0} \frac{u(\zeta)-\sum_{i=0}^{k-1}(\zeta-\zeta_0)^i a_i}
{(\zeta-\zeta_0)^k}, \qquad 0\leq k\leq\nu,
\end{equation}
and the~function
\[
v(\zeta)\deff \msmall{\frac{u(\zeta)-P^{(\nu)}(\zeta)}{(\zeta-\zeta_0)^\nu }},
\]
which satisfies the~conclusion of the Corollary.

\smallskip\qed

\newprg[REFL.norm]{The normal form of a holomorphic map}

Let $J$ be a Lipschitz-continuous almost-complex structure on $\rr^{2n}$ such that
$J|_{\rr^n}=J\st$ and let $u:\Delta^+\to \rr^{2n}$ be a $J$-holomorphic map
such that $u|_{\db_0\Delta^+}$ takes its values in $\rr^n$.

\begin{lem}
\label{van-order}
Suppose that $\norm{J-J\st}_{\calc^{Lip}}$ and $\norm{u}_{L^{1,p}(\Delta^+)}$ are small enough. 
Then there exist $\mu\in\nn$ and a holomorphic polynomial $P(\zeta)$
of degree $\leq\mu-1$ with real coefficients and $P(0)\not=0$ such that
\begin{equation}
\eqqno(nrm-frm1)
u(\zeta)=\zeta^\mu\cdot P(\zeta) +\zeta^{2\mu -1}v(\zeta) \quad\text{ for } \quad
\zeta\in \Delta^+,
\end{equation}
where the function $v(\zeta)$ is $\rr^n$-valued for real $\zeta$  and satisfies 
the estimates
\begin{equation}
\eqqno(est-v)
\norm{v}_{L^{1,p}(\Delta^+)}\le C\cdot \norm{u}_{L^{1,p}(\Delta^+)} \quad\text{ and }
\quad \norm{\zeta v}_{L^{2,p}(\Delta^+_r)}\le C_p\cdot r^{\alpha}\norm{u}_{L^{1,p}(\Delta^+)}
\end{equation}
for $0<r<1$, $p>2$ and $\alpha = 1-\frac{1}{2}$.
\end{lem}

\proof Extend $J_u\deff J\circ u$ by reflection to a complex structure $\tilde J_u$ 
on $E = \Delta\times \rr^{2n}$ as in \eqqref(ext-ref-j). Note that $\tilde J_u$ is 
Lipschitz-continuous. Consider $u$ as a $J_u$-holomorphic section of $E^+=\Delta^+\times
\rr^{2n}$ and then extend it to a $\tilde J_u$-holomorphic section $\tilde u$ of $E$ 
over $\Delta$ as in \eqqref(ext-v2). Notice that for $\im \zeta <0$ we have
\[
\left|\tilde J_u(\zeta)[\vect] -J\st [\vect]\right| = \left|-\overline{J_u(\bar\zeta)[\bar\vect]}
- J\st [\vect]\right| = \left|\left[J\st - J(u(\bar\zeta))\right][\bar\vect]\right|\le 
\norm{J}_{\calc^{Lip}}|u(\bar\zeta)|\,|\vect|.
\]
And since $|u(\bar\zeta)| = |\tilde u(\zeta)|$ we obtain for $\im \zeta <0$ the 
following estimate
\begin{equation}
\eqqno(est-tld-j)
|\tilde J_u - J\st|\le \norm{J}_{Lip}|\tilde u|.
\end{equation}
For $\im \zeta \ge 0$ the relation \eqqref(est-tld-j) is obvious. Using this remark and
the identity
$\dbar_{\tilde J_u}\tilde u = d\tilde u + \tilde J_u \scirc d\tilde u \scirc J\st 
= 0$ we obtain pointwisely almost everywhere 

\begin{equation}
|\dbar_{J\st} \tilde u |  = \bigl|d\tilde u + J\st \scirc d\tilde u \scirc
J\st\bigr| = \bigl| \msmall(J\st - \tilde J_u)\scirc d\tilde u \scirc J\st
\bigr| \leq \norm{ J }_{Lip} \cdot | \tilde u |  \cdot | d\tilde u |.
\eqqno(3.1.5)
\end{equation}
Since  $d\tilde u\in L^p_{loc}$ we conclude that our section $\tilde u$ satisfies 
the $\dbar$-inequality
\[
|\dbar_{J\st}\tilde u |\le h|\tilde u|,
\]
where $h\in L^p(\Delta)$. From Lemma \ref{thm3.1.1} 
we coclude that $\tilde u(\zeta) = \zeta^\mu \tilde v(\zeta)$ with some $\mu \in\nn$
and $\tilde v\in L^{1,p}_\loc(\Delta, \cc^n)$ for any $p<\infty$, $\tilde v(0)\not=0$ 
and we see that $\tilde v(\zeta)$ is $\rr^n$-valued for real $\zeta$. If $\mu =1$ we are done
by setting $v\deff \tilde v|_{\Delta^+}$. Otherwise $d\tilde u(\zeta)/\zeta^{\mu-1} = 
\mu\,\tilde v\,d\zeta + \zeta\,d\tilde v \in L^p_\loc(\Delta, \cc^n)$
for any $p<\infty$. Therefore Lemma follows now from the Corollary \ref{thm3.1.2}. Reality 
condition on coefficients of $P$ and $v|_{\db_0\Delta^+}$ follows from the equality 
$\overline{\tilde u(\bar \zeta)} = \tilde u(\zeta)$, the latter is the 
definition \eqqref(ext-v2) of how $u$ was extended. All what is left is to restrict
$\tilde u$ to $\Delta^+$ to obtain the first conclusion of the Lemma.

\smallskip Estimate \eqqref(est-v) were proved in Lemmas 2.1 and 2.3 of \cite{IS4} for 
$\tilde v$, $\tilde u$ in our present notations. One only remarks that 
$\norm{\tilde J_u-J\st}_{Lip (\Delta)} = \norm{J_u-J\st}_{Lip (\Delta^+)}$
and $\norm{\tilde u}_{L^{1,p}(\Delta)} = 2^{1/p}\norm{u}_{L^{1,p}(\Delta^+)}$, the same 
for $v$. 

\smallskip\qed

\begin{rema} \rm 
\label{r-cor3.1.3}
{\bf a)} Writing \eqqref(nrm-frm1) as
\begin{equation}
\eqqno(nrm-frm2)
u(\zeta) = \vect_0\zeta^{\mu} + O(|\zeta|^{\mu + \alpha}),
\end{equation}
\ie $\vect_0 = P(0)$, we shall call $\vect_0$ the {\slsf tangent vector} to $u$ at zero and 
\eqqref(nrm-frm2) the {\slsf normal form} of a $J$-complex half-disk attached to a totally real 
submanifold. The number $\mu$ we shall call the {\slsf order of vanishing} of $u$ at zero. Let us 
underline that \eqqref(nrm-frm2) comes from the same expression for the extension $\tilde u$
of $u$ by reflection, \ie from the fact that 
\[
\tilde u(\zeta) = \vect_0\zeta^{\mu} + O(|\zeta|^{\mu + \alpha}) \quad\text{ for }\quad 
\zeta\in \Delta .
\]
Let us stress that $\vect_0 \in \rr^n$ and that the term $O(|\zeta|^{\mu + \alpha})$ is 
$\rr^n$-valued for real $\zeta$. Indeed, writing $\zeta = \xi + i\eta$ we get from 
\[
 \tilde u(\xi) = \vect_0\xi^{\mu} + O(|\xi|^{\mu + \alpha})
\]
that $\vect_0 = \frac{\db u(\xi}{\db \xi}(0)$ is real and then this implies the same statement
for $O$.

\smallskip\noindent{\bf b)} It will be important for us that the differential of $\tilde u$ 
satisfies
\begin{equation}
\eqqno(nrm-frm3)
d\tilde u (\zeta)  = \mu \vect_0\zeta^{\mu - 1} + O(|\zeta|^{\mu-1+\alpha}),
\end{equation}
and this implies the same relation for $u=\tilde u|_{\Delta^+}$. To prove \eqqref(nrm-frm3)
write, using \eqqref(nrm-frm1) and \eqqref(est-v), the mapping $\tilde u$ in the form
\[
\tilde u(\zeta) = \zeta^{\mu} P(\zeta) + \zeta^{2\mu - 2}\left(\zeta \tilde v(\zeta)\right),
\]
and differentiate the rest. Taking into account that by \eqqref(est-v) we have that $d(\zeta 
\ti v) = O(|\zeta|^{\alpha})$ we obtain 
\[
d\left(\zeta^{2\mu -2}\zeta \tilde v(\zeta)\right) = \zeta^{2\mu -2}O(|\zeta{\alpha}) + 
\zeta^{2\mu -3}O(|\zeta|^{1+\alpha}) = \zeta^{2\mu -2}O(|\zeta|^{\alpha}).
\]
Now 
\[
d\tilde u(\zeta) = d\left(\zeta^{\mu}P(\zeta)\right) + d\left(\zeta^{2\mu -2}\zeta \tilde v(\zeta)
\right) = d\left(\zeta^{\mu}P(\zeta)\right) + O(|\zeta|^{2\mu -2 + \alpha}) = 
\]
\[
 = \mu \zeta^{\mu -1}(0) + O(|\zeta|^{\mu} + O(|\zeta|^{2\mu -2 + \alpha}) = \mu \vect_0\zeta^{\mu -1} 
 + O(|\zeta|^{\mu - 1 + \alpha}),
\]
because $2\mu - 2 \ge \mu -1$.

\smallskip\noindent{\bf c)} Finally this $\mu$ doesn't depend on the rectification map. 
Indeed, let $(B,W,J)$ and $u:(\Delta^+,\db_0\Delta^+)\to (B,W)$ are as above, where $u$
is $J$-holomorphic, then for a rectification diffeomorphism $\Psi^1:(B,W)\to (\rr^{2n},\rr^n)$ 
the $\Psi^1_*J$-holomorphic half-disk $u_1\deff\Psi^1\circ u$ has the form \eqqref(nrm-frm1) 
in the standard complex coordinates of $\cc^n$ with parameters $\mu_1$ and $P_1$. Given another 
rectification map $\Psi^2$ we repeat the same and get $u_2\deff\Psi^2\circ u$ in the form 
\eqqref(nrm-frm1) with parameters $\mu_2,P_2$. Composition $\Psi\deff \Psi^2\circ (\Psi^1)^{-1}$ 
is a $\calc^{2,\alpha}$-diffeomorphism, and sends $u_1$ to $u_2$. Therefore the orders of 
vanishing $\mu_1$ (resp. $\mu_2$) of $u_1$ (resp. $u_2$) at the origin are the same.

\smallskip\noindent{\bf d)} We shall use repeatedly the obvious observation (already 
used in the proof of Lemma \ref{van-order}) that $Lip_{\Delta}(\tilde J_u) = 
Lip_{\Delta^+}(J_u) \le Lip_B(J)\cdot Lip_{\Delta}(u)$ provided $u(\Delta^+)
\subset B$.
\end{rema}

\newsect[COMP]{Comparison of Tangent Analytic Half-Disks}

Our tool for the proof of Theorem \ref{b-pos-thm} is the following statement about the 
comparison of analytic disks. We shall state and prove it in  full generality.

\newprg[COMP.comp]{Comparison of touching half-discs}
This time we consider $J$-holomorphic mappings $u_k:(\Delta^+, $ $\db_0\Delta^+,0) 
\to (\rr^{2n}, \rr^n,0)$ that have the same  tangent vector at zero 
and the same order of vanishing,~\ie
\[
u_k(\zeta) = \vect_0 \zeta^{\mu} + O(|\zeta|^{\mu+\alpha}) \quad\text { for } \quad
k=1,2.
\]
In this case we shall say that $u_k$ \emph{touch} each other at zero.
Our task is to describe in the best possible way their difference. This can be done 
with the help of the following {\sl Comparison Lemma}.

\begin{lem} {\slsf (Comparison Lemma - I).}
\label{comp-lem-r}
Suppose that $J$ is Lipschitz-continuous. Then there exists a holomorphic function 
$\psi (\zeta) = \zeta+ O(\zeta^2)$, real-valued for real $\zeta$, 
an integer $\nu > \mu$ and an $L^{1,p}_{loc}$-regular $\cc^n$-valued function 
$w(\zeta)$, $\rr^n$-valued for real $\zeta$, such that for some $r>0$  one has
\begin{equation}
\eqqno(9.1.1)
u_2(\zeta) = u_1(\psi(\zeta)) + \zeta^{\nu} w(\zeta). 
\end{equation}
Moreover, 

\sli either $w(\zeta)$ vanishes identically and then $u_2$ is a reparametrization
of $u_1$,

\slii or, vector $\wect_0 \deff w(0)\not=0$ is orthogonal to $\vect_0$.
\end{lem}
\proof We work with the trivial bundle $E=\Delta\times \rr^{2n}$ equipped with linear complex
structures $\tilde J_k, k=1,2$, where $\tilde J_k$ stands for the extension by reflection of
$J_k\deff J\circ u_k$. The extended sections $\tilde u_1$ and $\tilde u_2$ of $E$ 
satisfy the  $\dbar_{\tilde J_k}$-equations  $\dbar_{\tilde J_k}\tilde u_k = 
(\partial_x + \tilde J_k\partial_y)\tilde u_k = 0$. Since $\tilde J_k$ are Lipschitz-continuous
the extensions $\tilde u_k$ are $\calc^{1,\alpha}$-regular. Without loss of generality we 
suppose that $\tilde u_1$ has no critical points, possibly except $0$. In this generality 
the comparison relation \eqqref(9.1.1) for $\tilde u_k$ was proved in \cite{IS4}. 
More precisely: there exists a holomorphic function $\psi(\zeta)$ in a neighbourhood of zero 
of the form $\psi (\zeta) = \zeta + O(\zeta^2)$ such that 
\begin{equation}
\eqqno(comp-h)
\tilde u_2(\zeta) = \tilde u_1(\psi(\zeta)) + \zeta^{\nu}\tilde w(\zeta)
\end{equation}
for some $\tilde w\in L^{1,p}(\Delta_r,E)$ with the properties as stated except of the reality 
conditions. Therefore all we need is to check that in our setting $\psi (\zeta)$ real-valued 
for real $\zeta$ and that $\tilde w(\zeta)$ takes values in $\rr^n$ for real $\zeta$. After
that we can take $w = \tilde w|_{\Delta^+_r}$. In order to achieve this we must examine the
proof in \cite{IS4}. Without loss of generality we can assume that $\vect_0 = e_1=(1,0,...,0)$.

\smallskip\noindent{\sl Step 1.} {\it The image $E_1$ of the differential $d\tilde u_1: T\Delta 
\to E$ is a well-defined $\tilde J_1$-complex line subbundle of the complex bundle $(E, \tilde 
J_1)$ over $\Delta\setminus\{0\}$. It extends to a $\tilde J_1$-complex line subbundle of $E$ 
of regularity $L^{1,p}$ over $\Delta$ such that $d\tilde u_1:T\Delta\to E_1$ is 
$L^{1,p}$-regular. }

\smallskip For the proof we refer to the {\slsf Claim 1} on the page 1175 of \cite{IS4}.
Since $\tilde u_1$ is $\rr^n$-valued when restricted to $\db\Delta^+$ we have that $E_1\cap 
\left(\db_0\Delta^+\times \rr^n\right)$ is a one-dimensional \emph{real} sub-bundle of 
$\db_0\Delta^+\times \rr^n$. Fix a $J_1$-hermitian metric $h^0$ on $E$ and set 
\[
h_{\zeta}(\vect , \wect) \deff \frac{1}{2}\left[h^0_{\zeta}(\vect , \wect) + 
\overline{h^0_{\bar\zeta}(\bar\vect , \bar\wect)}\right].
\]

This gives us a $\ti J_1$-hermitian
metric on $E$, which is invariant under the conjugation $\tau \deff (\tau_{\Delta}, 
\tau_{\cc^n})$ in $E=\Delta\times \cc^n$ as above. Indeed, one obviously has that 
\[
h_{\bar\zeta}(\bar\vect , \bar\wect) = \overline{h_\zeta (\vect , \wect)},
\]
and, moreover,
\[
\begin{split}
2& h\left(\ti J_1(\vect ), \wect \right) = h^0_\zeta \left(\ti J_{1,\zeta}(\vect ), \wect 
\right)+ \overline{h^0_{\bar \zeta} \left(\overline{\ti J_{1,\zeta}(\vect )}, \bar \wect 
\right)} = 
\\[2pt]
& = ih^0_\zeta (\vect , \wect) + \overline{h^0_{\bar\zeta}\left(-\ti J_{1,\bar\zeta}(\bar 
\vect ), \bar \wect\right)} = ih^0_\zeta (\vect , \wect) + i\overline{h^0_{\bar\zeta}
\left(\bar \vect , \bar \wect\right)} = i 2h(\vect ,\wect).
\end{split}
\]

\smallskip Take as $E_2$ the orthogonal complement to $E_1$ in $(E,\ti J_1, h)$.
Since $E_1$ is $\ti J_1$-complex the same is true for $E_2$. Moreover, since 
$d\ti u_1$ is invariant under conjugation the same is true for its image 
$E_1$ and therefore for $E_2$. Moreover, since $\tilde u_1$ is $\rr^n$-valued
 when restricted to $\db_0\Delta^+$ we have that 
$E_1\cap \left(\db_0\Delta^+\times \rr^n\right)$ is a one-dimensional \emph{real}
sub-bundle of $\db_0\Delta^+\times \rr^n$ with $E_1\cap \left(\{0\}\times \rr^n
\right) = \<e_1\>$. In addition we have that $\tilde J_1|_{\{\zeta\}
\times\rr^n} = J\st$ for real $\zeta$. Denote by $\Delta\times \cc\times \{0\}$ sub-bundle 
of $\Delta \times\cc^n$ with the fibre consisting 
of vectors of the form $(a, 0,\ldots,0)$ and  by $\Delta\times\{0\}\times \cc^{n-1}$ 
the subbundle of $\Delta\times \cc^n$ with the fibre consisting of vectors of the 
form $(0, a_2,\ldots,a_n)$. 

\smallskip With all these data at hand we can  construct an $L^{1,p}$-regular 
$(J\st, \ti J_1)$-linear isomorphism $\Phi: (\Delta\times \cc^n,J\st) \to (E, \ti J_1)$ 
such that:

\begin{itemize}
 \item the subbundle $\Delta\times\cc\times\{0\}\subset \Delta\times\cc^n$  is 
 mapped by $\Phi$ to $E_1$;

\smallskip\item $\db_0\Delta^+\times \rr^n$ is mapped by $\Phi$ to $\db_0\Delta^+\times \rr^n$;

\smallskip\item the real subbundle $\db_0\Delta^+\times\rr\times\{0\}$ is 
mapped by $\Phi$ to $E_1\cap \left(\db_0\Delta^+ \times \rr^n\right)$;

\smallskip\item $\Phi|_{\{0\}\times \rr^n}: \{0\}\times \rr^n \to \{0\}\times \rr^n$ is the 
identity map.

\smallskip\item $\Phi$ is invariant under the conjugation $\tau :(\zeta ,z)
\to (\bar\zeta, \bar z)$ on $\Delta\times \cc^n$ and on $E$.
\end{itemize}

\smallskip For $\zeta\in \Delta$ we denote the restriction $\Phi (\zeta , \cdot):
\{\zeta\}\times \cc^n \to E_{\zeta}$ by $\Phi_{\zeta}$.
On the trivial bundle $(\Delta\times \cc^n, J\st)$ we consider the natural 
coordinates $(\lambda, \zeta, w)$, where $\lambda  \in \Delta $, $\zeta \in \cc$, $w\in 
\cc^{n-1}$ and consider the following ``exponential'' map $\exp: \Delta \times 
\Delta \times \cc^{n-1} \to E$:
\begin{equation}
\eqqno(exp)
\exp : (\lambda, \zeta, w) \mapsto \big(\lambda, \tilde u_1(\zeta) + \Phi_{\lambda}[w]\big).
\end{equation}
This map $\exp$ is well-defined, $L^{1,p}$-regular in $\lambda$, $L^{2,p}$-regular in $\zeta$ 
and linear in $w$. In particular $\exp$ is continuous in $(\lambda,\zeta , w)$.
Moreover, for a fixed $\lambda \neq 0\in \Delta$ the linearization of
$\exp_{\lambda} := \exp(\lambda, \cdot,\cdot)$ with respect to variables 
$\zeta , w$ at $\zeta = \lambda, w=0$ is an isomorphism between $T_{\zeta }\Delta
\oplus \cc^{n-1}$ and $E_{\lambda}$. Therefore for $\lambda\neq 0$ the map $\exp_{\lambda}$ 
is an $L^{2,p}$-regular diffeomorphism of some neighbourhood $U_{\lambda}\subset\{\lambda\} 
\times \Delta\times\cc^{n-1}$ of the point $(\zeta =\lambda,0)$ onto some neighbourhood $V_{\lambda}$ 
of the point $\tilde u_1(\lambda)$ in $E_{\lambda}$.  It is important for us to make the following
\begin{rema} \rm
\label{exp-real}
For a real $\lambda\not=0$ mapping $\exp_{\lambda}$ is a diffeomorphism of intersections 
$U_{\lambda}\cap\left(\{\lambda\}\times \db_0\Delta^+\times \rr^{n-1}\right)$ and
$V_{\lambda}\cap \rr^n$. This follows from the fact that $\tilde u_1(\zeta)$ is $\rr^n$-valued
for real $\zeta$ and that $\db_0\Delta^+\times\rr^n$ is mapped by $\Phi$ to $\db_0\Delta^+\times\rr^n$.
\end{rema}

We need to estimate the size of $V_{\lambda}$. In order to do so let us consider the rescaled 
maps
\[
\tilde u_1^t(\zeta) := t^{-\mu}\, \tilde u_1(t\zeta),
\]
where $t\in (0,1]$. Notice that the family $\{\tilde u_1^t\}_{t\in (0,1]}$ is uniformly bounded 
with respect to the $L^{2,p}$-norm and its limit  $\lim_{t\searrow 0} \tilde u_1^t(\zeta)$ in 
$L^{2,p}$-topology is the map $\tilde u_1^0(\zeta) := \vect_0 \zeta^{\mu}$. Indeed, the 
$\calc^0$-convergence $\tilde u_1^t\rightrightarrows \tilde u_1^0$ is clear from the representation
\eqqref(nrm-frm2). To derive from here the $L^{2,p}$-convergence remark that $\tilde u_1^t$ is
$\tilde J_t$-holomorphic with respect to the structure $\tilde J_t(\zeta)\left[\cdot\right]
:=\tilde J_1(t\zeta)\left[t^{\mu}\cdot\right]$, 
$\zeta\in \Delta$, and that $\tilde J_t$ converge to $J\st$ in the Lipschitz norm. This implies the 
$L^{2,p}$-convergence.

\smallskip Now for $t\in [0,1]$ define the rescaled exponential maps
\begin{equation}
\exp^t_{\lambda}(\zeta, w) := \tilde u_1^t(\zeta ) + \Phi_{t\lambda}[w]. 
\eqqno(exp-t)
\end{equation}

\smallskip\noindent{\sl Step 2.} {\it 
There exist constants $c^*,c_1, \eps >0$ such that for every $\lambda\in \{ |\lambda|=\eps \}$ and 
$t\in [0,1]$ mapping $\exp^t_{\lambda}(\zeta, w)$ is an $L^{2,p}$-regular diffeomorphism of 
$U_{\lambda}= \{ (\zeta ,w): |\zeta -\lambda |<c_1, |w|<c_1\}$ and a neighbourhood $V_{\lambda}^t$ 
of $\tilde u_1^t(\lambda)$ in $E_{\lambda}$ which contains the ball $\{ |z - \tilde u_1(\lambda)| < c^*\}$. 
Moreover, the inverse maps $(\exp^t_{\lambda})\inv:V^t_{\lambda}\to U_{\lambda}$ are $L^{2,p}$-regular, 
their $L^{2,p}$-norms are bounded by a uniform constant independent of  $\lambda\in \{ |\lambda|=
\eps \}$ and $t$, and the dependence of $(\exp^t_{\lambda})\inv$ on $\lambda$ is $L^{1,p}$.
}

\smallskip 
This statement readily follows from the fact that $K=\{ |\lambda|=\eps \}\times \{ t
\in [0,1]\}$ is a compact, the function $\exp^t_{\lambda}(\zeta,w)$ is a local
$L^{2,p}$-regular diffeomorphism for every fixed $(\lambda,t)\in K$, depends
continuously on $t\in[0,1]$ and is $L^{1,p}$-regular on $\lambda\in\{|\lambda|=\eps\}$ with 
respect to the $L^{2,p}$-topology. In the sequel without loss of generality we may assume that 
$\eps = 1$. This can be always achieved by an appropriate rescaling.

\smallskip\noindent{\sl Step 3.} {\it 
For every $\lambda\in \Delta\bs\{0\}$ there exists a neighbourhood $V_{\lambda}\ni
\tilde u_1(\lambda)$ containing the ball $B(\tilde u_1(\lambda), c^*\cdot |\lambda|^{\mu})$ with 
the constant $c^*$ independent of $\lambda$, such that $\exp_{\lambda}$ is an 
$L^{1,p}$-regular  homeomorphism between some neighbourhood $U_{\lambda}$ of $(\lambda,0)$ in the 
fiber  $\{\lambda\}{\times}\cc^n$ and $V_{\lambda}$. 
}

Here by  $L^{1,p}$-regular homeomorphism we understand a  homeomorphism
which is  $L^{1,p}$-regular and its inverse is also  $L^{1,p}$-regular.
In order to prove this statement fix some $0<|\lambda|<\frac{1}{2}$ and
set $\tilde \lambda=\frac{\lambda}{|\lambda|}, \tilde \zeta = \frac{\zeta}{|\lambda|},
\tilde w = \frac{w}{|\lambda|^{\mu}}, t=|\lambda|$. Take $c^*$ as in Step2. Then we 
have a homeomorphism
\[
\exp^t_{\tilde\lambda} : \{ |\tilde\zeta - \tilde\lambda|<c_1, |\tilde
w|<c_1\} \to V_{\tilde\lambda}^t\supset \{ |\tilde\zeta - \tilde u_1^t(\tilde\lambda)|
<c^*\} .
\]
But $\exp^t_{\tilde\lambda}(\tilde\zeta ,\tilde w) =
t^{-\mu}\tilde u_1(t\tilde\zeta) + \Phi^{-1}(t\tilde\lambda)\tilde w = t^{-\mu}[
\tilde u_1(\zeta ) + \Phi^{-1}(\lambda)w] = t^{-\mu}\exp_{\lambda}(\zeta ,w)$ and this
map is a homeomorphism between $\{
|\frac{\zeta}{|\lambda|}-\frac{\lambda}{|\lambda|}|<c_1, \frac{|w|}{|\lambda|^{\mu}}<c_1\}$
and some $V_{\tilde\lambda}^t$ containing $\{ |\tilde\zeta - \tilde u_1^t(\tilde\lambda)|
<c^*\}$. Therefore $\exp_{\lambda}$ is a homeomorphism between
\[
\{ |\zeta - \lambda|<c_1|\lambda|, |w|<c_1|\lambda|^{\mu}\} \leftrightarrow
V_{\lambda}\supset \{ |\zeta - \tilde u_1(\lambda)|<c^*|\lambda|^{\mu}\}.
\]
The step is proved.

\smallskip\noindent{\sl Step 4.} {\it There exists an $L^{1,p}$-regular function $\psi(\lambda)$ 
in a neighbourhood of zero of the form $\psi (\lambda) = \lambda + O(\lambda^2)$, real-valued for 
real $\zeta$, such that 
\begin{equation}
\eqqno(comp-h1)
\tilde u_2(\lambda) = \tilde u_1(\psi(\lambda)) + \lambda^{\nu}\tilde w(\lambda)
\end{equation}
for some $\tilde w\in L^{1,p}(\Delta,E_2)$ which takes values in $\rr^{n-1}$
for $\lambda\in \rr$. Moreover $\wect_0 \deff \tilde w(0)$ is non-zero and orthogonal to $\vect_0$.
}

Since $\tilde u_2(\lambda) - \tilde u_1(\lambda) = O(|\lambda|^{\mu +\alpha})$ we obtain that 
$\tilde u_2(\lambda)$ $\in B(\tilde u_1(\lambda),c^*\,|\lambda|^{\mu})$ for $\lambda$ small 
enough. Define  $(\zeta(\lambda),W(\lambda)):=\exp_{\lambda}^{-1}(\tilde u_2(\lambda))$
where $\exp_{\lambda}^{-1}:V_{\lambda}\to U_{\lambda}$ is the local inversion of the map 
$\exp_{\lambda}$ which exists by Step 3. Set $\psi (\lambda):=\zeta (\lambda)$, 
$\tilde w(\lambda):=\Phi^{-1}(\lambda)W(\lambda)$.  We obtain the desired relation
\begin{equation}
\eqqno(tild-u1-u2)
\tilde u_2(\lambda) = \tilde u_1(\psi(\lambda)) + \tilde w(\lambda). 
\end{equation}
Notice that $\psi (\lambda)$ (resp. $\tilde w(\lambda)$) is real-valued (resp. $\rr^{n-1}$-valued)
for real $\lambda$. 
This follows from the fact mentioned in Remark \ref{exp-real} that our exponential map preserves 
$\rr^n$ for real $\lambda$. Notice that $\psi$ and $\tilde w$ are $L^{1,p}_{loc}$-regular but only 
in $\Delta\setminus \{0\}$. 

\smallskip\noindent{\slsf Step 5.} {\it The estimate of the $L^{1,p}$-norm of both of them near the origin
proved in the {\slsf Claim 5} on the page 1177 of \cite{IS4} gives us the $L^{1,p}$-regularity of 
both functions in $\Delta$. } 

\smallskip\noindent{\slsf Step 6.} {The term  $\tilde w$ in \eqqref(tild-u1-u2) can be presented 
in the form $\lambda^{\nu} \tilde w$ .} This is essentially  proved on the pages 1177-1179 of 
\cite{IS4}. The arguments are similar to that which were used in the proof Lemma \ref{van-order}, 
\ie rely on certain $\dbar$-inequalities, which provide the existence of such $\nu$ if $\tilde w
\not\equiv 0$. 

 Consider the pulled-back bundles $E'_1 := \psi^*E_1$ and $E'_2:= \psi^*E_2$ over $\Delta$ as 
the subbundles of the trivial bundle $E=\Delta\times \rr^{2n}$.  Equip $E$ with the complex
structure $\tilde J_1':=\psi^*\tilde J_1 = \tilde J_{u_1\circ\psi }$. Let $\pr'_2$ be the
projection of $E$ onto $E'_2$ parallel to $E'_1$. Then \eqqref(comp-h1) implies
\begin{equation}
\eqqno(pr-strix)
\pr'_2\big( ( \partial_{\xi} + \tilde J_1' \cdot \partial_{\eta}) \tilde w(\zeta) \big) =
\pr'_2\big( ( \partial_{\xi} + \tilde J_1' \cdot \partial_{\eta}) (\tilde u_2(\zeta) - 
\tilde u_1(\psi(\zeta))\big).
\end{equation}
The second term $(\partial_{\xi} + \tilde J_1' \cdot \partial_{\eta})u_1(\psi(\zeta))$
on the right hand side, of \eqqref(pr-strix) is the Cauchy-Riemann operator $\dbar_{\tilde J_1'}$ 
applied to the composition $\tilde u_1' := \tilde u_1 \circ \psi$.  Notice that for an 
$L^{1,p}$-function $\psi :\Delta\to \Delta$ and a $J$-holomorphic map $u:\Delta\to \cc^n$ 
we have
\begin{equation}
\eqqno(debar-psi)
\dbar_J(u\circ\psi) = du\circ \dbar\psi, 
\end{equation}
where $\dbar\psi$ is the standard $\dbar$-derivative of the function $\psi$. 
Indeed, the expression $\dbar_{J}(u \circ \psi)$ is the $J$-antilinear part of
the differential $d(u \circ \psi) = du \circ d\psi$. Since $du$ is
$J$-linear,  the antilinear part of $du\circ d\psi$ will be $du$ of
the antilinear part of $d\psi$ which is $\dbar\psi$. Relation \eqqref(debar-psi)
follows. In our case this gives

\begin{equation}
\dbar_{\tilde J_1'}(\tilde u_1\circ\psi ) = d\tilde u_1\circ \dbar\psi. \eqqno(3.6)
\end{equation}

\smallskip
Further, observe that $d\tilde u_1 \circ \dbar\psi$ takes values in the
pulled-back $E'_1 = \psi^*E_1$. Therefore \break $\pr'_2( \dbar_{\tilde J_1'}(\tilde u_1 \circ
\psi))$ vanishes identically. In order to estimate $( \partial_{\xi} 
+ \tilde J_1' \cdot \partial_{\eta})\tilde u_2$ subtract from it $(\partial_{\xi} + 
\tilde J_{u_2} \cdot \partial_{\eta})\tilde u_2 = \dbar_{\tilde J_2}\tilde u_2 =  0$ 
and obtain
\[
(\tilde J_1' - \tilde J_2)\db_{\eta}\tilde u_2 = (\tilde J_{u_1 \circ \psi} -  
\tilde J_{u_2})\cdot \partial_{\eta}u_2.
\]
The $L^{2,p}$-regularity for $\tilde u_2$ provides the $L^{1,p}$-regularity
of $\partial_{\eta}\tilde u_2$, whereas the Lipschitz condition on $J$ yields
the pointwise estimate
\[
\big| J \circ u_1 \circ \psi(\zeta) -  J\circ u_2(\zeta) \big|\leq
Lip(J) \cdot |w(\zeta)|,
\]
which implies the estimate 
\begin{equation}
\eqqno(pr-chast)
\big| \tilde J_{u_1 \circ \psi}(\zeta) -  \tilde J_{u_2}(\zeta) \big|\leq
Lip(J) \cdot |\tilde w(\zeta)|.
\end{equation}
Indeed, 
\[
\begin{split}
& \left|\tilde J_{u_1}(\zeta)[\vect] -\tilde J_{u_2}(\zeta) [\vect]\right| = 
\left|-\overline{J_u(\bar\zeta)[\bar\vect]} + \overline{J_{u_2}(\bar\zeta)[\bar\vect]}\right|
=  
\\[2pt]
& = \left|\left[J(u_2(\zeta)) - J(u_1(\bar\zeta))\right][\bar\vect]\right| 
\le \norm{J}_{\calc^{Lip}}|u(\bar\zeta)|\,|\vect|.
\end{split}
\]
Therefore the right hand side of \eqqref(pr-strix) admits the pointwise estimate
by $h(\zeta)\cdot|\tilde w(\zeta)|$ with some $h\in L^p(\Delta)$.

\smallskip %
Now, let us rewrite the left hand side $\pr'_2 \big ( ( \partial_{\xi} + J \circ u_1
\circ \psi \cdot \partial_{\eta}) \tilde w\big)$ of \eqqref(pr-strix) as a $\dbar$-type 
operator of $\tilde w$. Consider the restriction $\pr'_2|_{E_2} : E_2 \to E'_2$ of the
projection $\pr'_2$ onto $E_2$. Using the facts that $\psi(\zeta)$ is
continuous and $\psi(0) =0$ we conclude that $\pr'_2(\zeta) : (E_2)_{\zeta}
\to (E'_2)_{\zeta} = (E_2)_{\psi(\zeta)}$ is a bundle isomorphism over a
sufficiently small disk $\Delta_r$. Denote by  $(\pr'_2)\inv(\zeta)$ its inverse.
Setting $\ti w'(\zeta) := \pr'_2\tilde w(\zeta)$ we obtain a pointwise estimate
\begin{equation}
\eqqno(w-til-w)
1/C\cdot |\tilde w(\zeta)| \leq |\ti w'(\zeta)| \leq C\cdot |\ti w(\zeta)|
\end{equation}
with uniform constant $C$ for $\zeta$ in the disk $\Delta_r$.
Denote $\nabla_{\xi} :=  \pr'_2 \circ\partial_{\xi} \circ (\pr'_2)\inv$ and
$\nabla_{\eta}:= \pr'_2 \circ \partial_{\eta} \circ (\pr'_2)\inv$.  Using this
we obtain
\begin{equation}
\eqqno(nabla-h)
\pr'_2 \big ( ( \partial_{\xi} + \tilde J_1' \cdot\partial_{\eta}) \tilde w\big) = 
( \nabla_{\xi} + \tilde J_1' \cdot\nabla_{\eta}) (\pr'_2\tilde w) + H(\ti w)
\end{equation}
with some $L^p$-regular endomorphism $H$. Comparing \eqqref(nabla-h) and 
\eqqref(pr-chast)
with \eqqref(pr-strix) we conclude a pointwise differential inequality
\begin{equation}
\eqqno(debar-nab)
\big|( \nabla_x + \tilde J_1' \cdot \nabla_y) (\tilde w)\big| \leq h|\tilde w|
\end{equation}
with $J'_1 :=  J \circ u_1 \circ \psi$ and an  $L^p$-regular function $h$.
If we fix some $L^{1,p}$-trivialization $e_1,...,e_{n-1}$ of
$E_2^{'}$ and remark that in (any) such trivialization $\nabla_{\xi} = \db_{\xi}
+ R_{\xi}$
and $\nabla_{\eta} = \db_{\eta} + R_{\eta}$ with some $R_{\xi}, R_{\eta}\in 
L^p(\Delta , \endo (E_2^{'}))$ we obtain the following estimate
\begin{equation}
\big | (\partial_x + J'_1 \partial_y)\ti w \big| \leq h\cdot |\ti w|.
\eqqno(3.7)
\end{equation}

\smallskip Now we can apply  Lemma \ref{dbar-in1} and conclude that:

\begin{itemize}
\item either $w(\zeta)$ vanishes identically 
($\Rightarrow w\zeta)\equiv 0$)
and then the image $u_2(\Delta_r)$ is contained in $u_1(\Delta)$,

\item or, $w(\zeta) = \zeta^{\nu} f(\zeta)$ 
for some $f(\zeta)\in L^{1,p}(\Delta_r, \cc^{n-1})$ with $f(0) \neq 0$.
\end{itemize} 
 
The integer $\nu$ must be bigger than $\mu$ because
$u_2(\zeta)-u_1(\zeta)=o(|\zeta|^{\mu +\alpha})$.  Since the projection $w(\zeta)
\mapsto \ti w(\zeta) := \pr'_2(w(\zeta))$ is an $L^{1,p}$-regular
isomorphism, we obtain the same structure for $w$. Finally, observe
that $f(0)$ lies in the fiber $(E_2)_0$ which is orthogonal to 
$(E_1)_0 = \cc \vect_0$ since due to the condition $J(0)=J\st$ we were able to 
choose $\Phi_0 = \id$. Therefore we obtain

\begin{equation}
u_2(\zeta) = u_1(\psi (\zeta)) + \zeta^{\nu} w(\zeta), \eqqno(3.8)
\end{equation}

\noindent where $\nu >\mu$ and $\wect_0=w(0)$ is orthogonal to $\vect_0$.

\begin{rema} \rm
\label{trans-comp}
Let us turn the attention of the reader that we could choose as $E_2$ {\slsf any} 
$J_1$-complex subbundle of $E$ transverse to $E_1$, and our $w$ will take
its values in this $E_2$. This will affect the reparameterization function 
$\psi$.
\end{rema}

\bigskip
\smallskip\noindent{\slsf Step 7.} {\it There exists a holomorphic $\psi$ of the form 
$\psi (\zeta) = \zeta + O(\zeta^2)$ satisfying \eqqref(comp-h1), which is again real 
for real $\zeta$.}

\smallskip First consider the case when $\tilde w\equiv 0$. Write for $\zeta\in \Delta^+$ 
using the equation of $J$-holomorphicity
\[
0 = \dbar_{J}(u_1\circ \psi ) = du_1\circ \dbar\st \psi , 
\]
which implies that $\psi$ is holomorphic on $\Delta^+$. Since it is real on $\db_0\Delta^+$ 
it is holomorphic in a neighbourhood of zero. Otherwise we shall construct recursively 
complex polynomials $\psi_j$, $j=1,...,l$ obeying the following conditions:

\begin{itemize}
 \item $\psi_j(\zeta) = \zeta + O(\zeta^2)$ and $\psi_j(\zeta)$ is real for real $\zeta$;
 
\smallskip\item $u_2(\zeta) = u_1(\psi_j(\zeta)) + \zeta^{\nu_j}w_j(\zeta)$ for some
$w_j\in L^{1,p}(\Delta^+ , \cc^n)$, which is  $\rr^n$-valued on $\db_0\Delta^+$ and 
$\nu_j\in \nn$;

\smallskip\item $\mu <\nu_1<...<\nu_l$, for $1\le j<l$ vectors $w_j(0)$ are collinear with 
$\vect_0$, but $w_l(0)$ is not.
\end{itemize}

\smallskip Take $\psi_1(\zeta) = \zeta$. Let us see that there exists an integer $\nu_1 > 
\mu$ and $w\in L^{1,p}(\Delta, \cc^n)$, $w(0) \neq 0$ such that $u_1(\zeta) - 
u_2(\zeta) = \zeta^{\nu_1} v(\zeta)$. Notice that such $w(\zeta)$ is necessarily 
$\rr^n$-valued for real $\zeta$. Set $u=u_1-u_2$ and let us compute $\dbar_{J\circ u_1}(u) = 
(\partial_{\xi} + J(u_1(\zeta))\partial_{\eta})u(\zeta)$:
\[
\dbar_{J \circ u_1}(u) = (\partial_x + J(u_1)\cdot \partial_y) (u_1
-u_2) = (\partial_x + J(u_1)\cdot \partial_y) (u_1 -u_2) +
(\partial_x + J(u_2)\cdot \partial_y) u_2 =
\]
\[
= (J(u_2)\cdot
\partial_y - J(u_1)\cdot
\partial_y)u_2 = (J(u_1 +u) - J(u_1))\cdot
\partial_yu_2.
\]
By the Lipschitz regularity of $J$ and $\partial_yu_2 \in L^p(\Delta)$ we obtain a 
pointwise differential inequality
\[
|\dbar_{J \circ u_1}(u)(\zeta)| \leq h(z)\cdot |u(\zeta)|
\]
for some $h\in L^p(\Delta)$. Now we apply  Lemma \ref{dbar-in1} and, taking into account
that $u_1(\zeta) \not \equiv u_2(\zeta)$, we get that $u_1(\zeta) - u_2(\zeta) = \zeta^{\nu_1}
w_1(\zeta)$ with $w_1(\zeta)$ as stated. Since $u_1(\zeta) - u_2(\zeta) = O(|\zeta|^{\mu + \alpha}$
we conclude that $\nu_1 >\mu$.  If $w_1(0)$ occurred to be not proportional to $\vect_0$ we are 
done with $l=1$.

\smallskip 
Assume that we have constructed such sequences $\mu < \nu_0 < \nu_1 <\ldots
\nu_k$, and $w_1(\zeta), \ldots,$ $ w_k(\zeta)$, $\psi_1,...,\psi_k$ and that $w_1(0), 
\ldots, w_k(0)$ are proportional to $\vect_0$. We assume that the reality conditions
holds as well. 

\smallskip Remark that for a function $f$ of class $\calc^{1,\alpha}$ in the unit ball $B$ one has 
\[
f(x_0 + a) = f(x_0) + df(x_0)(a) + O(|a|^{1+\alpha}) \quad\text{ with } \quad
O(|a|^{1+\alpha})\le C\norm{f}_{\calc^{1,\alpha}(B)}|a|^{1+\alpha},
\]
where constant $C$ doesn't depend on $f$. Using this and the fact that 
$du_1(\zeta) = \mu \zeta^{\mu -1} \vect_0 + O(|\zeta|^{\mu -1 +\alpha})$ for
any $0 <\alpha <1$, see \eqqref(nrm-frm3), we can write
\[
u_1(\psi_k(\zeta) + a\zeta^m) = u_1(\psi_k(\zeta))  + du_1(\psi_k(\zeta))
(a\zeta^m) + O(\norm{u_1}_{\calc^{1,\alpha}(\Delta)}|\zeta|^{m + m\alpha}) = 
\]
\[ 
= u_1(\psi_k(\zeta)) + \mu \psi_k(\zeta)^{\mu-1}\cdot \psi_k'(\zeta)\cdot 
a\zeta^m \cdot  \vect_0 + O(|\zeta|^{m+\mu -1 +m\alpha}) =
\]
\[
=u_1(\psi_k(\zeta)) + \mu \zeta^{m+\mu-1} \cdot a\cdot  \vect_0 +
O(|\zeta|^{m+\mu}) \quad\text{ if } \quad m\alpha \ge 1.
\]
Take $m := \nu_k -\mu +1$, notice that $m\ge 2$ and since we can take $\alpha$ 
close to $\adyn$ the condition $m\alpha \ge 1$ will be satisfied. Now compute 
$a$ from the relation
\begin{equation}
\eqqno(comp-a)
\mu \zeta^{m+\mu-1} \cdot a\cdot  \vect_0 - \zeta^{\nu_k} w_k(0) =0.
\end{equation}
Notice that it will lie in $\rr^n$. Setting
\begin{equation}
\eqqno(set-psi)
\psi_{k+1}(\zeta) := \psi_k(\zeta) + a\zeta^m
\end{equation}
we obtain 
\[
u_2(\zeta) - u_1(\psi_{k+1}(\zeta)) = u_2(\zeta) - u_1(\psi_k(\zeta) + a\zeta^m) 
= u_2(\zeta) - u_1(\psi_k(\zeta)) - 
\]
\begin{equation}
\eqqno(get1)
- \mu \zeta^{m+\mu +1}a\vect_0 + O(\zeta^{m+\mu}) = \zeta^{\nu_k}w_k(\zeta) - 
\zeta^{\nu_k}w_k(0) + O(\zeta^{\nu_k+1}) = O(|\zeta|^{\nu_k +\alpha}).
\end{equation}
Exactly as in the starting case $k=1$ we obtain for $u_2(\zeta)-u_1(\psi_{k+1}(\zeta))$ 
a new $\nu_{k+1} >  \nu_k$ and a new $w_{k+1}(\zeta)$. Reality conditions obviously
hold. Compare the obtained presentations $u_2(\zeta) = u_1(\psi_j(\zeta)) + 
\zeta^{\nu_j} w_j(\zeta)$ with the decomposition \eqqref(comp-h1) and observe
that $\nu_j$ cannot be bigger than $\nu$. This implies that at some step we 
obtain $\nu_l = \nu$. At this step $w_l(0)$ is not proportional to 
$\vect_0$ and the recursive procedure halts.

\smallskip To make $w_l(0)$ orthogonal to $\vect_0$, \ie $w_l(0)\in (E_2)_0$
we must repeat the inductive step once again. Decompose $w_l(0) = w_l(0)^{||} + 
w_l(0)^{\bot}$, where $w_l(0)^{||}$ is proportional to $\vect_0$ and $w_l(0)^{\bot}$
is orthogonal. If $w_l(0)^{||}\not= 0$ take $a$ in \eqqref(comp-a) from the relation 
\begin{equation}
\eqqno(comp-a1)
\mu \zeta^{m+\mu-1} \cdot a\cdot  \vect_0 - \zeta^{\nu_l} w_l(0)^{||} =0.
\end{equation}
Then set
\begin{equation}
\eqqno(set-psi1)
\psi_{l+1}(\zeta) := \psi_l(\zeta) + a\zeta^m \quad\text{ with }\quad m = \nu_l -\mu +1
\end{equation}
and obtain  as in \eqqref(get1)
\begin{equation}
\eqqno(get2)
u_2(\zeta) - u_1(\psi_{l+1}(\zeta)) = \zeta^{\nu}w_l(\zeta) - 
\zeta^{\nu_l}w_l(0)^{||} + O(\zeta^{\nu_l+1}) \zeta^{\nu}w_{l+1},
\end{equation}
where $w_{l+1}(\zeta) \deff  w_l(\zeta) - w_l(0)^{||} + O(\zeta)$ and therefore 
$w_{l+1}(0)$ is orthogonal to $\vect_0$. Comparison Lemma is proved.

\smallskip\qed

\newprg[COMP.comp-m]{Comparison of meeting analytic half-disks} 

 Now we consider $J$-holomorphic mappings $u_k:(\Delta^+, $ $\db_0\Delta^+,0) 
\to (\rr^{2n}, \rr^n,0)$ for a Lipschitz-continuous that have the opposite  tangent 
vectors at zero and the  order of vanishing equal to $\adyn $,~\ie
\begin{equation}
\eqqno(tan-u1-u2)
\begin{split}
& u_1(\zeta) = \vect_0 \zeta + O(|\zeta|^{1+\alpha}),
\\[2pt]
& u_2(\zeta) = - \vect_0 \zeta + O(|\zeta|^{1+\alpha}).
\end{split}
\end{equation}
In that case we say that $u_k$-s meet each other at zero. For the proof of Theorem 
\ref{b-pos-thm} for meeting half-discs we need a version of 
the {\sl Comparison Lemma} for this case. Let us start with one preliminary remark.  
Take an almost complex structure $J$ in $\rr^{2n}$,
and define the new structure as 
\begin{equation}
\eqqno(j-minus)
J^-(z)[\vect] \deff -\overline{J(\bar z)[\bar\vect]}, \qquad\text { also set } \qquad
J^+ \deff J.
\end{equation}
We assume that $J|_{\rr^n} = J\st$ and therefore $J^-|_{\rr^n} = J\st$. In $\rr^{2n}$
we consider the $J\st$-complex coordinates $(z,...,z_n) = (x_1,...,x_n) + i (y_1,...,y_n)$,
where $\rr^n$ stands for the $J\st$-totally real subspace $\{y=0\}$. Consider half-spaces 
$\rr^{2n}_+ \deff \{y_1 \ge 0\}$ and $\rr^{2n}_- \deff \{y_1 \le 0\}$, ans cones  
$\calc^{\pm} =\{z\in R^{2n}_{\pm}: |y_1| >C\cdot |y'|\}$. Here $y'=(y_1,...,y_{n-1})$ and 
$C>0$ some constant to be clear from the context.

\begin{lem}
\label{cones}
There exists an almost complex structure $J^{\vee}$ on $\rr^{2n}$ which is invariant under 
the conjugation $\tau (z) = \bar z$ and such that $J^{\vee}|_{\calc^{\pm}} = J^{\pm}$. If,
in addition, $J$ was Lipschitz-continuous then $J^{\vee}$ is Lipschitz-continuous as well.
\end{lem}
\proof We start with a portion of the linear algebra. Let 
\begin{equation}
\eqqno(1.5.3)
\omega\st =\sum_{i=1}^n dx_i \land dy_i
\end{equation}
be the standard symplectic form in $\rr^{2n}$, $J\st$ the standard complex structure and 
$\jj_{\omega\st} $ the manifold of $\omega\st$-tamed almost complex structures. 
Notice that $J\st \in \jj_{\omega\st}$ and that for any $J \in \jj_{\omega\st}$ the operator 
$J + J\st$ is invertible. Indeed, otherwise we have $Jv = - J\st v$ for some non-zero 
$v\in \rr^{2n}$ and then
$0< \omega(v, Jv) = -\omega(v, J\st v) <0$, a contradiction.

\smallskip For $J \in \jj_{\omega\st}$  set $\sf L(J) \deff -(J - J\st) (J+
J\st)\inv$. The equivalent definitions are
\begin{equation}\eqqno(1.2.4)
\arraycolsep=1pt
\begin{array}{r@{\,}clcl}
\sf L (J) &=&  -(J - J\st) (J+J\st)\inv &=&
-\bigl((J - J\st)J\st\bigr) \cdot \bigl( (J+J\st) J\st \bigr)\inv \\
\vph&=&({\bf 1} + JJ\st )({\bf 1} - JJ\st  )\inv
     &=& ({\bf 1} - JJ\st)\inv ({\bf 1} + JJ\st).
\end{array}
\end{equation}

One easily checks that  
\begin{itemize}
    \item $W= \sf L(J)$ is $J\st$-antilinear, \ie $W J\st= - J\st W$.

    \item ${\bf 1} - W^\sft W \gg 0$ for $W=
\sf L(J)$. 
\end{itemize}

\smallskip
Therefore $\sf L$ maps $\jj_{\omega\st}$ into the set
\begin{equation}
\eqqno(1.5.5)
\scrw \deff
\{\, W \in \endo(V)\;:\; W J\st =- J\st W,\; {\bf 1} -W^\sft W \gg 0 \,\}.
\end{equation}

Explicit calculation show that $\sf L: \jj_{\omega\st} \to \scrw$ is a diffeomorphism,
the inverse map should be given by
\begin{equation}
J= \sf K(W) \deff J\st ({\bf 1} + W) ({\bf 1} - W)\inv=
J\st ({\bf 1} - W)\inv ({\bf 1} + W) =
J\st \frac{{\bf 1} + W}{{\bf 1} - W}
\eqqno(1.2.6)
\end{equation}

I.e., one checks that for $W\in \scrw$  operator $J= \sf K(W) = J\st \frac{{\bf 1} + W}{{\bf 1} - W}$
is well-defined, $J^2 = -{\bf 1}$ and $J \in \jj_{\omega\st}$. Finally one easily observes that 
$\scrw $ is convex as well as $W\st \deff \sf L(J\st)$ is zero. 




\smallskip The proof of our lemma is linearized. Let $W^{\pm}$ be the images of 
$J^{\pm}$ under $\sf L$. Take a cut-off function $\chi$ on the unit sphere $\ss^{2n-1}$ 
in $\rr^{2n}$ which is equal to $\adyn$ on $\calc^+\cap \ss^{2n-1}$ and zero on 
$\calc^+\cap \ss^{2n-1}$. Extend it independently of $r$ and set 
\[
W^{\vee} \deff \chi \cdot W^+ + (1-\chi ) W^-.
\]
Since $W^{\pm}|_{\rr^n}\equiv 0$ and are Lipschitz-continuous the resulting $W^{\vee}$
will be Lipschitz-continu\-ous. Indeed it is continuous at zero and, since derivatives of 
$\chi$ on $\rr^{2n}\setminus \{0\}$ grow as $1/r$ and $W^{\pm}$ decreases as $r$ we 
see that $W^{\vee}\in L^{1,\infty}(\rr^{2n}) = \calc^{lip}(\rr^{2n})$. All what is 
left is to set $J^{\vee} \deff \sf K(W^{\vee})$.

\smallskip\qed

\begin{lem} {\slsf (Comparison Lemma - II).}
\label{comp-lem-op}
Let $u_k:(\Delta^+, \db_0\Delta^+,0) \to (\rr^4,\rr^2,0)$ be two analytic half-disc 
meeting at zero. Almost complex structure $J$ is assumed to be Lipschitz-continuous 
and $J|_{\rr^n} = J\st$. Then there exists a holomorphic function $\psi$ of the form 
$\psi (\zeta) = \zeta+ O(\zeta^2)$, real-valued for real $\zeta$, 
an integer $\nu > 1$ and an $L^{1,p}_{loc}$-regular $\cc^n$-valued function 
$\tilde w(\zeta)$, $\rr^n$-valued for real $\zeta$, such that for some $r>0$  
\begin{equation}
\eqqno(comp-op)
\tilde u_2(-\zeta) = \tilde u_1(\psi(\zeta)) + \zeta^{\nu} \tilde w(\zeta).
\end{equation}
Moreover, 

\sli either $\tilde w(\zeta)$ vanishes identically and then $u_2(\Delta^+)\cup u_1(\Delta^+)$
is a $J$-complex disk,

\slii or, the vector $\wect_0 \deff \tilde w(0)\not=0$ can be chosen orthogonal to $\vect_0$.
\end{lem}
\proof We keep the notations of the previous lemma. Notice that, after shrinking of the 
disk (or, after rescaling) $u_1$ takes it values in $\rr^{2n}_+$ and $u_2$ in 
$\rr^{2n}_-$. Moreover, they take their values in the cones $\calc^+$ and $\calc^-$
correspondingly. After extension by reflection $\tilde u_1$ when restricted to 
$\Delta^-$ will take its values in $\calc^-$ and will be $J^-$-holomorphic there.
As a result $\tilde u_1$  as a map from $\Delta$ to $\rr^{2n}$ becomes
to be $J^{\vee}$-holomorphic. The same holds for $\tilde u_2(-\zeta)$ with 
the difference that it is $J^+$-holomorphic on $\Delta^-$ with values in $\calc^+$
and $J^-$-holomorphic on $\Delta^+$ with values in $\calc^-$. But anyway both 
are $J^{\vee}$-holomorphic on $\Delta$. 
Therefore the previous Lemma \ref{comp-lem-r} applies to $\tilde u_1(\zeta)$ and
$\tilde u_2(-\zeta)$ and gives \eqqref(comp-op).

\smallskip In the case when $\tilde w(\zeta) \equiv 0$ we have that 
$\tilde u_2(-\xi) \equiv \tilde u_1(\psi(\xi)) $ for $\xi \in (-1,1)$. This and 
the Cauchy-Riemann condition means that the (complex) tangents $T_{u_1(-\xi)}M_1$
and $T_{u_2(\psi(\xi))}M_2$ coincide, here $M_k$ are the half-discs parametrised by $u_k$. 
This means that $M_1\cup M_2 = u_2(\Delta^+)\cup u_1(\Delta^+)$ is a complex
$\calc^1$-manifold.

\smallskip\qed

\newsect[PERT]{Perturbations of Half-Disks and Proof of the Main Theorem}

We shall need to perturb analytic half-disks with estimates.

\newprg[REFL.pert]{Perturbations}

The item (\sliii of Theorem \ref{b-pos-thm} will be proved via certain perturbations.
But first let us state the following lemma.
\begin{lem}
\label{lip-cont}
Let $A$ be a Lipschitz-continuous $\End (\rr^{2n})$-valued function in the unit ball 
$B$ of $\rr^{2n}$ with $A(0)=0$ and let $u:(\Delta ,0)\to (B,0)$ be a $J$-holomorphic map 
with respect to a Lipschitz-continuous almost complex structure $J$ on $B$. Then for all 
integers $\nu$ and $\lambda$ satisfying $\nu \le \mu + \lambda -1$, where $\mu$ is order 
of vanishing of $u$ $\mu$ at zero, the function $\zeta^{-\nu}\cdot A(u(\zeta))
\cdot \zeta^{\lambda}$ is Lipschitz-continuous  and for every $0<r<1$ one has
the estimate 
\begin{equation}
\eqqno(est-A1)
Lip_{\Delta (r)}\big(\zeta^{-\nu}\cdot A(u(\zeta))\cdot \zeta^{\lambda}\big)\le 
C(p,r)\cdot Lip_B(A)\cdot\norm{u}_{L^{1,p}(\Delta)},
\end{equation}
\end{lem}
 
 This lemma was proved in \cite{IS4}, see Lemma 2.2 there, but the claim
there was that the estimate \eqqref(est-A1) holds also for $r=1$. This is not quite true. 
What was indeed proved there is this estimate for any $0<r<1$ as we state this here. At 
the same time due
to Lemma \ref{van-order} we have that $u$ is in $L^{2,p}_{loc}(\Delta)$ for any $2<p<\infty$
and therefore, after a rescaling we can assume that $u\in L^{2,p}(\Delta)\subset 
\calc^{lip}(\Delta)$. In the annulus $\Delta\setminus \Delta (r)$ the "twisting"
terms do not spoil the Lipschitz continuity of $\zeta^{-\nu}\cdot A(u(\zeta))\cdot 
\zeta^{\lambda}$  and  therefore we obtain the estimate
\begin{equation}
\eqqno(est-A2)
Lip_{\Delta}\big(\zeta^{-\nu}\cdot A(u(\zeta))\cdot \zeta^{\lambda}\big)\le 
C(p,r)\cdot Lip_B(A)\cdot\norm{u}_{\calc^{lip}(\Delta)},
\end{equation}
which holds true under the conditions of Lemma \ref{lip-cont} with an additional assumption
that $u$ is defined and $J$-holomorphic in a neighbourhood of $\bar\Delta$.
 
\begin{lem}
\label{disk-pert-r}
Let $J$ be a Lipschitz-continuous almost complex structure on $\rr^{2n}$ such that
$J|_{\rr^n} = J\st$ and let $u_0:(\Delta^+, \db_0\Delta^+,0)\to (\rr^{2n} , \rr^n , 0)$ 
be a $J$-holomorphic map. Let $\nu \ge 0$ be an integer and $\wect_0\in\rr^n$ be a vector. 
Then there exists a $L^{1,p}$-regular map $w :(\Delta^+, \db_0\Delta^+,0)\to (\rr^{2n}, 
\rr^n,0)$ with  $w(0)=\wect_0$ such that 
\begin{equation}
\eqqno(pert2)
u(\zeta) = u_0(\zeta) + \zeta^{\nu}\cdot w(\zeta)
\end{equation}
is $J$-holomorphic. Moreover, for any $2<p<\infty$ one has the following estimate 
\begin{equation}
\eqqno(l1p-est)
\norm{w}_{L^{1,p}(\Delta^+)}\le C_p|\wect_0|,
\end{equation}
with a constant $C_p$ independent of $\wect_0$.
\end{lem}
\proof As above,  we equip the trivial bundle $E^+=\Delta^+\times \rr^{2n}$ with 
the complex structure $J_{u_0}=J\circ u_0$, which is then extended to the complex 
structure $\tilde J_{u_0}$ on $E=\Delta\times \rr^{2n}$ by reflection, as in 
\eqqref(ext-ref-j). We extend $u_0$ (viewed as a section of $E^+$) to the disk 
$\Delta$ by reflection, denoting this extension as $\tilde u_0$. Recall that this 
extension $\tilde u_0$ is a $\tilde J_{u_0}$-holomorphic section of $E$. We look for a 
perturbation of $\tilde u_0$ in the form 
\begin{equation}
\eqqno(pert3)
\tilde u(\zeta) = \tilde u_0(\zeta) + \zeta^{\nu}\cdot \tilde w(\zeta),
\end{equation}
where  $\tilde w$ satisfies the reality condition  
\begin{equation}
\eqqno(real-cond-w)
\overline{\tilde w(\bar\zeta)} = \tilde w(\zeta).
\end{equation}  
The latter is equivalent the reality condition  
$\overline{\tilde u(\bar\zeta)}= \tilde u(\zeta)$ on $\tilde u$. Setting
$w\deff\tilde w|_{\Delta^+}$ we consider $w$ as a mapping $\Delta^+\to (\rr^{2n},J)$
and take $J_{u_0 + \zeta^{\nu}w}\deff J\circ (u_0 + \zeta^{\nu}w)$ as a new
complex structure on $E^+$. We extend this structure to $E$ by reflection 
and denote this extended structure as $\tilde J_u$. Our principal task is 
to find $\tilde u$ (or, equivalently to find $\tilde w$) as in \eqqref(pert3) such 
that it is $\tilde J_u$-holomorphic and $\tilde w$ satisfies the same type 
estimate as \eqqref(l1p-est), \ie $\norm{\tilde w}_{L^{2,p}(\Delta)} \le C_p\, 
|\wect_0|$. After that we restrict $\tilde w$ to $\Delta^+$ and obtain the 
required $\tilde w$.

\smallskip Now write the equation of holomorphicity for $\tilde u$:

\[
\begin{split}
 0 & = \dbar_{\tilde J_u}\tilde u = \db_{\xi}\tilde u + \tilde J_u\db_{\eta}\tilde u = 
\db_{\xi}\tilde u + \big(\tilde J_u - \tilde J_{u_0}\big)\db_{\eta}\tilde u + 
\tilde J_{u_0}\db_{\eta}\tilde u =
\\[2pt]
& = \db_{\xi}\tilde u_0 + \tilde J_{u_0}\db_{\eta}\tilde u_0 + \db_{\xi}(\zeta^{\nu}
\tilde w) + \tilde J_{u_0}\db_{\eta}(\zeta^{\nu}\tilde w) + \big(\tilde J_u - \tilde J_{u_0}
\big)\db_{\eta}(\tilde u_0 + \zeta^{\nu}\tilde w) =
\\[2pt]
& = \dbar_{\tilde J_{u_0}}(\zeta^{\nu}\tilde w) + \big(\tilde J_u - 
\tilde J_{u_0}\big)\big(\db_{\eta}\tilde u_0 + \db_{\eta}(\zeta^{\nu} \tilde w)\big).
\end{split}
\]
Here we used the $\tilde J_{u_0}$-holomorphicity of $\tilde u_0$, \ie that 
$\dbar_{\tilde J_{u_0}}\tilde u_0 =0$. Our goal is to solve the equation
\begin{equation}
\eqqno(pert3-r)
\dbar_{\tilde J_{u_0}}(\zeta^{\nu}\tilde w) =
\big(\tilde J_{u_0} - \tilde J_u\big)\big(\db_{\eta}\tilde u_0 + 
\db_{\eta}(\zeta^{\nu}\tilde w)\big), 
\end{equation}
with $\tilde w(\zeta)$ satisfying the reality condition \eqqref(real-cond-w). 
Multiplying this equation by $\zeta^{-\nu}$ we obtain
\begin{equation}
\eqqno(pert4-r)
\zeta^{-\nu}\dbar_{\tilde J_{u_0}}(\zeta^{\nu}\tilde w) = \zeta^{-\nu} 
\big(\tilde J_{u_0} - \tilde J_u\big)\big( \db_{\eta}\tilde u_0 + \db_{\eta}(\zeta^{\nu}
\tilde w)\big).
\end{equation}

\medskip\noindent{\slsf Step 1. Estimation of the left hand side of \eqqref(pert4-r).}
The left hand side of \eqqref(pert4-r) can be transformed as follows
\[
\zeta^{-\nu}\dbar_{\tilde J_{u_0}}(\zeta^{\nu}\tilde w) = \big(\db_{\xi} + \zeta^{-\nu}\tilde J_{u_0}
\zeta^{\nu}\db_{\eta}\big)\tilde w + \nu \zeta^{-\nu}\big(\zeta^{\nu -1} + \tilde J_{u_0}J\st \zeta^{\nu -1}
\big) \tilde w =
\]
\[
= \big(\db_{\xi} + \zeta^{-\nu}\tilde J_{u_0}\zeta^{\nu}\db_{\eta}\big)\tilde w
+ \nu \zeta^{-\nu}\big(1 + \tilde J_{u_0}J\st \big)\zeta^{\nu -1}\tilde w =: 
\dbar_{\tilde J^{(\nu )}_{u_0}}\tilde w + R^{(\nu )}\tilde w .
\]
Here 
\[
\dbar_{J^{(\nu )}_{u_0}} \deff \db_{\xi} + \tilde J^{(\nu )}_{u_0}\db_{\eta}
\]
is the $\dbar$-operator with respect to the ``twisted'' complex structure
$\tilde J^{(\nu )}_{u_0}\deff \zeta^{-\nu}\tilde J_{u_0}\zeta^{\nu}$ on $E$, and 
\begin{equation}
\eqqno(term-r)
R^{(\nu)}\deff \nu \zeta^{-\nu}\big(1 + \tilde J_{u_0}J\st \big)\zeta^{\nu -1}
\end{equation}
is the zero order term. The structure $\tilde J^{(\nu )}_{u_0}$ is Lipschitz-continuous 
due to Lemma \ref{lip-cont}. To apply the quoted lemma in our settings write $\zeta^{-\nu}
\tilde J_{u_0}\zeta^{\nu} = \zeta^{-\nu}\left(\tilde J_{u_0}-J\st + J\st\right)\zeta^{\nu} 
= \zeta^{-\nu}\left(\tilde J_{u_0}-J\st\right)\zeta^{\nu} + J\st$ and observe that the 
operator $A(z)\deff \tilde J_{u_0}(z)-J\st $ vanishes at zero. Therefore Lemma 
\ref{lip-cont} applies to $\tilde u_0$ (regarded as a map to $\rr^{2n}$). Remark also
that $\norm{\tilde J^{(\nu)}-J\st }_{\calc^{lip}(\Delta)}$ is a small as we wish.
Operator 
\begin{equation}
\eqqno(r-detail)
R^{(\nu)} = \nu \zeta^{-\nu}\big(1 + \tilde J_{u_0}J\st \big)\zeta^{\nu -1} = 
\nu \zeta^{-\nu}\big(\tilde J_{u_0} - J\st \big)\zeta^{\nu -1}
\end{equation}
satisfies the obvious pointwise estimate 
\begin{equation}
\eqqno(est1-r)
|R^{(\nu)}(\zeta)|\leq \nu Lip(J)\norm{\tilde u_0}_{\calc^{lip} (\Delta)},
\end{equation}
which in its turn implies the estimate
\begin{equation}
\eqqno(est2-r)
\norm{R^{(\nu)}}_{L^p(\Delta)} \leq \nu Lip(J)\norm{\tilde u_0}_{\calc^{lip}
(\Delta)}.
\end{equation}
Moreover, it is real in the sense that for $\tilde w$ satisfying 
$\overline{\tilde w (\bar\zeta)} = \tilde w(\zeta)$ one has: 
\[
\overline{R^{(\nu)}(\bar\zeta)\tilde w (\bar\zeta) } = 
\overline{\nu\bar\zeta ^{-\nu}\left(\adyn + \tilde J_{u_0}(\bar\zeta)J\st\right) \bar\zeta ^{\nu -1}\tilde w 
(\bar\zeta)} = \overline{\nu\bar\zeta ^{-\nu}\left[\bar\zeta ^{\nu -1}   + \tilde J_{u_0}(\bar\zeta)
J\st\bar\zeta^{\nu -1} \right]\tilde w (\bar\zeta) } = 
\]
\[
= \nu\zeta ^{-\nu}\left[\zeta ^{\nu -1}   + \tilde J_{u_0}(\zeta)J\st\zeta^{\nu -1}\right]
\overline{\tilde w (\bar\zeta)} = \nu\zeta^{-\nu}\left(\adyn  + \tilde J_{u_0}(\zeta)J\st\right)
\zeta^{\nu -1} \tilde w (\zeta)  = R^{(\nu)}(\zeta)\tilde w (\zeta),
\]
since $\bar J\st = -J\st$. In addition, operator $\dbar_{\tilde J^{(\nu )}_{u_0}}$  is real 
as well. Indeed, making the change $\zeta = \bar z$, where $z=x+iy$, we see that 
\[
\overline{\dbar_{J^{(\nu )}_{u_0}}(\zeta)\tilde w(\zeta)} = \overline{\db_{\xi} \tilde w(\zeta)}
+ \overline{\zeta^{-\nu}\tilde J_{u_0}(\zeta)\zeta^{\nu}\db_{\eta}\tilde w(\zeta)} = 
\db_x\overline{\tilde w(\bar z)} + 
\]
\[
 + z^{-\nu}\tilde J_{u_0}(z)z^{\nu}\overline{\db_y\tilde w(\bar z)} = \db_x\tilde w(z) + 
 z^{-\nu}\tilde J_{u_0}(z)z^{\nu}\db_y\tilde w(z) = \dbar_{J^{(\nu )}_{u_0}}(z)\tilde w(z), 
\]
\ie 
\[
\overline{\dbar_{J^{(\nu )}_{u_0}}(\bar\zeta)\tilde w(\bar\zeta)} =  
\dbar_{\tilde J^{(\nu )}_{u_0}}(\zeta)\tilde w(\zeta),
\]
as stated. Therefore the left hand side of \eqqref(pert4-r) has the form 
\begin{equation}
D_{\tilde J^{(\nu)}_{u_0},R^{(\nu)}}(\tilde w) \deff \dbar_{\tilde J^{(\nu )}_{u_0}}\tilde w + R^{(\nu )}
\tilde w,
\end{equation}
for the Lipschitz-continuous $\tilde J^{(\nu )}_{u_0}$ and real-valued $R^{(\nu)}$ with 
$\norm{R^{(\nu)}}_{L^p(\Delta)}$ small. Notice that $D_{\tilde J^{(\nu)}_{u_0},R^{(\nu)}}$
is also real.

\smallskip\noindent{\slsf Step 2.  Estimation of the right hand side of 
\eqqref(pert4-r).} {\it The right hand side
\begin{equation}
\eqqno(term-f)
F^{(\nu)}(\zeta,\tilde w) \deff \zeta^{-\nu} \big(\tilde J_{u_0} - \tilde J_{u_0+\zeta^{\nu}w}
\big)\big( \db_{\eta}\tilde u_0 + \db_{\eta}(\zeta^{\nu}\tilde w)\big)
\end{equation}
of \eqqref(pert4-r) admits the following estimates:

\begin{equation}
\eqqno(est1-f)
\norm{F^{(\nu)}(\zeta , \tilde w)}_{L^p(\Delta)}\leq C\cdot Lip(J) \norm{\tilde w}_{L^{1,p}(\Delta)},
\end{equation}
and 
\begin{equation}
\eqqno(est3-r)
\norm{F^{(\nu)}(\zeta , \tilde w_1) - F^{(\nu)}( \zeta , \tilde w_2) }_{L^p(\Delta)}
\leq C\cdot Lip(J)\norm{\tilde w_1 - \tilde w_2}_{L^{1,p}(\Delta},
\end{equation}
provided an priori bound on $\norm{\tilde w}_{L^{\infty}(\Delta)}\le \frac{1}{2}$ is imposed, 
see Remark \ref{induct-bnd}.} Notice that the substitution $\tilde w_2(\zeta) \equiv0$ 
in the second estimate gives the first one. Therefore we need to prove only the second estimate.
Set $\tilde u_k(\zeta) \deff \tilde u_0(\zeta) + \zeta^{\nu}\tilde w_k(\zeta)$ and write 
\begin{equation}
\eqqno(f-f1)
\begin{split}
& F^{(\nu)}(\zeta , \tilde w_1) - F^{(\nu)}( \zeta , \tilde w_2) = 
\zeta^{-\nu} \big(\tilde J_{u_0} - \tilde J_{u_1}\big)\big( \db_{\eta}
\tilde u_0 + \db_{\eta}(\zeta^{\nu}\tilde w_1)\big) - 
\\[2pt]
& - \zeta^{-\nu} \big(\tilde J_{u_0} - \tilde J_{u_2}\big)\big( \db_{\eta}
\tilde u_0 + \db_{\eta}(\zeta^{\nu}\tilde w_2)\big) = 
\zeta^{-\nu} \big(\tilde J_{u_2} - \tilde J_{u_1}\big)
\big( \db_{\eta}\tilde u_0 + \db_{\eta}(\zeta^{\nu}\tilde w_1)\big)
\\[2pt]
& + \zeta^{-\nu} \big(\tilde J_{u_0} - \tilde J_{u_2}\big)\big( \db_{\eta}(\zeta^{\nu}
\tilde w_1) - \db_{\eta}(\zeta^{\nu}\tilde w_2)\big).
\end{split}
\end{equation}
Now remark that $|\tilde J_{u_2} - \tilde J_{u_1}| \le Lip(J)|\zeta|^{\nu}|\tilde w_2 - 
\tilde w_1|$ as well as $|\tilde J_{u_2} - \tilde J_{u_0}| \le Lip(J)|\zeta|^{\nu}$,
both provided $\norm{\tilde w}_{L^{\infty}(\Delta)}\le \frac{1}{2}$. 
Estimate \eqqref(est3-r) follows.

\smallskip Furthermore, notice again that operator $F^{(\nu)}(\zeta , \cdot): 
L^{1,p}(\Delta) \to L^p(\Delta)$ is real. The proof is analogous to that of the reality
of $R^{[\nu)}$:
\[
\overline{F^{(\nu)}(\bar\zeta , \tilde w(\bar \zeta))} = \overline{\bar\zeta^{-\nu} 
\big(\tilde J_{u_0}(\bar\zeta) - \tilde J_u(\bar\zeta)\big)\big( \db_{\eta}
\tilde u_0(\bar\zeta)+ \db_{\eta}(\bar\zeta^{\nu}\tilde w(\bar\zeta)\big)} = 
\]
\[
= \zeta^{-\nu} \big(\tilde J_{u_0}(\zeta) - \tilde J_u(\zeta) \big)\big( \db_{\eta}
\overline{\tilde u_0 (\bar\zeta)} + \db_{\eta}(\zeta^{\nu}\overline{\tilde 
w(\bar\zeta)})\big) = 
\]
\[
= \zeta^{-\nu} \big(\tilde J_{u_0} - \tilde J_u\big)\big( \db_{\eta}
\tilde u_0(\zeta) + \db_{\eta}(\zeta^{\nu}\tilde w (\zeta))\big) = F^{(\nu)}(\zeta , 
\tilde w(\zeta)).
\]

\smallskip  Therefore our problem is reduced to solving the following equation
\begin{equation}
\eqqno(syst1-r)
\begin{cases}
D_{J^{(\nu)}_{u_0}}\tilde w = F^{(\nu)}(\zeta , \tilde w),\cr
\tilde w(0) = \wect_0,\\[3pt]
\overline{w(\bar\zeta)} = \tilde w(\zeta).
\end{cases}
\end{equation}

We follow the proof of Theorem 6.1 in \cite{IS4}, insuring the reality condition. 
Let's start with the following remark.  Let $J$ be a continuous complex
structure on the trivial bundle $E=\Delta\times \rr^{2n}$, standard over $\beta_0$.
Consider the following operators $L^{1,p}(\Delta, \rr^{2n}) \to L^p
(\Delta , \rr^{2n})$ for $p>1$:
\begin{equation}
\dbar_Jw \deff \frac{\db w}{\db x} + J\frac{\db w}{\db y} \qquad\text{ and } \qquad
D_Jw\deff \dbar_Jw + Rw,
\end{equation}
where $R$ is a matrix valued function from $L^p(\Delta)$.
\begin{prop}
\label{morrey1}
If the  norms $\norm{J-J\st}_{\calc^0(\Delta)}, \norm{R}_{L^p(\Delta)}$ 
are sufficiently small then:

\sli there exists a bounded linear operator $T_{J,R}^0:L^p(\Delta)\to L^{1,p}(\Delta)$ such that $(T_{J,R}^0w)(0) = 0$ 

\quad for every $w\in L^{1,p}(\Delta)$ and which is right inverse to  $D_J$, \ie is such that 

\quad $(\dbar_J + R)\circ T_{J,R}^0\equiv \id$;

\slii if , moreover, $J(\zeta)=J\st$ and $R(\zeta)$ is real for real $\zeta$, 
then $T_{J,R}^0$ is  real provided $J$

\quad satisfies 
\begin{equation}
\eqqno(im-j)
\overline{J(\bar\zeta)\left[\bar \vect\right]} = - J(\zeta) \left[\vect\right].
\end{equation}
\end{prop}
\proof \sli For $J=J\st$ and\/ $R=0$ the operator in question is constructed using 
the standard Cauchy-Green operator $ T_{CG}$. Namely, we set $T_{J\st,0}^{\,0}(u) = 
T_{CG}u - (T_{CG}u)(0)$, where 
\begin{equation}
\eqqno(cg-op)
\left(T_{CG}w\right)(z) = \frac{1}{2\pi i } \int\limits_{\Delta}\frac{w(\zeta)}{\zeta - z}
d\zeta\wedge d\bar\zeta .
\end{equation}
This gives a  bounded operator  $T_{J\st,0}^{\,0}:L^p(\Delta , \rr^{2n})\to L^{1,p}
(\Delta ,\rr^{2n})$ with the desired properties. For general $J$ and $R$ the operator 
$T_{J,R}^{\,0}$ can be constructed as a perturbation series:
\[
T_{J,R}^{\,0} \deff \left(\dbar_J+R\right)^{-1} = \left(\dbar_{J\st} + \dbar_J 
- \dbar_{J\st} + R\right)^{-1} = \dbar_{J\st}^{-1}\left(\id + (\dbar_J 
- \dbar_{J\st} + R)\circ \dbar_{J\st}^{-1}\right)^{-1} = 
\]
\begin{equation}
\eqqno(pert-op)
= T_{J\st,0}^{\,0}\circ\sum_{n=0}^{\infty} (-1)^n\big((\dbar_J - \dbar_{J\st} 
+ R)\circ T_{J\st,0}^{\,0}\big)^n.
\end{equation}
The series converges in the appropriate operator norm provided the norm of the 
perturbation $\dbar_J - \dbar_{J\st} + R$ is sufficiently small. The latter is  
estimated by $\norm{J-J\st}_{\calc^0(\Delta)} + \norm{R}_{L^p(\Delta)}$.

\smallskip\noindent\slii Notice that $T_{J\st ,0}$ is real. The reality of $T_{J,R}^0$
follows from the reality of $R$ and the condition \eqqref(im-j) on $J$ as before.

\smallskip\qed


\smallskip We turn back to the proof of Lemma \ref{disk-pert-r} and are going to apply 
the Newton's method of successive approximations using the estimates and Proposition 
\ref{morrey1} as above for $J=\tilde J^{(\nu)}_{u_0}$ and $R = R^{(\nu)}$. Set 
\[
\tilde w_1(\zeta) = \wect_0 - T_{\tilde J^{(\nu)}_{u_0},R^{(\nu)}}^0\big(D_{\tilde J^{(\nu)}_{u_0},R^{(\nu)}}\wect_0\big).
\]

Notice that $\tilde w_1$ is a solution of the following system
\begin{equation}
\eqqno(syst2)
\begin{cases}
D_{J^{(\nu)}_{u_0},R^{(\nu)}}\tilde w_1 = 0, \cr
\tilde w_1(0) = \wect_0, \\[3pt]
\overline{\tilde w_1(\bar \zeta)} = \tilde w_1(\zeta).
\end{cases}
\end{equation}

The reality of $\tilde w_1$ follows from that one of $D_{\tilde J^{(\nu)}_{u_0},R^{(\nu)}}$ and  
$T_{\tilde J^{(\nu)},R^{(\nu)}}^0$. Furthermore, set
\begin{equation}
\eqqno(newton)
\tilde w_{n+1} = T_{J^{(\nu)}_{u_0},R^{(\nu)}}^0\big[F^{(\nu)}(\zeta, \tilde w_n)\big] + \tilde w_1,
\end{equation}
where 
\begin{equation}
\eqqno(ef-nu-n)
F^{(\nu)}(\zeta , \tilde w_n) \deff \zeta^{-\nu}\big(\tilde J_{u_0} - \tilde J_{u_0+\zeta^{-\nu}
w_n}\big) \zeta^{\nu}\big( \db_{\eta}u_0 + \db_{\eta}(\zeta^{\nu}\tilde w_n)\big),
\end{equation}
and $w_n \deff \tilde w_n|_{\Delta^+}$. Notice that all $\tilde w_n$ stay real.

\bigskip 
Estimates \eqqref(est1-f) and \eqqref(est3-r) ensure the convergence of the
iteration process. Indeed, as it was explained there at the beginning of the proof 
using dilatations we can suppose that $Lip(J)$ as well as 
$\norm{u_0}_{\calc^{lip}(\Delta)}$  are as small as we wish, \ie less than
some $\eps>0$ to be specified later in the process of the proof. This implies also 
that $\norm{R^{(\nu)}}_{L^p(\Delta)}\leq \eps$. Finally,  $|\wect_0|$ 
will be supposed also small enough. 

\smallskip\noindent{\slsf Step 3.  $L^{1,p}$-bound on $\tilde w_n$.} {\it One
has the following estimates:
\begin{equation}
\eqqno(newton-r1)
\norm{\tilde w_{n+1}}_{L^{1,p}(\Delta)} \le C_p\cdot |\wect_0|,
\end{equation}
and 
\begin{equation}
\eqqno(newton-r2)
\norm{\tilde w_{n+1} -\tilde w_n}_{L^{1,p}(\Delta)} \le C_p\cdot Lip_B(J)\cdot
\norm{\tilde w_{n+1} -\tilde w_n}_{L^{1,p}(\Delta)}.
\end{equation}}

\medskip We start from $\tilde w_1$:
\[
\begin{split}
&\norm{\tilde w_1}_{L^{1,p}(\Delta)} \le a_p |\wect_0| + \norm{T_{\tilde J^{(\nu)}_{u_0},
R^{(\nu)}}^0\big(D_{\tilde J^{(\nu)}_{u_0},R^{(\nu)}}\wect_0\big)}_{L^{1,p}(\Delta)} 
\le 
\\[2pt]
& \le a_p |\wect_0|+ C\, \norm{R^{(\nu)}}_{L^p(\Delta)}|\wect_0| \le C_p\cdot |\wect_0|,
\end{split}
\]
by \eqqref(est2-r). Here $a_p=\pi^{1/p}=\norm{1}_{L^p(\Delta)}$. Furthermore, due to 
\eqqref(est1-f), we get
\begin{equation}
\eqqno(case-l=0)
\norm{\tilde w_{n+1}}_{L^{1,p}(\Delta)} \le C|\wect_0| + C \norm{F^{(\nu)}(\zeta, 
\tilde w_n)}_{L^p(\Delta)} \le C|\wect_0| +C\cdot Lip(J) \norm{\tilde w_n}_{L^{1,p}(\Delta)}.
\end{equation}
Now one proves by induction that 
\[
\norm{\tilde w_{n+1}}_{L^{1,p}(\Delta)} \le C|\wect_0|\sum_{k=0}^l(C\eps)^k + (C\eps)^{l+1}
\norm{\tilde w_{n-l}}_{L^{1,p}(\Delta)}. 
\]
with the case $l=0$ being the inequality \eqqref(case-l=0) above, where $\eps = Lip(J)$. 
This implies 
\[
\norm{\tilde w_{n+1}}_{L^{1,p}(\Delta)} \le \frac{C|\wect_0|}{1-C\eps} + (C\eps)^n
\norm{\tilde w_1}_{L^{1,p}(\Delta)} \le C_p\cdot |\wect_0|,
\]
as stated.

\smallskip As for the second estimate write
\[
\begin{split}
& \norm{\tilde w_{n+1} -\tilde w_n}_{L^{1,p}(\Delta)} = \norm{T_{J^{(\nu)}_{u_0},R^{(\nu)}}
^0\big[F^{(\nu)}(\zeta, \tilde w_n)-F^{(\nu)}(\zeta, \tilde w_{n-1})\big]}_{L^{1,p}(\Delta)}\le
\\[2pt]
&  \le \norm{F^{(\nu)}(\zeta, \tilde w_n)-F^{(\nu)}(\zeta, \tilde w_{n-1})}_{L^p(\Delta)}
\le C\cdot Lip(J)\norm{\tilde w_1 - \tilde w_2}_{L^{1,p}(\Delta},
\end{split}
\]
by \eqqref(est3-r), and this gives us \eqqref(newton-r2).

\begin{rema} \rm
\label{induct-bnd}
Notice that we proved inductively that $\norm{\tilde w_n}_{L^{\infty}(\Delta)}\leq 
\frac{1}{2}$ in order to use estimates \eqqref(est1-f) and \eqqref(est3-r). 
Indeed, fix some $p_0>2$ and get from here  our inductive assumption that
$\norm{\tilde w_n}_{L^{1,p_0}(\Delta)}\leq \frac{1}{2H_{p_0,\alpha_0}}$ in order to keep
$\norm{\tilde w_n}_{L^{\infty}} \leq \frac{1}{2}$. Since the constant $C$ in \eqqref(newton-r1)
depends on $p$ we fix $p_0$ and then choose $|\wect_0|$ to be appropriately 
small. 
\end{rema}

\medskip\noindent{\slsf Step 4.} {\sl Convergence of the approximation.}
Write
\[
\norm{\tilde w_{n+1}-\tilde w_n}_{L^{1, p}(\Delta)}\leq C
\norm{F^{(\nu)}(\zeta,\tilde w_n) - F^{(\nu)}(\zeta,\tilde w_{n-1})}_{L^p(\Delta)}\leq
\]
\begin{equation}
\eqqno(6.20)
 \leq  C\cdot\eps \norm{\tilde w_{n}-\tilde w_{n-1}}_{L^{1,p}(\Delta)}
\end{equation}
by \eqqref(est3-r) with $\eps >0$ as small as we wish. 
This gives us the desired convergence of 
\[
\tilde w_n = \sum_{k=1}^n\left(\tilde w_k - \tilde w_{k-1}\right ) +\wect_0
\]
in $L^{1,p}(\Delta)$  to a limit $\tilde w$, which obeys the relation
\begin{equation}
\eqqno(newton-fin)
\tilde w = T_{J^{(\nu)}_{u_0},R^{(\nu)}}^0\big[F^{(\nu)}(\zeta,\tilde w)\big] + \tilde w_1,
\end{equation}
and therefore is a solution $\tilde w$ of \eqqref(syst1-r). From \eqqref(newton-r1) 
we get also \eqqref(l1p-est). 

\smallskip 
The operator for the solution of the system \eqqref(syst1-r) constructed in 
the proof is continuous  with respect to the $L^{1,p}$-norm. Indeed, 
let $\tilde w'$ and $\tilde w''$ be solutions 
of \eqqref(syst1-r) with initial data $\tilde w'(0)=\wect'_0$ 
and $\tilde w''(0)=\wect''_0$. Then from \eqqref(newton) we have
\[
\norm{\tilde w' - \tilde w''}_{L^{1,p}(\Delta)}\leq |\wect'_0-\wect''_0| + 
C\norm{F^{(\nu)}(\zeta, \tilde w') - F^{(\nu)}(\zeta,\tilde w'')}_{L^p(\Delta)}\leq
\]
\[
\leq  |\wect'_0-\wect''_0| + C\,\eps \norm{\tilde w' -\tilde w''}_{L^{1,p}(\Delta)}.
\]
And therefore
\begin{equation}
\label{cont-dep1}
\norm{\tilde w' - \tilde w''}_{L^{1,p}(\Delta)}\leq \frac{1}{1-C\eps} |\wect'_0-\wect''_0|,
\end{equation}
\ie solution $w$ of \eqqref(syst1-r) continuously depend on the initial data 
$\tilde w(0)=\wect_0$. In particular we have the uniqueness of the solution for 
the given initial data $\wect_0$ and also 
\begin{equation}
\label{cont-dep2}
\norm{\tilde w}_{L^{1,p}(\Delta)}\leq \frac{1}{1-C\eps} |\wect_0|.
 \end{equation}

Lemma is proved.

\smallskip\qed

\newprg[B-POS.proof]{Proof of the Positivity of Boundary Intersections.} 

Let us prove Theorem \ref{b-pos-thm}.

\smallskip\noindent (\sli Our complex half-discs $C_k = u_k(\Delta^+)$ are tangent to 
each other at $p_0$, \ie $\mu_1=\mu_2=1$ and $\vect_0^1$, $\vect^2_0$ are collinear. 
After shrinking, if necessary, we can assume that $u_1(\Delta^+)$ and $u_2(\Delta^+)$ 
are embedded and after a rectification we can assume that $W=\rr^2$ and $p_0=0$. 
Rescaling the parametrizations of $u_k$ and making a real rotation of $\rr^2
\subset \cc^2$, \ie applying a complex linear map which sends $\rr^2$ to $\rr^2$,
we can assume that $\frac{\db u_1}{\db \zeta}(0) = e_1 = \pm \frac{\db u_2}{\db \zeta}(0)$.
Therefore we have two possibilities. 

\smallskip\noindent{\slsf Case 1.} {\it $C_k$-s touch each other at zero, \ie 
$\vect_0^1=\vect^2_0 = e_1$.} Extend $u_k$ by reflection and denote these extensions 
as $\tilde u_k$. Applying Lemma  \ref{comp-lem-r} we can write
\begin{equation}
\eqqno(comp-eq+)
\tilde u_2(\zeta) -\tilde u_1(\psi (\zeta)) = \zeta^{\nu}\tilde w(\zeta),
\end{equation}
where $\tilde w(0) = e_2$ and $\psi(\zeta) = \zeta + O(\zeta^2)$ is a holomorphic 
reparametrization, real for real $\zeta$, and $\nu >1$. 

\smallskip Consider $\tilde u_k$ as mappings to $\cc^2$. They are not $J$-holomorphic but do 
satisfy the comparison relation \eqqref(comp-eq+). Denote by $\sss^3_r$ the sphere of radius 
$r$ centred at the origin. For $r>0$ small enough circles $\gamma_i(r) \deff \tilde u_k(\Delta)
\cap \sss^3_r$ are immersed. This follows from the form \eqqref(nrm-frm3) of the differentials
of $\tilde u_k$, \ie from the fact that 
\[
d\ti u_k(\zeta) = e_1 + O(|\zeta|^\alpha).
\]
Equation \eqqref(comp-eq+) shows that the curve $\gamma_2(r)$ stays in the tubular $\rho =
2r^{\nu}$-neighbourhood of $\gamma_1(r)$ and winds $\nu$ times around $\gamma_1(r)$. This shows
that the linking number $l(\gamma_1(r),\gamma_2(r))$ is $\nu >1$, \ie $\bindo (u_1,u_2) 
\deff \ind_0(\ti u_1,\ti u_2) = \nu > 1$.

\smallskip\noindent{\slsf Case 2.} {$C_k$-s meet each other at zero, \ie $\vect_0^1= - 
\vect^2_0 = e_1$.} This time by Lemma \ref{comp-lem-op}  we have 
\begin{equation}
\eqqno(comp-eq-)
\tilde u_2(-\zeta) -\tilde u_1(\psi (\zeta)) = \zeta^{\nu}\tilde w(\zeta),
\end{equation}
Disk $\tilde u_2(\Delta)$ coincides with $\tilde u_2(-\Delta)$ and therefore the same arguments
which were used in the first case apply. Item (\sli  of Theorem \ref{b-pos-thm}
is proved.

\begin{rema} \rm
\label{pnt-2}
Notice that it is the intersection of $\tilde u_k$-s at zero which was proved to be 
positive. Let us see that $\tilde u_k$-s intersect only at zero, provided one is not 
a reparametrization  of another, \ie if $w(0)\not=0$. Indeed, 
we can assume that $\vect_0 =e_1$ and $\wect_0\deff w(0) = e_2$. Denote by $z_1,z_2$
the standard complex coordinates in $\cc^2$. Mapping $\ti u_1(\psi(\zeta))$  
rename back to $\ti u_1(\zeta)$. Since the first coordinate of $\ti u_1$ satisfies 
$z_1 = \ti u_1^1(\zeta ) = \zeta + O(|\zeta|^{1+\alpha})$ we see that $\zeta = 
z_1 + O(|z_1|^{1+\alpha})$. Consider first the touching case and write \eqqref(comp-eq+) 
for the second coordinates of $u_1$ and $u_2$: 
\[
\begin{split}
&\ti u_2^2(z_1 + O(|z_1|^{1+\alpha})) - \ti u_1^2(z_1 + O(|z_1|^{1+\alpha})) = 
\\[2pt]
& = (z_1 + O(|z_1|^{1+\alpha}))^{\nu} + O(z_1 + O(|z_1|^{1+\alpha})^{\nu +\alpha})=
z_1^{\nu} + O(|z_1|^{\nu + \alpha}).
\end{split}
\]
This means that as the functions of the first coordinate the second coordinates of
$\ti u_k$ are different by a term which vanishes only at zero. 

\smallskip In the meeting case repeat the same argument with  $\ti u_2(-\zeta)$ on the
pace of $\ti u_2(\zeta)$ and prove that $\ti u_1(\Delta)\cap \ti u_2(-\Delta) = \{0\}$.
But, again, $\ti u_2(-\Delta) = \ti u_2(\Delta)$. The assertion is proved.
\end{rema}

\smallskip\noindent (\slii We proved that if $C_1$ and $C_2$ are tangent (both touching 
or meeting) then the intersection index is $\nu \ge 2$. This implies the statement of
the item (\slii of the theorem.

\smallskip\noindent (\sliii  Along the proof of the part (\sli of our theorem we remarked 
that point $p_0=0$ is an isolated point of intersection of $C_1=\tilde u_1(\Delta)$ and 
$C_2= \tilde u_2(\Delta)$. Shrinking $\Delta$ (or rescaling) we can suppose that $p_0$ is the
only intersection point of $C_1$ and $C_2$, see Remark \ref{pnt-2}.  After a rectification 
in an appropriate coordinates we have $\vect_0^1 = e_1$ and $\vect_0^2=\pm e_1$. 
By the comparison relation(s) we can write
\[
\tilde u_2(\pm\zeta ) = \tilde u_1(\zeta) + \zeta^{\nu} w(\zeta),
\]
where we denoted  $\tilde u_1(\psi(\zeta))$ back as $\tilde u_1(\zeta)$. Notice that 
$\nu \ge 2$ since the order of vanishing of $u_k$-s at zero are $=1$.
Moreover, $w(0) = e_2$ and therefore we have 
\[
\tilde u_2(\pm\zeta) = u_1(\zeta ) + \zeta^{\nu} e_2 + w_1(\zeta),
\]
with $w_1(\zeta) = O(|\zeta|^{\nu + \alpha})$. Using Lemma \ref{disk-pert-r} 
we find a $J$-holomorphic perturbation of $\tilde u_2(\pm \zeta)$ in the following form 
\[
u_2^{prt}(\zeta) = \tilde u_2(\pm\zeta) + \zeta^{\nu -1} \eps e_2 + \tilde w_2(\zeta),  
\]
where $\tilde w_2(\zeta) = O(|\zeta|^{\nu -1 + \alpha})$. Therefore 
\[
u_2^{prt}(\zeta ) - \tilde u_1(\zeta) = (\eps \zeta^{\nu -1} + \zeta^{\nu} )e_2 + 
\tilde w_3(\zeta),
\]
where $\tilde w_3(\zeta) = O(|\zeta|^{\nu -1 +\alpha})$. Notice that 
all these relations are valid for the extensions of our maps by reflection as 
they are written, and they all satisfy the reality condition. For $w_3$ this writes as 
$\ti w_3(\bar\zeta) = \overline{\ti w_3(\zeta)}$. Notice that 
$u_2^{prt}(\zeta ) - \tilde u_1(\zeta)$ has zero of order $\nu -1$ at zero. Moreover, 
being a small perturbation of $\tilde u_2(\pm \zeta)$, its index of intersection with 
$\ti u_1$ stays to be $\nu$. Let $\zeta_0$ be the (unique) root of $u_2^{prt}(\zeta )
- \tilde u_1(\zeta)$ different from zero. By the reality properties of our 
functions we see that $\bar\zeta_0$ should be also a root. Therefore 
$\zeta_0$ is real. Proceeding by induction we find a perturbation $u^{prt}_2$ 
of $\tilde u_2(\pm\zeta)$ such that $\ti u_2^{prt}$ has all intersections with 
$\ti u_1$ real and transverse. 

\smallskip\noindent\sliv Since perturbations $u^{prt}$ are small we have that the
linking number of their intersections with small sphere stays constant, \ie is 
$\nu$. In the case of transverse (positive!) intersections this number is equal to
the number of intersection points. Theorem \ref{b-pos-thm} is proved.

\smallskip\qed

\smallskip 
\begin{prop}
\label{point-int}
The index of intersection of two tangent $J$-holomorphic mappings  $u_k:(\Delta^+,
\db_0\Delta^+,0)\to (\rr^{2n}, \rr^n,0)$ such that one is not a reparametrization 
of another doesn't depend on the rectifying map.
\end{prop}
\proof Let $\Psi = \Psi_2\circ \Psi_1^{-1}$ be the change of 
the rectifying maps. Since $d\Psi (0)$ is $J\st$-linear it can be regarded as a 
complex $2\times 2$-matrix
in the standard basis of $\cc^2$. But $\Psi (\rr^n)\subset \rr^n$ and therefore this 
matrix has real coefficients. Denote by $u_k$ our half-disks after the first rectification
and by $v_k$ after the second, and by $\ti u_k$ and $\ti v_k$ their extensions by 
reflection. Since $v_k = \Psi(u_k)$  the reality of $d\Psi (0)$ implies that $\tilde v_k$ is 
a small perturbation of $\Psi (\tilde u_k)$. Therefore $\bindo (v_1,v_2) = \ind_0(\ti v_1,
\ti v_2) = \ind_0(\ti u_1,\ti u_2) = \bindo (u_1,u_2)$.

\smallskip\qed

\newsect[TOUCH]{Application I: Boundary Intersection of Two Analytic Disks in $\cc^2$}

\newprg[TOUCH.touch-r]{Index of intersection of two tangent real disks in $\rr^4$} 

By a smoothly embedded real half-disk in $\rr^4$ we mean a $\calc^{\infty}$-embedding of 
$u: \Delta^+\to\rr^4$ and denote $du(0)[\db_x]$ as $\vect_0$. $M \deff u(\Delta^+)$ 
is naturally oriented by the canonical orientation of $\Delta^+$. Let $M_1$, $M_2$ 
be two smoothly embedded half-disks in $\rr^4$ intersecting at $p_0 = u_1(0)=u_2(0)$.
We say that $M_k$ are \emph{tangent} at $p_0$ if 

\smallskip\sli their oriented tangent planes $T_{p_0}M_k$ coincide; 

\smallskip\slii tangent vectors $\vect_0^1 = du_1(0)[\db_x]$ and $\vect_0^2 = 
du_2(0)[\db_x]$ are collinear;

\smallskip\sliii we assume everywhere that the curves $\db_0M_2$ and $\db_0M_1$ 
are tangent to each other up 

\quad to a finite order, see Example \ref{inf-touch}. 

\begin{rema} \rm 
\label{inv-orient1}
The case when the tangents coincide but orientations are opposite can be treated analogously.
We shall not do that.
\end{rema}

After an appropriate reparametrisations of $u_k$-s we can assume that $\vect_0^2 = \pm 
\vect_0^1$. Applying an appropriate orientation preserving diffeomorphism of $\rr^4$ we 
can locally in a neighbourhood of $p_0$ bring $M_1$ to the form $M_1=\{(x_1,y_1,0,0) 
\in \rr^4: x_1^2 + y_1^2< 1, y_1\ge 0\}$ with $p_0=0$. Now $M_2$ will be realised in 
the form of a graph $M_2=\{(x_2,y_2) = \phi (x_1,y_1)\}$ over some domain $D\subset 
\rr^2_{x_1,y_1}$ which is diffeomorphic to a half-disc, and such that its boundary 
$\db D$ is tangent to $x_1$ axes at zero. Here $\rr^4 = \rr^2_{x_1,y_1}\oplus \rr^2_{x_2,y_2}$. Shrinking $M_1$ and $M_2$ we can assume that either $D=M_1=\Delta^+$, 
or $D=\Delta^-$. Set  $M_1^{\pm}\deff \pm M_1$. Mapping $\phi =(\phi_1,\phi_2): 
M_1^{\pm}\to \rr^2_{x_2,y_2}$ is 
supposed to be of class $\calc^{\infty}$ up to the edge of $M_1$, vanishing at zero.  
Set $\wect_0^k\deff du_k(0)[\db_y]$.
After an appropriate reparametrisation of $M_1$ we can assume that $\wect_0^1 =
\db_{y_1}$. Denote by $\pr_1(\wect^2_0)$ the orthogonal projection of $\wect_0^2$ 
to $\rr^2_{x_1,y_1}$. If $\pr_1(\wect^2_0)=0$ then $\wect_0^2\perp T_0M_1$ and this 
is not our case.

\begin{itemize}

\item If $\pr_1(\wect^2_0)$ is directed {\slsf  up} in $\rr^2_{x_1,y_1}$, 
\ie $\pr_1(\wect^2_0) = \alpha \vect_0^1 + \beta\wect_0^1$ with $\beta >0$,
we say that $M_2$ touches $M_1$ at $p_0$, in this case necessarily $\vect_0^2= 
\vect_0^1$;

\item in the opposite case we say that $M_2$ meets $M_1$ each other at $p_0$, 
in this case $\vect_0^2 = -\vect_0^1$. 
\end{itemize}

\smallskip These two cases need to be treated separately.

\smallskip\noindent{\slsf Case 1. $M_k$ touch each other at zero.}
In this case, after a shrinking, we can assume that $M_2$ is a graph 
over $M_1^+$. Our assumption that the curves $\db_0M_2$ and $\db_0M_1$ are 
tangent up to a finite order means that the Taylor series  of $\phi(x_1,0)$ 
at zero are not identically zero, \ie 
\begin{equation}
\eqqno(eq-m2)
\phi(x_1,0) = T(x_1) + O(|x_1|^{d+1}), 
\end{equation}
where $\deg T = d\in \nn$. Under  these conditions we can find one more diffeomorphism which
will bring our configuration to the following:

\begin{equation}
\eqqno(tot-r1)
\text{\sli } M_1 \text{ is still of the form } \{(x_1,y_1)\in \rr^2_{x_1,y_1}: x_1^2+y_1^2< 1, 
y_1\ge 0\};
\end{equation}
\[
\text{\slii } \phi_2(x_1,0) = 0, \text{ that is } M_2 \text{ is attached to } \rr_{x_1}\times 
\rr_{x_2} \text{ by its edge}.
\]

This can be done as follows. Write 
\begin{equation}
\eqqno(eq-phi)
\phi (x_1, 0) = (r_1x_1^d, r_2x_1^d) + O(|x_1|^{d+1}), \quad\text{where } r_1^2+ r_2^2\not=0, 
\text{ and } d\ge 2.
\end{equation}

\begin{defi}
Number $d$ will be called the order of tangency of $\db_0M_1$ and $\db_0M_2$ at zero.
\end{defi}

Consider the following set in $\rr^4$:
\[
W' = \{\left(x_1, 0, \alpha\left[(r_1x_1^d, r_2x_1^d) + O(|x_1|^{d+1})\right]\right): x_1
\in (-1,1), \alpha\in \rr\}.
\]
Its geometric meaning is this: for a point $(x_1,0)\in \db_0M_1$ on the $\rr^2_{x_1,y_1}$-plane 
we take the line $l_{x_1}$ in the $\rr^2_{x_2,y_2}$-plane joining $0$ with the point $(r_1x_1^d, 
r_2x_1^d) + O(|x_1|^d)$. Notice that $\left(x_1,0,(r_1x_1^d, r_2x_1^d) + O(|x_1|^d)\right)\in 
\db_0M_2$, \ie  both $M_k$, $k=1,2$ are attached to $W'$. The corresponding values of the 
parameter $\alpha$ is $0$ for $\db_0M_1$ and $1$ for $\db_0M_2$. 
Observe now that $W'$ is the subset of the following smooth surface in $\rr^4$:
\begin{equation}
W = \left\{\left(x_1,0, \alpha \left[(r_1, r_2) + O(|x_1|)\right]\right): x_1 \in (-1,1), 
\alpha\in \rr\right\}.
\end{equation}

This can be seen by changing parameter $\alpha$ to $\alpha x_1^d$.  We see that
\begin{equation}
\Psi:(t, \alpha) \to \left(t,0, \alpha \left[(r_1, r_2) +  O(|t|)\right]\right)
\end{equation}
is a parametrization of $W$. Since $\frac{\partial \Psi}{\partial t}(0) = (1,0,*,*)$
and $\frac{\db \Psi}{\db\alpha}(0) = (0,0, r_1, r_2)$ this $\Psi$ is an embedding,
\ie $W$ is smooth in a neighbourhood of the origin. Notice furthermore that 
\begin{equation}
\eqqno(m1-m2-w)
\Psi^{-1}(M_1\cap W) = \{\alpha =0\} \qquad\text{ and } \qquad \Psi^{-1}(M_2\cap W) = 
\{\alpha =t^d\}.
\end{equation}

Extend $\Psi$ to a neighbourhood of zero in $\rr^3$ as follows
\begin{equation}
 \eqqno(ext-psi)
 \Psi (t,\tau , \alpha ) = \left(t,\tau , \alpha \left[(r_1, r_2) +  O(|t|)\right]
 \right), 
\end{equation}
and notice that this extension stays to be an embedding which maps $\{(t,\tau , 0): t^2+\tau^2
< 1,\tau >0\}$ to $M_1$. Indeed $\frac{\db \Psi}{\db \tau}(0) = (0,1,0,0)$. Now extend $\Psi$ 
to a diffeomorphism of a neighbourhood of zero in $\rr^4$ arbitrarily (preserving orientation)
and denote still by $\Psi (t,\tau , \alpha , \beta)$ the local diffeomorphism obtained.
From \eqqref(m1-m2-w) we see that $\Psi^{-1}$ brings $M_1$ and $M_2$ to the form as in \eqqref(tot-r1).

\smallskip 
Assuming that $M_k$-s are brought to the form as in \eqqref(tot-r1) we extend $M_1$ to a disk 
$\widetilde{M}_1 = \{(x_1,y_1)\in \rr^2_x: x_1^2+y_1^2< 1\}$ and $M_2$ 
extend as a graph $\{ (x_2,y_2) = \tilde\phi (x_1,y_1)\}$ over $\widetilde{M}_1$, where 
\begin{equation}
\eqqno(gr-tild)
\tilde\phi_1 (x_1, y_1) = \phi_1 (x_1, - y_1) \quad\text{ and }\quad 
\tilde\phi_2 (x_1, y_1) = - \phi_2 (x_1, - y_1) \quad\text{ for } \quad y_1 < 0.
\end{equation}
$\widetilde{M}_2$ is an embedded Lipschitz surface in $\rr^4$ intersecting the smooth surface 
$\widetilde{M}_1$ at the origin only. Therefore the index of intersection of $\widetilde{M}_1$ 
with $\widetilde{M}_2$ at zero is well defined as the algebraic number of points of 
intersections of $\widetilde{M}_1$ with a small smooth perturbation of $\widetilde{M}_2$. 
\begin{defi}
\label{index2}
Call this number the intersection index of  half-disks $M_1$ and $M_2$ at $p_0$. 
\end{defi}

This index is obviously independent of the diffeomorphisms and perturbations involved and
can be any integer number, both positive and negative.

\medskip\noindent{\slsf Case 2. $M_k$ meet at zero.} In this case $M_2$ is a graph over 
$M_1^-\deff -M_1^+$. This is again after shrinking and deforming  $M_2$ if necessary.  
Let $\phi = (\phi_1,\phi_2): M_1^-\to \rr^2_{x_2,y_2}$ defines $M_2$ as a graph. Again 
we assume that $\phi(r_0(x_1))$ has finite order of vanishing at zero, \ie 
\[
\phi(x_1) = (r_1x_1^d, r_2x_2^d) + O(|x_1|^{d+1}) \quad\text{ with }\quad 
 r_1^2+r_2^2\not=0.
\]

Literally repeating constructions from the Case 1 we construct a surface $W$, having 
the same equation as before, and such that $M_k$ are both attached to $W$. 
After that we bring $M_2$ to the form as in \eqqref(tot-r1) with the only difference 
that $\phi$ is defined on $M_1^-$ and not on $M_1^+$. Then extend $M_1$ to the 
disc $\widetilde{M}_1$ as there and $M_2$ to a graph $\widetilde{M}_2$ 
over $\widetilde {M}_1$ given be the same relations as \eqqref(gr-tild).
Finally define the index of intersection of $M_1$ with $M_2$ as that of 
$\widetilde{M_1}$ with $\widetilde{M}_2$. 

\newprg[COMP.touch-c]{Index of intersection of two tangent complex disks in $\cc^2$} 

By an embedded analytic (or complex) half-disk $C$ in $\cc^2$ we mean the image of an 
$\calc^{\infty}$-embedding of $u:\Delta^+\to \cc^2$ which is holomorphic in the interior 
of $\Delta^+$. 
Consider two embedded analytic half-disks $u_k:\Delta^+\to \cc^2$ such that $u_k(0)=0$ for 
$k=1,2$ and the tangent vectors $\vect_0^k=du_k(0)[\db_{\xi}]$ to both $C_k=u_k(\Delta^+)$ 
at zero are collinear to $\pm e_1$. By holomorphicity we see that $\wect_0^k\deff du_k(0)
[\db_{\eta}] = \pm ie_1$. The degree of tangency of $C_1$ and $C_2$ is still
assumed to be finite, see Example \ref{inf-touch}.  We want to prove that the index 
of intersection of $C_1=u_1(\Delta^+)$ with $C_2=u_2(\Delta^+)$ at zero is positive, 
in fact it is $\ge 2$, as it should be.

\smallskip\noindent{\slsf Case 1. $C_k$-s touch each other at zero.} 
Denote by $\pr_1$ the standard projection onto the $z_1$-coordinate plane in 
$\cc^2$ and $\pr_2$ onto the $z_2$-plane. After shrinking, if necessary, we can assume 
that the projection  $(\pr_1 \circ u_k)(\Delta^+)$ of both of them is the same in the 
case of touching, denote it as $D\subset \cc_{z_1}$. 

Applying the Riemann mapping theorem to $D$ we can assume 
that $D=\Delta^+$. Consider the map 
\begin{equation}
\eqqno(map-down)
\pr : (z_1,z_2) \to \left(z_1, z_2 - (\pr_2\circ u_1)(z_2)\right).
\end{equation}
Note that $\pr \circ u_k$ are both an analytic half-disks, touching the direction $e_1$ at zero
and that $\pr\circ u_1$ takes it values in $\cc_{z_1}$-plane, in fact its image is $\Delta^+$ 
itself. Moreover, since the Riemann mapping function extends diffeomorphically to a neighborhood 
of the boundary all what we deed up to now are diffeomorphisms in a neighborhood of zero. 


From now on $C_1 = \Delta^+ \subset \cc_{z_1}$ and $C_2$ is a graph over $C_1$, \ie 
$C_2 = \{(z_1,z_2): z_1\in \Delta^+, z_2 = \phi (z_1)\}$ where $\phi$ holomorphic 
in the interior of $\Delta^+$ and smooth up to the boundary. Extend $\phi$ to a 
smooth function on $\Delta$. Since Cauchy-Riemann relations for $\phi$ are valid
up to the boundary we see that all partial derivatives of $\tilde\phi$ at zero 
containing a derivative with respect to $\bar z_1$ vanish. Therefore 
\begin{equation}
\eqqno(index-phi)
\tilde \phi (z_1) = az_1^d + O(z_1^{d+1})
\end{equation}
in a neighborhood of zero. Notice that all coefficients of the Taylor series at
zero are determined by the values of $\phi$ on $\Delta^+$. In particular such are 
$d$ and $a$. \eqqref(index-phi) restricted to $\Delta^+_r$ can be 
considered as a comparison relation \eqqref(9.1.1) for analytic half-disks
$u_1(\zeta) = (\zeta ,0)$ and $u_2(\zeta) = (\zeta , a\zeta^d + O(\zeta^{d+1}))$.
These half-disks do not satisfy the reality condition on $\db_0\Delta^+$ and 
cannot be reflected right away. But since they are defined on the whole of 
$\Delta^r$, and the relation \eqqref(index-phi) holds there, we can apply 
the reasoning of the proof of Theorem \ref{b-pos-thm} to conclude that the 
"natural" index of boundary intersection here should be $d$.

\smallskip 
To prove this repeat the construction for the surface $W$ from the previous subsection to  
$M_k=C_k$.  Notice that the surface $W$ is totally real in our (complex) case. 
Indeed, $\frac{\partial \Psi}{\partial t}(0) = (1,0,*,*)$ and $\frac{\db \Psi}{\db\alpha}
(0)  = (0,0, r_1, *)$ are linearly independent over $\cc$, where $\Psi$ is the 
 diffeomorphism  constructed there. We  take $\Psi^{-1}$ as a ``rectification map''. 
 The structure $J\deff \Psi^{-1}_*J\st$ might be non  standard  on $\rr^2_{t,\alpha}$, 
 but the latter is $J$-totally real and, as  it was noticed in  the Remark 
 \ref{red-rem} we can make $J|_{\rr^2_{t,\alpha}}$ standard by one more 
 diffeomorphism which is identity on  $\rr^2_{t,\alpha}$. 

 \smallskip Analytic half-disks $A_k\deff \Psi^{-1}(C_k)$ are both attached 
 to $\rr^2_{t,\alpha}$  and we observe that both constructions used in this 
 paper, (\ie for real half-discs of the present section and for complex ones 
 attached to a totally real submanifold), give the same index of intersection  
 of $C_1$ with $C_2$. This proves  the positivity of our index. To compute it
 explicitly   remark that 
 \begin{equation}
 \eqqno(a1-a2-1)
 A_1\cap \rr^2_{t,\alpha} = \{\alpha =0\} \text{ and }A_2\cap \rr^2_{t,\alpha} = 
 \{\alpha = t^d\}. 
 \end{equation}
 Use $t+i\tau$ and $\alpha + i\beta$ as complex coordinates in $\cc^2$ - 
 the source of $\Psi$. Denote by $u_k:\Delta^+\to \cc^2$ the $J$-holomorphic 
 parametrizations  of $A_k$. By Comparison Lemma \ref{comp-lem-r} we 
 know that for some holomorphic  $\psi (\zeta)  =\zeta + O(\zeta^2)$, real 
 for real $\zeta$, we have 
 \begin{equation}
 \eqqno(a1-a2-2)
u_2(\zeta) = u_1(\psi (\zeta)) + \zeta^{\nu}w(\zeta), 
\end{equation}
where $\nu >\mu =1$ in our case. Comparing \eqqref(a1-a2-2) with \eqqref(a1-a2-1) we see
that $\nu$ can be nothing but $d$. In the proof of  Theorem \ref{b-pos-thm}
we had shown that the index of intersection is equal to $\nu$ and therefore is equal
to $d$ in our case. In any event the index of boundary intersection of complex disks 
is proved to be equal to the degree of tangency of these disks.
 
\smallskip Finally, index is equal to $\adyn$ if the lowest degree term of the 
Taylor expansion of $\phi$ is one, \ie $C_1$ and $C_2$ intersect transversally.
Case 1 is proved.

\smallskip\noindent{\slsf Case 2. $C_k$-meet at zero.} 
The projection of $C_1$ to the $z_1$-axes denote by $D^+$ and that of $C_2$ by $D^-$. Let $S^{\pm}$ be circles touching $\db D^{\pm}$ at zero from the interior. Shrinking $D^{\pm}$
up to the discs bounded by these circles, and therefore shrinking  $C_1$ and $C_2$, we can assume that each of $\db D^{\pm}$ contain  an arc of the corresponding circle, denote
them still as $S^{\pm}$. Find an automorphism of the Riemann sphere sending $S^+$ 
to $S^-$ and interior of $D^+$ to the exterior of $D^-$. This translates $C_1$ to a graph 
over a domain with boundary containing $S^-$. In order to achieve all this one may need
to shrink $D^+$ once more. Finally take a biholomorphism sending $S^-$ to the 
interval $(-1,1)\times \{0\}$. We realised $C_1$ and $C_2$ as a graphs over $\Delta^+$
and $\Delta^-$. Repeat the construction of the surface $W$ from the previous subsection 
for the Case 2 and use the Comparison Lemma \ref{comp-lem-op} for meeting complex 
discs attached to the totally real submanifold to conclude that the intersection
index is equal to the order of tangency of $C_1$ and $C_2$. Theorem \ref{touch-thm}
is proved.

\begin{exmp} \rm
\label{inf-touch}
Analytic half-discs in $\cc^2$ can touch each other up to an infinite order at their
boundary point. Set, for example, 
\[
C_1 \deff \Delta^+\times\{0\}  \quad\text{ and }\quad C_2  \deff 
\left\{\left(z_1,e^{-\frac{1}{\sqrt[3]{z_1}}}\right):z_1\in \Delta^+\right\}.
\]
Notice that the function $e^{-1/w}$ is of class $\calc^{\infty}$ in the sector 
\[
S^+ \deff \left\{w \in \Delta^+: 0\le \arg w \le \frac{\pi}{3}\right\}\setminus \{0\}.
\]
Moreover all its derivatives tend to zero exponentially when $w\to 0, w\in S^+$. 
Making substitution $w(z)=\sqrt[3]{z}$ we preserve the latter property since 
all derivatives of $w(z)$ have at most polynomial growth at zero.
\end{exmp}

\newsect[ADJ]{Application II: Adjunction Formula}

\newprg[B-POS.dbl]{Doubling Along a Totally Real Submanifold} 

In order to be able to write the Adjunction Formula for complex curves with boundary on
a totally real submanifold we shall need the notion of a homological self-intersection 
of the Schottky double of such curve. But first recall the construction of the Schottky 
double. Let $(C,\d C) $ be a smooth connected compact complex curve with boundary. Take two 
copies $C^{\pm}$ of it with opposite orientations. Glue them together along the boundary
using the identity map. The compact complex curve $C^d= C^+\sqcup_{\d C}C^-$ thus 
obtained  is called the \emph{Schottky double} of $C$. $C^d$ admits a natural 
anti-holomorphic involution $\tau_C :C^d\to 
C^d$ interchanging $C^+$ with $C^-$.  Assume we are given a complex vector bundle 
$E$ over $C$ and a real sub-bundle $F$ of $E|_{\d C}$, \ie for every $\zeta\in \db C$ 
the fibre $F_\zeta$ is a totally real subspace of $E_\zeta$ of maximal dimension. 
Denote by $E^{-}$ the copy of $E$ over $C^{-}$ with the opposite complex structure. 
Glue $E^+=E$ and $E^-$ over the boundary by an anti-complex involution fixing $F$. 
Denote by $E^d$ the complex vector bundle over $C^d$ thus obtained. $E^d$ admits an 
anti-holomorphic involution $\tau_E$ over $\tau_C$, \ie it sends $E_\zeta$ anti-linearly 
to $E_{\tau_C(\zeta)}$ fixing $F_\zeta$ for $\zeta\in \db C$.

\smallskip Let us pass to complex dimension two. We consider an almost complex 
surface $(X,J)$, \ie a $4$-dimensional real manifold equipped with an almost 
complex structure $J$. We fix a $J$-totally real surface $W$, which is supposed 
to be closed as the subset of $X$. Let $(C,\db C)$ be as above and $u:\bar C 
\to X$ be a continuous map, which is $J$-holomorphic and primitive in $C$, such
that $u(\db C)\subset W$ and $u(C)\subset X\bs W$. We assume that $J\in 
\calc^{1,\alpha}$ and $W\in \calc^{2,\alpha}$. Then $u$ is of class 
$\calc^{2,\alpha}$ up to $\db C$. Denote $u(C)$ as $M$ and consider $M$ as a 
$J$-complex curve in $X$ with 
boundary on $W$. Assume in addition, that $M$ is immersed in $X$ near the 
boundary. That is $M$ has no cusps on $W$, at most self-intersections, but not 
necessarily transverse. Fix a $J$-Hermitian metric $h$ on $X$. Our
curve $M$ meets the boundary $\db B_\delta$ of every  tubular 
$\delta$-neighbourhood $B_\delta$ of $W$ transversely for $\delta>0$ 
small enough. This is so because $M$ intersects $W$ Bott-transversely, 
\ie $\angle (T_xW,T_xM)>0$ at every point $x\in \db M$. Denote by 
$B_{\eps ,2\eps} \deff B_{2\eps} \bs \bar B_{\eps}$ the layer,  where 
$\eps >0$ is taken small enough. In particular, $M$ should intersect 
every $\db B_{\delta}$ for $0<\delta < 3\eps$ transversely. Take two
copies $X^\pm_{\eps}$ of $X\setminus \bar B_\eps$, $X^+_{\eps}$ with 
the same orientation as $X$ and $X^-_{\eps}$ with the opposite. Glue 
them along $\db B_\eps$ by the identity map and let $X^d_{\eps} \deff 
X^+_{\eps} \sqcup_{\db B_\eps}X^-_{\eps}$ be the manifold obtained. 
We call it an $\eps$-double of $X$ if $W$ is clear from the context, 
or as $\eps$-double of the pair $(X,W)$. Set $M^+_{\eps} \deff M
\cap X^+_{\eps}$ and denote by $M^-_{\eps}\subset X^-_\eps$  the same 
Riemann surface but with the opposite orientation. Set $M^d_{\eps} = 
M^+_\eps \sqcup_{\db M^+_{\eps}} M^-_\eps$ and call it an $\eps$-double 
of $M$. $M^d_\eps$ naturally realises as a surface in $X^d_\eps$. 
Moreover, we shall prove (but not use this) that $X^d_\eps$ admits
an almost complex structure $J^d_\eps$ extending $J$ from $X^+_\eps$
such that $M^d_\eps$ is $J^d_\eps$-complex. $X^d_{\eps}$ admits a natural
involution $H$, which arises from the natural identification of 
$X^+_{\eps}$ with $X^-_{\eps}$. Notice that $H$ exchanges $M^+_{\eps}$ 
and $M^-_{\eps}$. 

\begin{defi}
\label{m-double}
We call the pair $(X^d_{\eps},M^d_{\eps})$ an $\eps$-double of $(X,M)$. 
\end{defi}

\begin{rema} \rm Let us make a few remarks about this object. 

\smallskip\sli $M^d_\eps$ is naturally equipped with a complex structure,
which is not isomorphic to the 

\quad Schottky double of $M$, but to a small deformation of it. 
Nevertheless this is

\quad appropriate for us, since wee need only self-intersection of it in 
$H_2(X^d_\eps,\zz)$.

\smallskip\slii If $x\in M^+_{\eps}\cap \ti M^+_{\eps}$ is a point of the 
transverse intersection and $y = H(x)$ is the

\quad   symmetric point then indices of intersection at $x$ and $y$ are
equal.

\smallskip\sliii Let $\Sigma$ be any embedded compact surface in $X^d_{\eps}$, 
set $\Sigma^{\pm}_{\eps} = \Sigma\cap X^{\pm}_{\eps}$. We say that 
$\Sigma$ 

\quad is "symmetric" if  $H(\Sigma^+_{\eps}) = \Sigma^-_{\eps}$.
Our $M^d_{\eps}$ is symmetric. Denote by $[M^d_{\eps}]^2$ the 
homo-

\quad logical self-intersection of $M^d_{\eps}$ in  $X^d_{\eps}$,
\ie the self-intersection of its homology class 

\quad $[M^d_{\eps}] $ in $H_2(X^d_\eps , \zz)$. Or, which is the same,
the intersection of $M^d_{\eps}$ with its generic 

\quad perturbation. The latter can be made symmetric, if needed.

\end{rema}

\begin{figure}[h]
\begin{minipage}[b]{.40\textwidth} 
\caption{Surgery.}
\label{Fig-chir} 
Vector $\vect$ is a direction in the $\rr^2_{x_1,x_2}$-plane. Tubular
$\eps$-neighbourhood $B_\eps$ of the line $\<\vect\>$ removed. The 
half-hyperbola is swept by half-circles from $C_1^+$ to $C_2^+$ and 
then to $C_3^+$, each for a different $\vect_1,\vect_2,\vect_3$. 
We patch the half-hyperbola to our half-discs 
$\calk_1^+$ and $\calk_2^+$ by $C_1^+$ and $C_2^+$ correspondingly. 
After doubling with respect to $\db B_\eps$ we obtain a handle filled
from bottom by discs with boundaries $C_1^-\cup C_1^+$ and from the
top by $C_3^+\cup C_3^+$, which after patching became the doubled 
$\calk_1^+$ and $\calk_2^+$.
\vspace{5pt} 
\end{minipage}%
\hfill%
\begin{minipage}[b]{78 true mm}
\centering
\includegraphics[width=68 true mm]{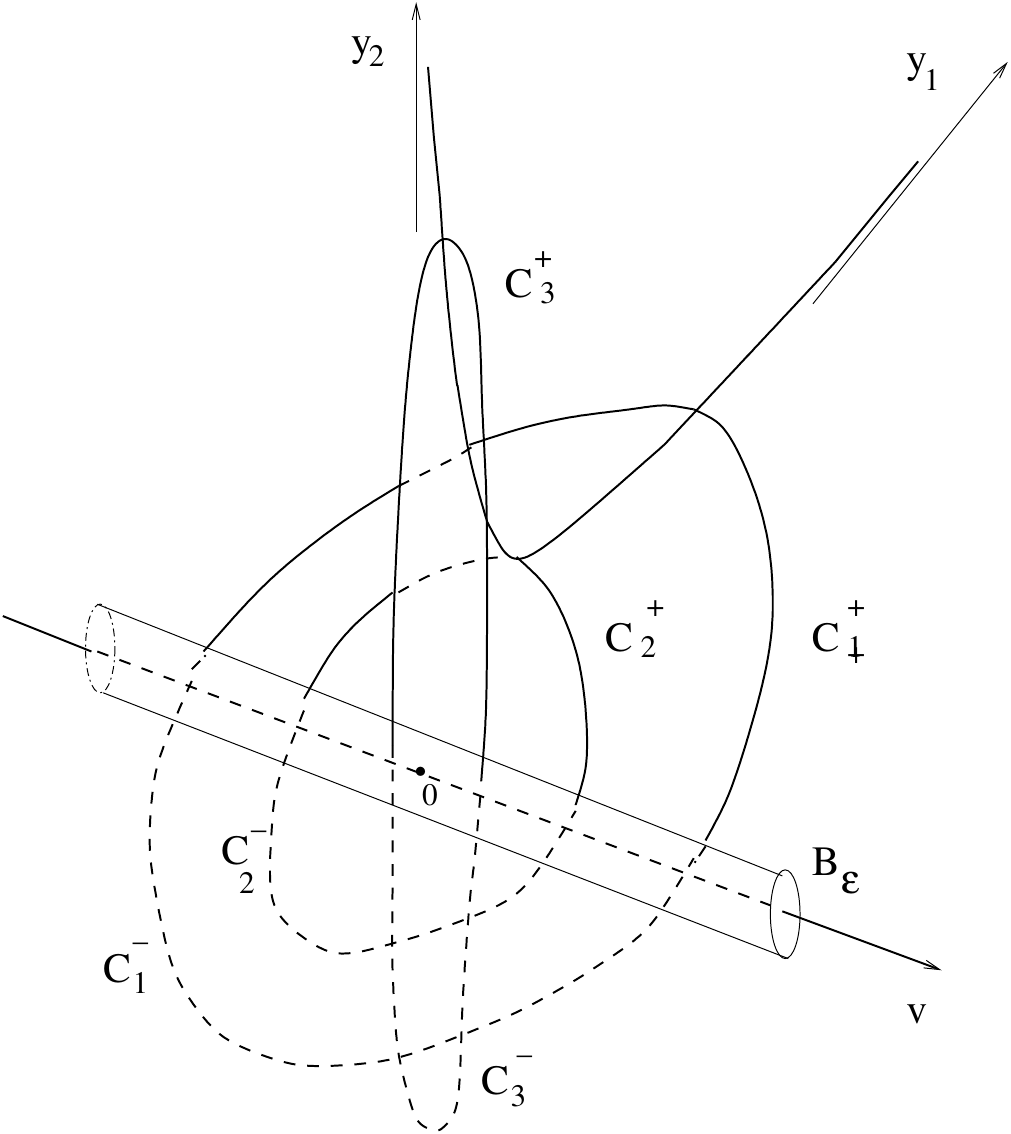}
\end{minipage}
\end{figure}

\smallskip We shall need to perform the following elementary 
surgery near the transverse boundary self-intersection points of 
$M$. Let $x\in M\cap W$ be such a point. In appropriate coordinates 
we can represent $W$ as $\rr^2$ and $M$ in a neighbourhood of  
$x=0$ as the half-cross 
\[ 
\calk^+\deff \{(z_1,z_2)\in \Delta^2: \im z_1\ge 0, \im z_2\ge 0, 
z_1\cdot z_2 = 0\},
\]
where the branches of $M$ are the coordinate half-discs $\calk^+_1 = \Delta^+
\times \{0\}$ and $\calk^+_2 = \{0\}\times \Delta^+$. Take $\delta >0$ small enough 
and replace $\calk^+$ by the "half-hyperbola" 
\[
H^+ \deff \{(z_1,z_2)\in \Delta^2: \im z_1 \ge 0, 
\im z_2 \ge 0, z_1\cdot z_2 = - \delta \}, 
\]
which is glued to the half-cross $\calk^+$ near the outer boundary of the 
latter. By the outer boundary of the half-disc we mean $\db^{\text{ot}}\Delta^+
\deff \db\Delta^+\setminus (-1,1)$, and then $\db^{\text{ot}}\calk^+\deff 
\db^{\text{ot}}\calk_1^+\cup \db^{\text{ot}}\calk_2^+$. Denote by 
$\widetilde{M}^{\delta}$ the surface obtained. It is not difficult to see that 
for $0<\eps \ll \delta$ the double $\left(\widetilde{M}^{\delta}\right)^d_\eps $ 
is different from the double $M_\eps^d$ by a handle with top and bottom filled by 
the discs. The latters are the doubles of the half-discs
$\calk_1^+\setminus B_\eps$ and $\calk_2^+\setminus B_\eps$. As well as that
the closed surface $S$ thus obtained is contained in a neighbourhood of $\db B_\eps$
in $X^d_\eps$. But the latter is locally diffeomorphic to the disc times the annulus
and has zero second homology group. We proved that 
\begin{equation}
\eqqno(homol)
\left[\left(\widetilde{M}^\delta\right)^d_\eps\right] =\left[M^d_\eps\right]
\qquad\text{ in }\qquad H_2(X_\eps^d,\zz).
\end{equation}

\begin{prop}
\label{eps-surg}
Let $M$ be a compact $J$-complex curve in $X\setminus W$ attached by 
its boundary to $W$ and without cusps on $W$. Then:

\smallskip\sli the number $[M^d_{\eps}]^2$  doesn't depend on $\eps > 0$ 
taken small enough;

\smallskip\slii this number doesn't changes under the surgery described above, provided
$0<\eps \ll\delta$.
\end{prop}

\proof (\sli We can perturb $M$  away from $W$, in such a way
that every cusp point will be replaced by that number of transverse self-intersection 
points equal to the cusp index of this point. This perturbation do not change 
$[M^d_\eps]^2$. Indeed, if $\widetilde{M}$ is such a perturbation then it coincides with 
$M$ near $W$. This implies that $\ti M^d_\eps$ coincides with $M^d_\eps$ near 
$\db B_\eps$. Which means that the homology classes of $\ti M^d_\eps$ and 
$M^d_\eps$ are the same. Assertion (\sli is proved.

\smallskip\noindent (\slii follows 
from the observation made above in \eqqref(homol). In the sequel we can assume that $M$ has only transverse 
self-intersection away from $W$. On $W$ they can be arbitrary.

\smallskip\qed 

\smallskip 
Now we are in position to introduce the following notion. 

\begin{defi}
We define the homological self-intersection number $[M^d]^2$ of the Schottky 
double of $M$ as 
\begin{equation}
[M^d]^2 = \lim_{\eps \searrow 0} [M^d_{\eps}]^2,
\end{equation}
or, due to what was just proved $[M^d]^2 = [M^d_{\eps}]^2 $ for $\eps >0$ small enough.
\end{defi}

\newprg[ADJ. masl-i]{Maslov Index}

\smallskip Let $(C,\db C)$ be a smooth compact connected complex curve with a 
non-empty boundary, $E$ a complex line bundle over $C$ and a totally real 
subbundle $F$ of $E|_{\db C}$. Take a section $\nu$ of $F$  along $\db C$ with 
non-degenerate zeros and extend it to a section of $E$ over $C$ again with 
non-degenerate zeros.  When saying that $\nu$ has non-degenerate zeroes on 
the boundary we mean that $\nu$ should be such that for every boundary zero
$x\in \db M$ of $\nu$ there exists a disc $M_x$ in the double $C^d$ centred at $x$,
which extends a corresponding half-disc neighbourhood $M^+_x$ of $x$, and 
an extension of $\nu$ to $M_x$ from $M_x^+$ as a section of $E^d$, such that 
this extension has a \emph{non-degenerate} zero at $x$.
\begin{defi}
\label{mslv-ind2}
Define the {\slsf Maslov index} $\bfmu_C(E,F)$ of the pair $(E,F)$ over $C$ 
as twice the number of interior zeros plus the number boundary zeros, all counted 
with signs.
\end{defi}

Consider now the complex bundle of rank $n$ on $C$ and a totally real subbundle 
$F$ of $E|_{\db C}$. Let $\det(E) := \wedge^n_\cc E$ and $\det(F):= \wedge^n_\rr F$
be the corresponding determinant bundles. The Maslov index $\bfmu_C(E,F)$ is 
defined as the Maslov index of the corresponding determinant bundles, \ie 
$\bfmu_C(E,F)\deff\bfmu_C(\det(E),\det(F))$.



\begin{prop}
Maslov index $\bfmu_{\db C} (E,F)$ can be computed with the help of the 
Schottky double as follows:
\begin{equation}
\eqqno(mas-dbl)
\bfmu_C (E,F) = \bfc_1(E^d)[C^d].
\end{equation}
\end{prop}
\proof Since $\bfc_1(E^d)$ is, by definition $\bfc_1(\det(E^d))$ we need 
to prove this for line bundles only. Take $\nu$ as in Definition \ref{mslv-ind2}
and double it to a section $\nu^d$ of $E^d$ over $C^d$. Therefore every inner 
zero of $\nu$, \ie zero in the 
interior part of $C$, produces a zero of its mirror $\bar\nu$ of the
same sign. To the contrary with inner case, a boundary zero of $\nu$
stays to be a  single zero of $\nu^d$. And the sum of zeroes with signs of $\nu^d$ 
is known to ve the first Chern class.

\smallskip\qed


\noindent{\bf 1.} For example, in the case of $(E,F) = (TC,T\db C)$, since $\bfc_1(TC^d)[C] = 
\bfchi (C^d) = 2\bfchi (C)$, one has 
\begin{equation}
\eqqno(mas-tang)
\bfmu_C(TC, T\db C) = 2\bfchi (C) = 4-4g-2\sigma,
\end{equation}
where $g$ is the genus of $C$ and $\sigma$ is the number of the boundary 
circles.

\smallskip\noindent{\bf 2.} If $E=E_1\oplus E_2$ and 
$F=F_1\oplus F_2$ then the additivity of the Chern class implies the additivity 
of the Maslov index
\begin{equation}
\eqqno(masl-add)
\bfmu_C (E,F) = \bfmu_C (E_1,F_2) + \bfmu_C (E_2,F_2).
\end{equation}

\smallskip\noindent{\bf 3.} Let $(X,J)$ be an almost complex manifold and $W$ a 
$J$-totally real submanifold 
of $X$ of maximal dimension. We assume that $W$ is a closed subset of $X$. 
Let $u:(C,\db C)\to (X,W)$ be a $J$-holomorphic map. Smoothness assumptions on 
$J$ and $W$ are as everywhere in this paper. 
Consider the induced bundles $E\deff u^*TX$ on $C$ and $F\deff u^*TW$ on $\db C$.
Setting $M=u(C)$ we denote the Maslov index $\bfmu_C (u^*TX, u^*TW) \deff \bfmu_C(E,F)$ 
as $\bfmu_M(TX,TW)$. Assume that $M$ has no cusps. Then the normal bundle 
$N_M$ is well defined and, moreover, the normal bundle $N^W_{\db M}$ of 
$T\db M$ in $TW$ is well defined as well. $N^W_{\db M}$ is a real sub-bundle of 
$N_M$ and the formula \eqqref(masl-add) applies. Combining \eqqref(mas-tang) with 
\eqqref(masl-add) we obtain 
\begin{equation}
\eqqno(masl-add-n)
\bfmu_M(TX,TW) = \bfmu _M(N_M,N^W_{\db M}) + 4-4g-2\sigma .
\end{equation}



\begin{rema} \rm
\label{non-or1}
{\bf a)} 
It is instructive to distinguish two cases. If $F$ is orientable then the 
section $\nu$ can be taken without zeros on $\db C$ at all. As the result 
tha Maslov index is even in this case. If $F$ is 
non-orientable, a vector field with values in $F$ on $\db C$ must vanish,
and requirement is that it intersects the zero section transversely. And 
then the Maslov index will be the sum of boundary zeros plus the twice 
sum of inner zeros, all computed with signs.

\smallskip\noindent{\bf b)} Maslov index is well defined, even if $(C,u)$ has cusps 
(also on the boundary), and is invariant under perturbations both inner and 
boundary, provided the latter ones keep $M$ being attached to $W$. 
\end{rema}

\begin{prop}
\label{eps-mas}
Under the assumptions of Proposition \ref{eps-surg} assume, moreover, that
$M$ has no cusps \emph{at all}. Then $[M^d_\eps]^2$
is equal to $\mu_{\d M}(N_M,N^W_{\db M}) + 4\delta^{(i)}$ for $\eps > 0$ small
enough.
\end{prop}
\proof Maslov index $\mu (N_M,N^W_{\db M})$ is equal to 
$c_1(N_M^d)[M^d]$ by \eqqref(mas-dbl), \ie to the self-intersection number 
of $M^d$ in $N^d_M$. Now we shall compare it with the self-intersection 
number of  $M^d_{\eps}$ in $N_{M^d_{\eps}}$ for $\eps \sim 0$. Take a 
smooth section  $\nu$ of the normal to $M$ bundle $N_M$, \ie orthogonal 
to $TM$ at every point, which takes values in $TW$ at points of $\db M$ 
and is orthogonal to $T\db M$. In addition, $\nu$ is takrn to have 
non-degenerate zeroes only.

\smallskip   Chern number  
$c_1(N_M^d)[M^d]$ is the twice number of inner zeros of $\nu|_{\db M}$ 
plus number of its boundary zeros, all with signs. This can be 
reformulated as follows. Our section $\nu$ of $N_M$ naturally 
defines the section $\nu^d$ of $N_M^d$ over $M^d$ and the number of 
zeros of this $\nu^d$ counted with signs is actually $c_1(N_M^d)[M^d]$.  

\smallskip  On the other hand the very same $\nu$ defines naturally
the sections $\nu^+ = \nu$ of $N_M^+$ over $M^+_{\eps}$ and $\nu^- = 
\bar \nu $ of $N_{M^-_{\eps}}$ over $M^-_{\eps}$. One can glue
these sections  to a section $\nu^d_{\eps}$ of $N_{M^d_{\eps}}$ 
over $M^d_{\eps}$ in $X^d_{\eps}$ as follows. Consider two cases.

\begin{figure}[h]
\centering
\includegraphics[width=119 true mm]{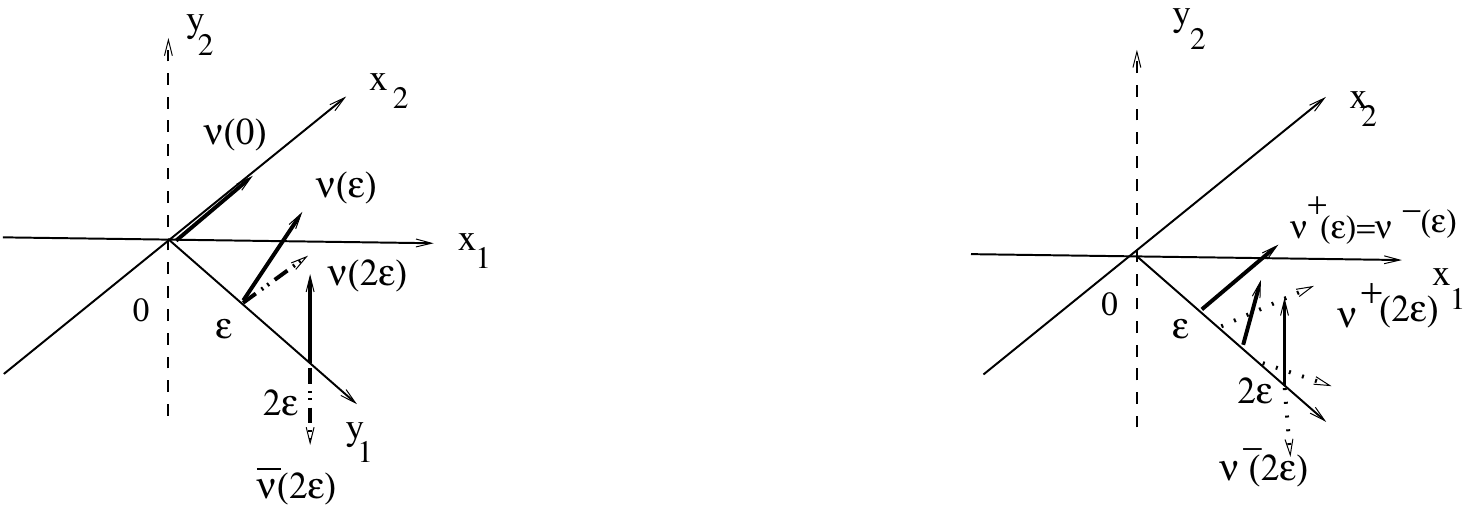}
\caption{Gluing $\nu^+$ and $\nu^-$. Case 1.}
\label{Fig-rot-nu1}  

\vspace{3pt}
\begin{minipage}{0.72\linewidth} 
\setlength{\parindent}{1.5em}
\normalsize
\noindent On the left the vector field $\nu$ near the point
$x=0$ is described. Our curve $M$ is $\Delta^+ = \{(x_1,y_1):
y_1 \ge 0\}$ is attached to $W$, which is $\rr^2$ with coordinates
$x_1,x_2$. Dependence of $\nu$ on variable $x_1$ is ignored on this
picture. At zero $\nu$ is $\db /\db x_2$. On the right we make 
our vector field rotate when $y_1$ changes from $2\eps$ to $\eps$
quicker  in order to achieve the value $\db/\db x_2$ at $y_1=\eps$.
And then running from $\eps$ to $2\eps$ the vector field achieves the value 
$\bar\nu (2\eps)$. 
\end{minipage}
\end{figure}

\smallskip\noindent{\slsf Case 1.}{ \it $\nu$ had no zeros on $\db M$. }
Notice that such a section always exists if $W$ is orientable. Then we proceed 
as explained on the Picture \ref{Fig-rot-nu1}. As the result the obtained 
$\nu^d_{\eps}$ has no zeros on $M\cap\db B_{\eps}$ and therefore 
has the same number of zeros on $M^d_{\eps}$ with the same signs 
as $\nu^d$ on $M^d$.

\smallskip\noindent{\slsf Case 2.}{ \it $\nu$ vanishes at some point 
$x\in \db M$.}  Then $\nu^d_{\eps}$
will vanish with with the same sing on some close point on $\db B_{\eps}$.
Consequently, again the number of zeros of $\nu^d_{\eps}$ is the same 
as of $\nu^d$. 

\begin{figure}[h]
\centering
\includegraphics[width=119 true mm]{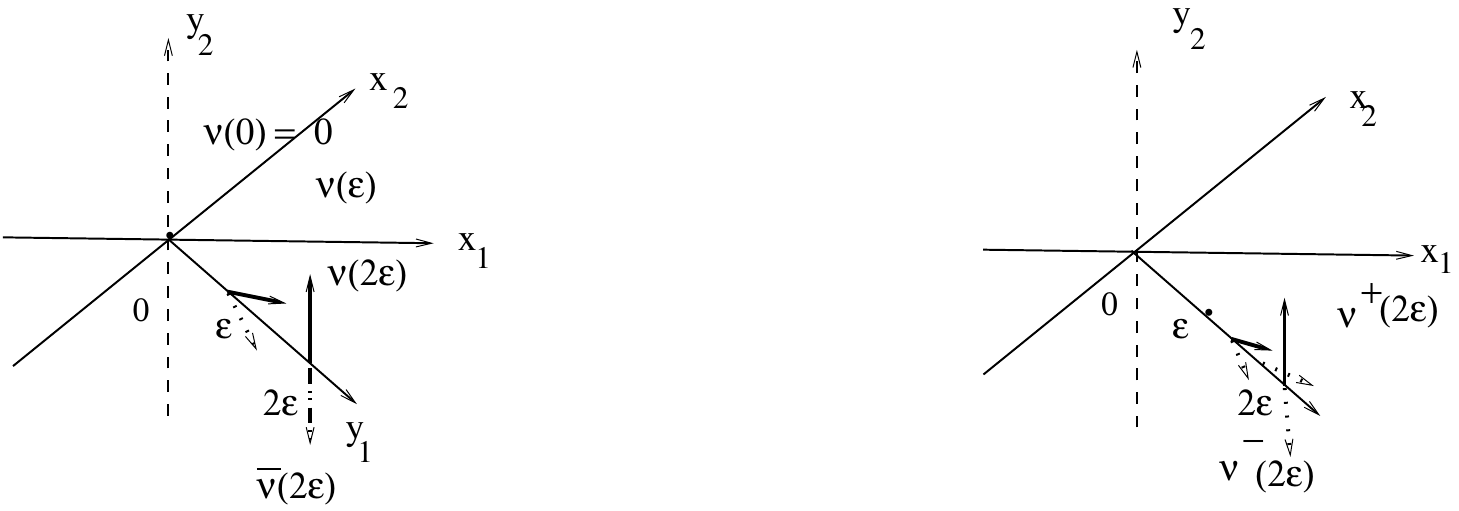}
\caption{Gluing $\nu^+$ and $\nu^-$. Case 2.}
\label{Fig-rot-nu2}  

\vspace{3pt}
\begin{minipage}{0.72\linewidth} 
\setlength{\parindent}{1.5em}
\normalsize
\noindent This picture is similar to the previous one, 
only $\nu$ vanishes at zero. 
\end{minipage}
\end{figure}

\smallskip 
We proved that twice the number of inner zeroes of $\nu$ plus the number of its boundary 
zeroes is equal to the number of zeroes of $\nu^d$, all with signs. The first is the 
normal Maslov index in question and the second  is the self-intersection number 
of $M^d_\eps$ in $N_{M^d_\eps}$. The homological self-intersection number of 
$M^d_\eps$ is its self-intersection number in $N_{M^d_\eps}$ plus twice the number 
of its self-intersection points. These points come in symmetric pairs and therefore 
their number is equal to the twice the number of self-intersection points of $M$.
All together we obtain that $[M^d_\eps]^2 = \mu (N_M, N_{\db M}^W) + 4\delta^{(i)}$.
Proposition is proved. 

\smallskip\qed

\newprg[ADJ. masl-c]{Maslov Class}
We shall need to interpret the Maslov index as the value of the certain relative 
cohomology class from $\sfh^2(X,W;\zz)$ on $M\in \sfh_2(X,W;\zz)$. This class 
is called the \emph{Maslov class}. Let $X$ be a paracompact Hausdorff locally 
contractible topological space. Notice that for such topological spaces we have 
equivalence of singular  and  \v{C}ech cohomology theories. In addition, by the 
classical result of Milnor \cite{Mi} such  space $X$ has homotopy type of a 
$CW$-complex.  Moreover, let $Y$ be a closed subset of $X$ which is also a locally
contractible topological space. We assume that $Y$ has finite number of connected 
components (this is, in fact, not needed) and denote by $X/Y$ the quotient space.
By the latter, with some ambiguity, we understand the contraction of connected 
components of $Y$ separately. More precisely: connected components $Y_1,...,Y_N$ 
of $Y$ are contracted to the \emph{distinct} points $y_1,...,y_N$ of $X/Y$. The 
quotient space again $X/Y$ has homotopy type of a $CW$-complex. Next, let $E$ be a 
complex vector bundle over $X$ of rank $n$, $F$  a real vector sub-bundle of the 
restriction $E|_Y$ of  rank $n$ such that for every point $y\in Y$ the vector space 
$F_y$ is a totally real subspace of $E_y$. Consider the complex line bundles 
$\det(E) := \wedge^n_\cc E$ and $\det^2(E) := \big( \det(E) \big)^{\otimes2}$, as 
well as  real line bundles $\det(F):= \wedge^n_\rr F$ and $\det^2(F) := \big( \det(F)\big)^{\otimes2}$. Then $\det^2(F)$ is a trivial vector bundle over $Y$: 
$\det^2(F)\cong Y \times\rr$. Fix a global non-vanishing section $s$ of  $\det^2(F)$ 
over $Y$. Then $s$ is a non-vanishing section $\det^2(E)$ over $Y$. By the theorem 
of Tietze-Urysohn, see \cite{U}, section $s$ extends as a non-vanishing section 
$\ti s$ of $\det^2(E)$ to some open neighbourhood $U$ of $Y$ in $X$.  This implies
that $\det^2(E)$ is (isomorphic to) a pull-back of some complex line bundle $L$ 
over the quotient space $X/Y$ wrt.\ the quotient map $q:X \to X/Y$, such that 
the section $\ti s$ defines a local 
trivialisation of $L$ over $U/Y$. The latter is a neighbourhood of $y_1,..,y_N$.
Or, better $\ti s$ is a couple of sections $\ti s_1,...,\ti s_N$ each over 
a neighbourhood $U_i/Y_i$ of $y_i$.

\smallskip The classifying space $BU(1)$ for $\bfu(1)$ and the total  space 
$EU(1)$ of the universal $U(1)$-bundle are the infinite dimensional projective space 
and the infinite dimensional sphere correspondingly:
\[
\cc\pp^\infty \deff \lim_{\substack{\longrightarrow}}\cc\pp^n \qquad\text{ and }
\qquad \sph^\infty \deff \lim_{\substack{\longrightarrow}}\sph^n.
\]
Our bundle $L$ is the pull-back of the associate universal line bundle $\calo 
(1)$ by some continuous map $f:X/Y\to BU(1)$. The first Chern class 
$\bfc_1(L)\in \sfh^2(X/Y;\zz)$  is the pull-back by $f$ of the universal Chern class
$\bfc_1(\calo(1))$ and the pull-back of $\bfc_1(L)\in \sfh^2(X/Y;\zz)$ by the 
quotient map $q:X\to X/Y$ lies in $\sfh^2(X,Y;\zz)$. Indeed, $q$ sends relative 
modulo $Y$ cycles in $X$ to cycles in $X/Y$. Therefore by duality a cocycle $\sigma$ 
representing a class from $\sfh_2(X/Y;\zz)$ lifts to a cocycle from $\sfh_2(X,Y;\zz)$.

\begin{defi} 
The obtained cohomology class $q^*\bfc_1(L)\in \sfh_2(X,Y;\zz)$ is called the 
\emph{Maslov class} of the pair $(E,F)$ and is denoted as $\bfmu (E,F)$.
\end{defi}

Let $(X_1, Y_1)$ be another pair of paracompact Hausdorff locally contractible topological 
spaces, and $f: (X_1, Y_1) \to (X, Y)$ a continuous map. Then the pull-back gives us a complex
vector bundle $E_1 := f^*E$ over $X_1$ and a totally real subbundle $F_1 := f^*F$ over $Y_1$. 
Since the constructions used in the definition of the Maslov cohomology class are functorial, 
\ie compatible with the pull-back we can conclude the formula
\begin{equation} 
\eqqno(f^*mu)
\bfmu(f^*E,f^*F) = f^*\bfmu(E,F) \in \sfh^2( X_1,Y_1; \zz).
\end{equation}

Let $(X,Y)$ be as above and let $\sigma \in \sfh_2( X,Y; \zz)$ be any relative
homology class represented by a cycle $\sum_j a_j [u_j]$. Here the sum is finite, $a_j$ 
are integer coefficients, and $u_j: \triangle \to X$ are continuous maps from the 
standard $2$-simplex $\triangle$. Using the definition of the boundary of a singular 
chain one can show that there exists a compact oriented surface $C$ with boundary 
consisting of a union of circles (possibly empty), its triangulation $C= \cup_j C_j$, 
a continuous map $u: C \to X$, and orientation preserving homeomorphisms $\phi_j: C_j 
\to \triangle$ such that $u|_{C_j} = u_j \circ \phi_j$. Moreover, $u$ maps $\db C$ to 
$Y$, so that we have a continuous map of pairs $u: (C, \db C) \to (X,Y)$. This means 
that the class $\sigma \in \sfh_2( X,Y; \zz)$ is a direct image $u_*[C]$ of the 
fundamental homology class $[C] \in \sfh_2(C, \db C; \zz)$. Therefore by \eqqref(f^*mu)
we can conclude the formula 
\begin{equation}  
\eqqno(mu.[S])
\bfmu(E,F)(\sigma) = \bfmu(E,F)(u_*[C]) = u^*\bfmu(E,F)[C] = \bfmu(u^*E,u^*F)[C] 
\in \zz.
\end{equation}
This gives us the value of the Maslov cohomology class $\bfmu(E,F)$ on relative 
homology classes from $\sfh_2( X,Y; \zz)$. 



\begin{prop}
\label{masl-cl-prop}
Let $X$ be an almost complex surface and $W$ be a totally real surface in 
$X$, closed as a subset of $X$. Then for every compact complex curve 
$M = u(C)$ with boundary on $W$ one has 
\begin{equation}
\eqqno(masl-cl-ind)
\bfmu (TX,TW)[M] = \bfc_1(E^d)[C^d],
\end{equation}
here $E=u^*TX$.
\end{prop}
\proof  Since $\db M$ has only finitely many connected components we can assume
that $W$ has only finitely many connected components too. Therefore  the previous 
considerations apply where $W$ plays the role of $Y$. Consider the bundles $E\deff u^*TX$
and $F\deff u^*TW$ over $C$ and $\db C$ respectively. 

\smallskip Contract each boundary circle $\gamma_i$ of $C$ to a separate point  $y_i'$. 
Let $C'$ be the closed oriented surface obtained in this way and $q : C\to C'$ the contraction 
map. For a complex vector bundle $E$ over $C$ and a totally real subbundle $F$ of the restriction 
$E|_{\db C}$ the line bundle $\det^2(E)$ admits a distinguished homotopy class of
trivi\-ali\-sations over each boundary circle $\gamma_i$ induced by the bundle $\det^2(F)$. This 
allows us to construct a line bundle $L$ over $C'$ with an isomorphism  $\det^2(E) \cong q^*L$ 
compatible with the trivialisations of $\det^2(E)$ over each boundary circle $\gamma_i$ induced 
by $\det^2(F)$. Then the Maslov class $\bfmu(E,F)\in \sfh^2(C, \db C; \zz)$ is equal to the 
pull-back of the Chern class $c_1(L)$ wrt.\ the map $q :C \to C'$. 

\smallskip Let $C''$ be a copy of $C'$ with the reversed orientation. By $y'_i\in C'$ 
we denoted the point on $C'$ arising from the contraction of a boundary circle $\gamma_i$ of $C$. 
Denote by  $y''_i$ the corresponding point on $C''$. Let $C^\times$ 
be the topological space obtained from the disjoint union $C' \sqcup C''$ by identifying each point 
$y'_i$ on $C'$ with the corresponding point $y''_i$ on $C'$,  and denote by $y^\times_i$ the arising 
point on $C^\times$. Then $C^\times$ is a {\slsf nodal surface} and each $y^\times_i$ is a nodal 
point. Notice that we can obtain  $C^\times$ from the Schottky double $C^d$ 
contracting each circle $\gamma_i \subset C^d$ to the point $y^\times_i$. We denote this contraction
map by $q^\times : C^d \to C^\times$.

\smallskip  Further, we denote by $\tau^\times: C^\times \to  C^\times$ the 
canonical involution interchanging $C'$ with $C''$. It is not difficult to see 
that our construction of the Maslov cohomology class $\bfmu(E,F)\in \sfh^2(C, \db C; \zz)$ 
is compatible with the doubling construction and the involutions. In particular, the
$\tau_E$-image of the bundle $\det^2(E)$ over $C$  and the $\tau^\times$-image of
the line bundle $L'$ over $C'$ will be the bundle $\det^2(E^-)$ over $C^-$  and the line 
$L''$ bundle over $C''$. $L'$ and $L''$ naturally glue to a line bundle $L^\times$
over $C^\times$, and the latter admits an anti-holomorphic involution $\tau_L$ 
which sends $L'$ to $L''$. The totally real subbundle $F$ is $\tau_E$-invariant, 
the involution $\tau_L : L' \to L''$ is anti-linear, and we still have the isomorphism
$\det^2(E^-) \cong q^*L''$. This gives us a canonical isomorphism $\det^2(E^d) \cong 
(q^\times)^*L^\times$ where $q^\times : C^d \to C^\times$ is the contraction as map above. 
Therefore we obtain
\begin{equation}
\eqqno(tau-class)
\tau^*\bfmu(E,F) = -\bfmu(E^-,F), 
\end{equation}
where $\bfmu(E,F)\in \sfh^2(C, \db C; \zz)$ and $\bfmu(E^-,F)\in \sfh^2(C^-, \db C^-; \zz)$
and $\tau \deff (\tau_S, \tau_E)$. Notice that we have the sign reversion on the level
of \emph{cohomology classes}, but the equality
\begin{equation}
\eqqno(tau-index)
\bfmu(E,F)\cdot [C] = \bfmu(E^-,F)\cdot [C^-]
\end{equation}
for the Maslov \emph{indices}. Formula \eqqref(tau-class) holds true because $E^-$ has 
the opposite complex structure, and \eqqref(tau-index) because of the reversion of the 
orientation. Further, the standard properties of the Chern class give us the equality
\[
\begin{split}
2 \bfmu(E,F)\cdot [C]\  &= \bfmu(E,F)\cdot [C] + \bfmu(E^-,F)\cdot [C^-]
\ = c_1\big( \det^2(E^d) \big) \cdot [C^d]
\\[3pt]
&= 2 c_1\big( \det(E^d)\big ) \cdot [S^d]
\ = 2 c_1(E^d) \cdot [S^d].
\end{split}
\]
Therefore we conclude the following  formula for the Maslov index
\begin{equation}
\bfmu(E,F)[C] = c_1(E^d) [C^d],
\end{equation}
which is \eqqref(masl-cl-ind). Proposition is proved.

\smallskip\qed

\begin{corol}
\label{masl-cl-inv}
Let a pair of almost complex surface and a totally real submanifolds $W$ be as 
above and let $M$ be a complex curve attached to $W$ without cusps on $W$. 
Let $\widetilde{M}^\delta$ be the result of surgery as above. Then 
\begin{equation}
\eqqno(ind-surg)
\bfmu_{\widetilde{M}^\delta}(TX,TW) = \bfmu_M(TX,TW).
\end{equation}
\end{corol}
\proof From \eqqref(mas-dbl) and \eqqref(masl-cl-ind) we see that the Maslov index
both on the left and right hand side of \eqqref(ind-surg) is equal to the value 
of the Malov class $\bfmu (TX,TW)$ on $[\tilde M^\delta]$ and  $[M]$ respectively. 
But since $(\tilde M^\delta ,\db \tilde M^\delta)$ is homotopic to $(M,\db M)$ 
(just let $\delta\searrow 0$) they define the same homology class in $\sfh_2(X,W;\zz)$.
Corollary is proved.

\smallskip\qed

\newprg[ADJ.proof]{Proof of the Adjunction Formula}

Let us start with the following preliminary statement assuming that 
$M$ is embedded, \ie $\varkappa^{(i)}=\delta^{(i)} = \delta^{(b)} = 0$.

\begin{thm} {\slsf (Adjunction Formula - I).} 
\label{adj-2} 
Let $M$ be a smooth connected compact $J$-complex curve in $X$ 
with boundary on $W$  of genus $g$ with $\sigma$ boundary circles. Then 
\begin{equation}
\eqqno(adj-frm2)
2g + \sigma = \frac{[M^d]^2- \mu_{\db M} (TX,TW)}{2} + 2. 
\end{equation}
Here  $\mu_{\db M} (TX,TW)$ is the Maslov index of $TW$ along 
$\db M$, and $[M^d]^2$ is the homological self-intersection of the Schottky 
double of $M$.
\end{thm}
\proof  Since $TX|_M = TM\oplus N_M$ and $TW|_{\db M} = T\db C\oplus
N^W_{\db M}$ we have that 
\begin{equation}
\eqqno(sum-mas)
\mu_{\db M} (TX,TW) = \mu (TM,T\db M) + \mu (N_M,N^W_{\db M}).
\end{equation}
Therefore we derive from \eqqref(mas-tang)  that 
\[
\mu_{\db M} (TX,TW) = 4-4g-2\sigma + \mu (N_M,N^W_{\db M}).
\]
By Proposition \ref{eps-mas}  we have that $\mu (N_M,N^W_{\db M})
= [M^d]^2$ for embedded $M$. Consequently  
\[
\mu_{\db M} (TX,TW) = 4-4g-2\sigma + [M^d]^2.
\]
And this is equivalent to \eqqref(adj-frm2).

\smallskip\qed 

\smallskip 
Now let us remove the condition on $M$ to be embedded. But we shall assume
that $M$ has no cusps on the boundary. That is we are going to prove Theorem
\ref{adj-1} from the Introduction.

\smallskip Perturb $M$ in a neighbourhood of its singular points (also 
on $W$) in order to have 
only double  transverse self intersections. On points in $\db M$ 
the perturbation is made as in Theorem \ref{b-pos-thm} (\sliii, \ie 
keeping all intersections on $W$.  As the result we obtain 
a $\ti J$-complex curve $\ti M$ (with $\ti J$ being a perturbation of $J$) 
such that: 

\begin{itemize}
\item the sum of the boundary intersection indices stays the same, but
all intersections are double and transverse;

\item the sum $\delta^{(i)} + \varkappa^{(i)}$ for $M$ is equal 
to $\delta^{(i)}$ for $\widetilde{M}$ and again, all intersections
are double and transverse.
\end{itemize}

Denote $\ti M$ and $\ti J$ still as $M$ and $J$. Homological self-intersection 
and Maslov index $\mu_{\db M}(TX,TW)$ $[M^d]^2$ do not change. Therefore the 
formula \eqqref(adj-frm1) stays intact, only $\delta^{(i)} + \varkappa^{(i)}$ will 
be replaced by $\delta^{(i)}$. 
What is left to prove is 
\begin{equation}
\eqqno(adj-3)
2g + \sigma = \frac{[M^d]^2- \mu_{\db M} (TX,TW)}{2} + 2 - 
\delta^{(b)} - 2\delta^{(i)},
\end{equation}
where $M$ is now immersed and has only double transverse self-intersections,
both inner and boundary. Replacing $[M^d]^2$ by the normal Maslov index plus
$4\delta^{(i)} $, as in  Proposition \ref{eps-mas}, we are left with   
\begin{equation}
\eqqno(adj-4)
2g + \sigma = \frac{\mu (N_M,N^W_{\db M})- \mu_{\db M} (TX,TW)}{2} + 2 -
\delta^{(b)}
\end{equation}
to prove. Replace every inner self-intersection point by a handle, 
\ie the cross $\{z_1\cdot z_2 = 0\}$ replace by $\{z_1\cdot z_2 = \delta\}$. 
Maslov index and homological self-intersection will not change again. Indeed, 
let $C$ be the curve parametrising $M$ and $\ti C$ the curve parametrising 
$\ti M$. Then $C $ is obtained from $\ti C$ contracting certain circles, \ie 
there exists a projection $\sigma :\ti C \to C$ contracting circles to nodes. 
Mapping $u:C\to X$, which parametrises $C$, is homotopic to $u\circ \pi :
\ti C\to X$. Due to the homotopy invariance of Maslov index we get the needed
assertion.

\smallskip Genus will increase by 
$1$. As for the normal Maslov index remark that the vector field $\nu =
(- \bar z_2, \bar z_1)$, which stays normal to $M$ along the perturbation
has two negative zeroes on the cross, and no zeroes on the handle. 
I.e., $\mu (N_M,N^W_{\db M})$ increases by four and the formula \eqqref(adj-3)
stays intact.

\smallskip Therefore we need to prove \eqqref(adj-3) for $M$ embedded to $X\setminus W$
with at most 
transverse self-intersections on $W$. Again by Proposition \ref{eps-mas}
in the case of absence of inner self-intersections we have that $\mu (N_M,N^W_{\db M})
=[M^d]^2$ and, since $\delta^{(b)} + 2g + \sigma -2 = -\chi (M)$, we rewrite \eqqref(adj-3)
as 
\begin{equation}
\eqqno(adj-5)
-\chi (M) = \frac{[M^d]^2- \mu_{\db M} (TX,TW)}{2}.
\end{equation}

In order to prove this last formula perform the surgery as described before the 
Proposition \ref{eps-surg} at every point of boundary self-intersection. The 
homological self-intersection of the double will not change due to item (\slii 
of this proposition. Maslov index doesn't change as well as the Euler 
characteristic, because the surface obtained by the surgery 
retracts to the initial one. We reduced our task to the embedded case, which is
served by formula \eqqref(adj-frm2). Theorem is proved.

\smallskip\qed

\newprg[ADJ.cplx]{Almost Complex Structure on the Double}
It is worth to notice that  $X^d_{\eps}$ admits an almost complex structure
keeping $M_\eps^d$ complex. Let us show this even if we didn't use this 
in this paper.

\begin{prop}
The $\eps$-double $X^d_{\eps}$  of an almost complex surface  $(X,J)$ 
relative to the totally real surface $W$ admits an almost complex 
structure $J^d_{\eps}$, which is equal to $J$ on $X^+_{\eps}$. Moreover, 
let $M$ is a $J$-complex curve embedded to $X\setminus W$, attached by 
with boundary $\db M$ to $W$ and having at most self-intersections at 
points of $W$. Then $J^d_\eps$ can be constructed such that the double 
$M^d_{\eps}$ will be $J^d_{\eps}$-complex.
\end{prop}
\proof On $X^+_{\eps}$ we do not change the almost complex structure. 
On the layer $B^-_{\eps , 2\eps}$ in $X^-_{\eps}$ change the structure $J$ 
as follows. Using the fact that $B^-_{\eps , 2\eps}\deff \bigcup_{\eps <\delta 
<2\eps}\db B_\delta$ we take for every $x\in B_{\eps , 2\eps}$ he unique
$\delta$ such that $x\in \db B_\delta$ and set $J^-= J$ on the tangent to 
$\db B_{\delta}$ complex space through $x$ and $J^-=-J$ on the orthogonal 
to $\db B_\delta$ line. Notice that $J^-$ matches with $J$ on $\db B_\eps$ 
in $X^d_{\eps}$ and therefore we 
got an almost complex structure in an $\eps$-neighbourhood $B^-_{\eps ,2\eps}
\cup \db B_{\eps} \cup B^+_{\eps , 2\eps}$ of $\db B_\eps$ in $X^d_{\eps}$ 
extending $J$. 

\medskip Denote by $X^+$ the manifold 
$X\setminus W$ and by $X^-$ the
same manifold with the opposite orientation, all this agrees with $X^{\pm}_{\eps}$ 
in the sense that $X^{\pm}_{\eps}$ is naturally a domain in $X^{\pm}$. 
Take some Morse exhaustion $\rho$ of $X^+=X\setminus W$, having 
$\db B_{\delta}$ as its level sets for $0<\delta <3\eps$. If $M$ is smooth 
then we can make $\rho|_M$ Morse as well.  Associated with the $h$ Riemannian metric
we denote $g$, \ie $g = \re h$. This metric 
naturally extends to $X^-$. Take the same exhaustion function $\rho$ on 
$X^-$ and define an almost complex structure on $X^-$ as follows: for a 
regular point $x\in \{\rho = c\}$ in $X^-$ set $J^- = - J$ on the orthogonal 
to $\{\rho = c\}$ $J$-complex line (notice that this makes sense, since 
$J$ exists also on $X^-$), and $J^- = J$ on the tangent. This defines us 
an almost complex structure on $X\-\setminus \{\text{ critical points of } 
\rho \}$, and this structure  gives $X^-$ its orientation. Perturb $J^-$ 
near every regular point $\rho|_M$ in such a way that $J^-|_{M^-}
= - J|_{M^-}$. We get an almost complex structure on $X^d_{\eps}\bs 
\{\text{ finite set in } X^-_{\eps}\}$.

\smallskip Take any point $x$ in this finite set. Denote by $F^{\pm}_r$ 
the distribution of $J$-complex planes on the sphere $\ss^3_r$ such
that $J^-|_{F^{\pm}_r} = \pm J|_{F^{\pm}_r}$. Notice that $F^{\pm}$
are well defined ($J^-$ cannot be equal to $J$ on two distinct complex
planes). Set $F^{\pm}_r = F^{\pm}|_{\ss^3_r}$. Since these distributions 
are trivial we can deform them when passing from some fixed $r$ to 
$r/2$ to the constant ones, say $F^-|_{r/2} = \{z_2 = 0\}$ and 
$F^+|_{r/2} = \{z_1 = 0\}$. If our point $x$ belongs to $M$ then 
find coordinates in a neighbourhood of it such that $M = \{z_2=0\}$.
Now extend $J^-$ to the ball $B_{r/2}(x)$ as equal to $\pm J$ on 
$F^{\pm}$. 

\smallskip\qed

\newsect[CSP]{Smoothing of Boundary Cusps}

\newprg[CSP.pert]{Higher order estimates}

For the proof of the Theorem \ref{b-pert-thm} we shall need a stronger version of Lemma \ref{disk-pert-r},
which includes the $L^{2,p}$-estimates on the place of $L^{1,p}$-ones. Let us start with an another version 
of Proposition \ref{morrey1}. 





\begin{prop}
\label{morrey2}
Under the conditions of Proposition \ref{morrey1} assume that $J$ is
Lipschitz-continuous, $R\in L^{1,p}(\Delta)$, and the norms $\norm{J-J\st}_{\calc^{Lip}
(\Delta)}$, $\norm{R}_{L^{1,p}(\Delta)}$ are sufficiently small. Then there 
exists  a bounded linear operator $\wt T_{J,R}^{\,0}: L^{1,p}(\Delta , \rr^{2n})\to
L^{2,p}(\Delta , \rr^{2n})$ such that:

\sli $\wt T_{J,R}^{\,0}(u)(0) = 0$ and $\wt T_{J,R}^{\,0}$ is the right inverse 
to $\dbar_J + R$, \ie $ (\dbar_J + R)\circ \wt T_{J,R}^{\,0}\equiv \id$;

\slii it satisfies $\norm{\wt T_{J,R}^{\,0}u}_{L^{1,p}(\Delta)} \le C_p\, 
\norm{u}_{L^p(\Delta)}$ for $u\in L^{1,p}(\Delta)$;

\sliii if, moreover, $J(\zeta)=J\st$ and $R(\zeta)$ is real for real $\zeta$ 
then $\wt T_{J,R}^{\,0}$  is real provided $J$ 

\quad satisfies \eqqref(im-j).
\end{prop}
\proof First let us construct the operator $\wt T_{J,R}^{\,0}$. denote by $L^{k,p}_0(\Delta ,
\cc^n)$ the closure of $\cald (\Delta ,\cc^n)$ in the spaces $L^{k,p}(\Delta, \cc^n)$. 
Here $\cald (\Omega ,\cc^n)$ is the space of $\calc^{\infty}$-maps $\Omega\to \cc^n$ 
with compact support contained in $\Delta$. We use the following construction from 
\cite{GT} to construct an {\slsf extension operator}
\begin{equation}
\eqqno(ext-oper)
\ext :   L^{1,p}(\Delta) \to L^{1,p}_0(\Delta (2)).
\end{equation}
Fix a point $z_0\in \db\Delta$ and take a diffeomorphism of a neighbourhood 
of $z_0$ to a neighbourhood $\Omega$ of zero, which sends locally $\db\Delta$ 
to $\rr$. Wlog we assume that $\Omega = (-\delta , \delta)\times (-\eps, \eps)$ 
is a rectangle. Denote $\Omega^{\pm} = \Omega\cap \Delta^{\pm}$. Define an
extension of an $L^{1,p}$-function $u$ from $\Omega^+$ to $\Omega$ as follows:
\begin{equation}
\eqqno(l1p-ext)
\tilde u(x,y) = 
\begin{cases}
u(x,y) \text{ for } y\ge 0,\cr
-3u(x,-y) + 4u\left(x,-\frac{y}{2}\right) \text{ for } y<0.
\end{cases}
\end{equation}
It is not difficult to see, see Theorem 7.25 in \cite{GT}, that for $u\in L^{1,p}(\Omega^+)$ 
one has $\norm{\tilde u}_{L^{1,p}(\Omega)} \le C\cdot \norm{u}_{L^{1,p}(\Omega^+)}$.
Then glue these extensions together using a partition of unity, see Lemma 6.37 there
and  gets a bounded operator as in \eqqref(ext-oper), \ie 
\begin{equation}
\eqqno(ext-oper-l1p)    
\norm{\tilde u}_{L^{1,p}(\Delta} \le C_p\norm{u}_{L^{1,p}(\Delta (2))}.
\end{equation}
Moreover, it is not difficult to see that for $u\in L^{1,p}(\Omega^+)$ one has 
also
\begin{equation}
\eqqno(ext-oper-lp)  
\norm{\tilde u}_{L^p(\Omega)} \le C\cdot \norm{u}_{L^p(\Omega^+)}.
\end{equation}

\begin{exmp} \rm The construction above works for any domain with sufficiently 
regular boundary, see \cite{GT}. In the case of a disc this construction can be
performed globally as follows: define 
\[
\hat u(r,\phi) = 
\begin{cases}
u(r\phi) \text{ for } r\le 1\cr
-3u(2-r,\phi) + 4u(1-r/2,\phi),
\end{cases}
\]
and set $\ti u (r,\phi)=\chi (r)\hat u(r,\phi)$, where $\chi$ is a bump-function equal to 
$1$ for  $r\le 1+\eps$ and $0$ for $r\ge 1+2\eps$. 
\end{exmp}

\begin{defi}
\label{wt-CG-def}
We define the extension of the Cauchy-Green operator
\begin{equation}
\eqqno(CG-ext1)
\wt T_{CG}: L^{1,p}(\Delta) \to L^{2,p}(\Delta),
\end{equation}
as $\wt T_{CG} \deff T_{CG}\circ \ext$ restricted back to $\Delta$.
\end{defi}
Notice that from \eqqref(ext-oper-lp)  we see that $\widetilde{T}_{CG}$ extends 
to a bounded linear operator $L^p(\Delta) \to L^{1,p}(\Delta )$. 
Set $\wt T_{J\st,0}^{\,0}(u) := \wt T_{CG}u - (\wt T_{CG}u)(0)$ in order 
to have that  $\wt T_{J,R}^{\,0}(u)(0)=0$. This gives us the construction 
of $\tilde T^0_{CG}$ as the right inverse to $\dbar_{J\st}$ both on 
$L^{1,p}(\Delta)$ and $L^p(\Delta)$. Moreover, $\wt T^0_{CG}$ is obviously
real.

\smallskip 
The rest of the proof here is essentially the same as in Proposition \ref{morrey1} 
after we replace 
the standard Cauchy-Green operator $T^0_{CG}$ by its extension $\wt T^0_{CG}$ 
just constructed. For the general 
$J$ and $R$ operator $\wt T_{J,R}^{\,0}: L^{1,p}(\Delta)\to L^{2,p}(\Delta)$ is constructed 
using the same series \eqqref(pert-op) with the operator $T^0_{J\st ,0}$ replaced 
by $\wt T^{\,0}_{J\st,0}$. In this case a sufficient condition for the convergence is 
that $\norm{J-J\st}_{\calc^{Lip}(\Delta)} + \norm{R}_{L^{1,p}(Delta)}$ is small 
enough, which is supposed in the hypotheses of the proposition.

\smallskip\qed

\begin{rema} \rm 
\label{tildet-ext)}
The estimate \eqqref(ext-oper-lp) means that $\wt T_{J,R}^0$ extends by continuity 
to a bounded linear operator from the
dense subspace $L^{1,p}(\Delta)$ of $L^p(\Delta)$ to the whole of $L^p(\Delta)$.
\end{rema}

\medskip 
Pushing forward the constructions from the proof of  Lemma \ref{disk-pert-r} with the 
operator $T^0_{\tilde J^{(\nu)},R^{(\nu)}}$ from Proposition \ref{morrey1} replaced by 
the operator $\wt T^0_{\tilde J^{(\nu)},R^{(\nu)}}$ from the Proposition \ref{morrey2} 
we can improve the statement in the case when the initial map $u_0$ has a 
cusp at the origin, \ie $\mu \ge 2$ or, if $\mu \ge 1$ but $\nu =0$, and, moreover 
$J$ is of class $\calc^{1,\alpha}$. In the following lemma we assume that a 
$J$-holomorphic map $u_0:(\Delta^+,\db_0\Delta^+,0)\to (\rr^{2n},\rr^n,0)$ takes its
values in the unit ball $B$ of $\rr^{2n}$. Moreover, modulo a dilatation $J_{\delta}(z)
\deff J(\delta^{\mu}z)$, $u_{\delta}(\zeta) \deff \frac{1}{\delta^{\mu}}u_0(\delta \zeta)$
we can assume that $\calc^{1,\alpha}$-norm of $J$ in $B$ is as small as we need.

\begin{lem}
\label{cusp-pert-h} 
Let $J$ be an almost complex structure of class $\calc^{1,\alpha}$ on $\rr^{2n}$ such
that $J|_{\rr^n} = J\st$ and let $u_0:(\Delta^+, \db_0\Delta^+, 0) \to (\rr^{2n}, \rr^n,0)$
be a $J$-holomorphic map. Suppose  that 
\begin{equation}
\eqqno(mu-nu-alpha)
2\mu - 2 + (\alpha -1)\nu\ge 0 
\end{equation}
where $\mu$ is
the vanishing order of the mapping $u_0$ at the origin. Then  there exists $w:(\Delta^+,\db_0\Delta^+,0)\to (\rr^{2n}, \rr^n,0)$, which belongs to $L^{2,p}(\Delta )$ for all 
$p>2$, such that 
\begin{equation}
u(\zeta) = u_0(\zeta) + \zeta^{\nu}w(\zeta)
\end{equation}
is $J$-holomorphic and satisfies the estimate 
\begin{equation}
\eqqno(pert5)
\norm{w}_{L^{2,p}(\Delta^+)} \leq C_p\,|\wect_0| + C_p\norm{J}_{\calc^{1,\alpha}}
|\wect_0|^{\alpha}.
\end{equation}
\end{lem}

\proof We continue the constructions from the proof of Lemma  \ref{disk-pert-r}, 
with the difference that in the appropriate places $L^p$-norms will be replaced by 
$L^{1,p}$-ones and $L^{1,p}$-norms by $L^{2,p}$-ones. We keep the 
notations from that proof. In particular we first extend $u_0$ (in fact the section 
$u_0:\Delta^+\to E^+$) to the section $\tilde u_0:\Delta\to E$ by reflection. 
\begin{rema} \rm 
{\slsf a)} Take  $\tilde w$ which satisfies the reality condition 
$\overline{\tilde w(\bar \zeta)}= \tilde w(\zeta)$ and set $w=\tilde w|_{\Delta^+}$.
Then for $J(u_0 + \zeta^{\nu}w)$ extended to $E$ by reflection, \ie for $\tilde J(u_0 + 
\zeta^{\nu}w)$, we use the notation $J(\tilde u_0 + \zeta^{\nu}\tilde w)$. 

\smallskip\noindent{\slsf b)} Notice that for such $\ti w$ one has 
$\norm{\ti w}_{L^{1,p}(\Delta)} = 2\norm{w}_{L^{1,p}(\Delta^+)}$. As well as 
for $w\in L^{1,p}(\Delta^+)$, which is real valued for real $\zeta$, its extension 
to $\Delta$ by reflection $\ti w$ also belongs to $L^{1,p}(\Delta)$, provided $p>1$.

\smallskip\noindent{\slsf c)} Finally, that we can differentiate $J(\tilde u_0 + 
\zeta^{\nu}\tilde w)$ as follows 
\begin{equation}
\eqqno(diff-j-u)
\db_{\xi}J(\tilde u_0(\zeta) + \zeta^{\nu}\tilde w(\zeta)) =
\begin{cases}
 \nabla_zJ(u_0(\zeta) + \zeta^{\nu}w(\zeta))\left(\db_{\xi}u_0(\zeta) + 
 \db_{\xi}(\zeta^{\nu}w(\zeta) \right)    \,\, +\cr
 + \, \nabla_{\bar z}J(u_0(\zeta) + \zeta^{\nu}w(\zeta))\left(\db_{\xi}\bar u_0(\zeta)
 + \db_{\xi}(\bar\zeta^{\nu} \bar w(\zeta)\right)   \quad\text{ if } \im\zeta \ge 0,\\[6pt]
 - \nabla_z\overline{J(u_0 (\bar\zeta)+ \bar \zeta^{\nu}w(\bar\zeta))}
 \left(\overline{\db_{\xi}u_0(\bar\zeta)} + \overline{\db_{\xi}(\bar\zeta^{\nu}
 w(\bar\zeta)} \right)    \,\, +\cr
 + \, \nabla_{\bar z}\overline{J(u_0(\bar\zeta) + \bar\zeta^{\nu}w(\bar\zeta))}
 \left(\overline{\db_{\xi}\bar u_0(\bar\zeta)} + \overline{\db_{\xi}(\bar\zeta^{\nu}
 \bar w(\zeta)}\right)   \quad\text{ if } \im\zeta \le 0,
\end{cases}
\end{equation}
and the same for $\db_{\eta}$. One the real axis the derivative $\db_{\eta}$
is taken on each side separately. Since we need the $L^p$-norm of these derivatives
only, this will be sufficient.

\end{rema}

Therefore it will be sufficient to perform all computations for $\im\zeta\ge 0$. 
We will further study the expressions $R^{(\nu)}$ and $F^{(\nu)}$ defined in 
\eqqref(term-r) and \eqqref(term-f). First of all we need to establish an 
$L^{1,p}$-estimates analogues to $L^p$-ones as in \eqqref(est2-r), \eqqref(est1-r)
and \eqqref(est3-r) for $R^{(\nu)}$ and $F^{(\nu)}$.  

\smallskip\noindent{\slsf Step 5.} 
{\it $L^{1,p}$-estimate of $R^{(\nu)} = \nu \zeta^{-\nu}\big(\tilde J_{u_0} - 
J\st\big)J\st \zeta^{\nu -1}$.}  If $\nu =0$ we see that $R^{(0)}\equiv 0$ and 
there is nothing to prove. Otherwise, we see from \eqqref(mu-nu-alpha) that 
$\mu $ must be $\ge 2$. Now from \eqqref(est-A2) we get the following estimate of the 
Lipschitz constant of $R^{(\nu)}$:

\begin{equation}
\eqqno(e-r-nu2a)
Lip_{\Delta }(R^{(\nu)})\leq C\,Lip_B(J)\,\norm{\tilde u_0}_{\calc^{lip}(\Delta)}.
\end{equation}
Further, using the obvious pointwise estimate
\begin{equation*}
\big|R^{(\nu)}(\zeta) \big| 
\le \nu\cdot |\zeta|^{-1} \cdot Lip_B(J) \cdot |\tilde u_0(\zeta)|
\le \nu\cdot |\zeta|^{-1} \cdot Lip_B(J) \cdot |\zeta|^\mu \cdot 
    \norm{\tilde u_0}_{\calc^{lip}(\Delta)},
\end{equation*}
we conclude that the matrix valued function $R^{(\nu)}(\zeta)$ vanishes at $\zeta=0$.
Therefore we obtain the estimate
\begin{equation}
\eqqno(e-r-nu3)
\norm{R^{(\nu)}}_{\calc^{lip}(\Delta)} = \norm{R^{(\nu)}}_{L^{1,\infty}(\Delta)}\leq C\, Lip_B(J)\,
\norm{\tilde u_0}_{\calc^{lip}(\Delta)}.
\end{equation}
Notice that \eqqref(e-r-nu3)  implies for all $p>2$ the  estimate
\begin{equation}
\eqqno(e-r-nu4)
\norm{R^{(\nu)}\tilde w}_{L^{1,p}(\Delta )} \leq C\,Lip_B(J)
\norm{\tilde u_0}_{\calc^{lip}(\Delta)}\norm{\tilde w}_{L^{1,p}(\Delta)},
\end{equation}
for any $\tilde w\in L^{1,p}(\Delta , \cc^n)$. 

\medskip\noindent{\slsf Step 6.}  
{\sl $L^{1,p}$-estimates of $F^{(\nu)}(\zeta,w) = \zeta^{-\nu} \big(\tilde J_{u_0} - 
\tilde J_{u_0 +\zeta^{\nu}w}\big)\big( \db_{\eta}\tilde u_0 + \db_{\eta}(\zeta^{\nu}
\tilde w)\big)$. We shall prove the following estimates under the assumption that 
$2\mu -2 + (\alpha -1)\nu \ge 0$:}
\begin{equation}
\eqqno(est2)
\norm{F^{(\nu)}(\zeta, \tilde w)}_{L^{1,p}(\Delta)}\leq C\cdot \norm{J}_{\calc^{1,\alpha}
(\Delta)}\; \norm{\tilde u_0}_{L^{2,p}(\Delta)}\left(\norm{\tilde w}_{L^{2,p}(\Delta)}
+ \norm{\tilde w}_{L^p(\Delta)}^{\alpha}\right),
\end{equation}
and 
\begin{equation}
\eqqno(est4)
\begin{split}
& \norm{F^{(\nu)}(\zeta, \tilde w_1) - F^{(\nu)}(\zeta,\tilde w_2) }_{L^{1,p}(\Delta)}
\leq  
\\[2pt]
& \le \norm{J}_{\calc^{1,\alpha}(\Delta)}\; \norm{\tilde u_0}_{L^{2,p} \Delta)}
\left(\norm{\tilde w_1 -\tilde w_2}_{L^{2,p}(\Delta)} + \norm{\tilde w_1 - \tilde w_2}_{L^p(\Delta)}^{\alpha}\right).
\end{split}
\end{equation}
Again, as in the proof of Lemma \ref{disk-pert-r} we need to prove only the second estimate.
First, we decompose the difference $F^{(\nu)}(\zeta,\tilde w_1) - F^{(\nu)}(\zeta,
\tilde w_2)$ into the terms which will be estimated separately. As before we use the
notation $\tilde u_1(\zeta) := \tilde u_0(\zeta) + \zeta^\nu \tilde w_1(\zeta)$ and 
$\tilde u_2(\zeta) := \tilde u_0(\zeta) + \zeta^\nu \tilde w_2(\zeta)$. Rewrite
\eqqref(f-f1) as 
\begin{align}
F^{(\nu)}(\zeta,\tilde w_1)\, - &\,F^{(\nu)}(\zeta,\tilde w_2) =
\notag \\
\eqqno(I.1-I.2)
& = \zeta^{-\nu} \big(J(\tilde u_2) - J(\tilde u_1) \big)\, \db_\eta \tilde u_0 
\;+\; \zeta^{-\nu} J(\tilde u_0) \;\db_\eta \big( \zeta^\nu\,( \tilde w_1 -
\tilde w_2) \big) \;+\; 
\\
\eqqno(I.3-I.4)
&  + \zeta^{-\nu} \big(J(\tilde u_2) - J(\tilde u_1) \big)\, \db_\eta ( \zeta^\nu\, 
\tilde w_1)\;+\;  \zeta^{-\nu} \, J(\tilde u_2) \;\db_\eta \big( \zeta^\nu\,
( \tilde w_2 -\tilde w_1) \big).
\end{align}
We denote the terms in \eqqref(I.1-I.2) by $I_1$ and $I_2$, and in \eqqref(I.3-I.4) 
by $I_3$ and $I_4$. 

\medskip\noindent{\slsf Estimation of the term $I_1 =  \zeta^{-\nu}\big(J(\tilde u_2) - 
J(\tilde u_1) \big)\, \db_\eta \tilde u_0$.} The obvious pointwise estimate
\[
|I_1(\zeta)| \leq Lip_B(J)\cdot |\tilde w_2(\zeta) -\tilde w_1(\zeta)|\cdot 
|\db_{\eta}\tilde u_0|
\]
implies
\begin{equation*}
\norm{I_1}_{L^\infty(\Delta)} \le  Lip_B(J)\, 
\norm{\tilde w_2 -\tilde w_1}_{L^{\infty}(\Delta)} \, 
\norm{d \tilde u_0}_{L^\infty(\Delta)}.
\end{equation*}
Therefore, in order to estimate $L^{1,p}$-norm of $I_1$ its is sufficient 
to estimate the $L^p$-norm of its derivative. Write 
\[
\begin{split}
& \frac{\db I_1}{\db \xi} = -\nu \zeta^{-\nu -1} \big(J(\tilde u_0 + \zeta^{\nu}\tilde w_2) - J(\tilde u_0 + \zeta^{\nu}\tilde w_1)\big)\db_{\eta}
\tilde u_0 \, +
\\[2pt]
& + \zeta^{-\nu}\left[\nabla_zJ(\tilde u_2)\cdot\left(\frac{\db\tilde u_0}{\db \xi}
+\nu\zeta^{\nu -1}\tilde w_2\right)- \nabla_zJ(\tilde u_1)\cdot\left(\frac{\db
\tilde u_0}{\db \xi}+\nu\zeta^{\nu -1}\tilde w_1\right)\right]\db_{\eta}
\tilde u_0 \, + 
\\[2pt]
& + \zeta^{-\nu}\left[\nabla_zJ(\tilde u_2)\zeta^{\nu}\frac{\db\tilde w_2}{\db \xi}
- \nabla_zJ(\tilde u_1)\zeta^{\nu}\frac{\db\tilde w_1}{\db \xi}\right]\db_{\eta}
\tilde u_0 \, + 
\\[2pt]
& + \zeta^{-\nu}\left[\nabla_{\bar z}J(\tilde u_2)\cdot\left(\frac{\db
\overline{\tilde u_0}}{\db \xi} +\bar\zeta^{\nu}\frac{\db\overline{\tilde w_2}}
{\db \xi}\right)- \nabla_{\bar z}J(\tilde u_1)\cdot \left(\frac{\db
\overline{\tilde u_0}}{\db \xi} +\bar\zeta^{\nu}\frac{\db\overline{\tilde w_1}}
{\db \xi}\right)\right]\db_{\eta}
\tilde u_0 \, + 
\\[2pt]
& + \zeta^{-\nu} \big(J(\tilde u_2) - J(\tilde u_1) \big)\, \db^2_{\xi \eta} 
\tilde u_0 \quad =: \quad A_1+A_2+A_3+A_4+A_5.
\end{split}
\]
Now we can estimate the terms $A_k$. First, in the case $\nu =0$ we have that $A_1=0$.
For $\nu \ge 1$ write pointwisely
\[
\big|A_1\big| \le \nu |\zeta|^{-\nu -1}Lip_B(J)|\zeta|^{\nu}|\tilde w_2-\tilde w_1||\db_{\eta}
\tilde u_0 |. 
\]
Notice that in our settings $\tilde w_k(0) = \wect_0$ and therefore 
\[
\tilde w_k(\zeta) = \wect_0 + \int\limits_0^1\frac{d}{dt}\tilde w_k(t\zeta) dt = 
\wect_0 + \int\limits_0^1\left(\frac{\db \tilde w_k}{\db \zeta}\,\zeta + 
\frac{\db \tilde  w_k}{\db \bar\zeta}\,\bar\zeta\right)dt.
\]
This gives us 
\begin{equation}
\eqqno(w2-w1)
|\tilde w_2 - \tilde w_1| \le \int\limits_0^1\left(\Big\vert\frac{\db \tilde w_2}{\db \zeta} - 
\frac{\db \tilde w_1}{\db \zeta}\Big\vert + \Big\vert\frac{\db \tilde w_2}{\db \bar\zeta} - 
\frac{\db \tilde w_1}{\db \bar\zeta}\Big\vert\right)dt \cdot |\zeta| \le 
\norm{\tilde w_2-\tilde w_1}_{\calc^1(\Delta)}\cdot |\zeta|,
\end{equation}
which implies that 
\begin{equation}
\eqqno(est-a1)
\norm{A_1}_{L^p(\Delta)} \le \nu Lip_B(J)\, \norm{\tilde w_2 - \tilde w_1}_{L^{2,p}(\Delta)}
\, \norm{\tilde u_0}_{L^{1,p}(\Delta)}.
\end{equation}

The second term is treated as follows. The pointwise estimation yields
\[
\begin{split}
& \big|A_2\big| \le |\zeta|^{-\nu}\Big\vert \big[\nabla_zJ(\tilde u_2) - 
\nabla_zJ(\tilde u_1)\big] \Big\vert \Big\vert\frac{\db\tilde u_0}{\db \xi}\Big\vert 
|\db_{\eta}\tilde u_0 |\, +
\\[2pt]
& + |\zeta|^{-\nu}\nu \Big\vert \nabla_zJ(\tilde u_2)\zeta^{\nu-1}\tilde w_2 - 
\nabla_zJ(\tilde u_1)\zeta^{\nu-1}\tilde w_1\Big\vert |\db_{\eta}\tilde u_0 |\le
\\[2pt]
& \le |\zeta|^{-\nu }\norm{J}_{\calc^{1,\alpha}(B)} |\zeta|^{\nu\alpha}
|\tilde w_2-\tilde w_1|^{\alpha}|\zeta|^{2\mu -2}\norm{\tilde u_0}_{L^{1,p}(\Delta)}^2\, +
\\[2pt]
& + \nu|\zeta|^{-\nu}\Big\vert \nabla_zJ(\tilde u_2)\zeta^{\nu -1}(\tilde w_2 - \tilde w_1)
\Big\vert |\zeta|\norm{\tilde u_0}_{L^{1,p}(\Delta)} \, +
\\[2pt]
& + \nu|\zeta|^{-\nu}\Big\vert \big[\nabla_zJ(\tilde u_2) - \nabla_zJ(\tilde u_1)\big]
\zeta^{\nu -1}\tilde w_1\Big\vert |\zeta|\norm{\tilde u_0}_{L^{1,p}(\Delta)} \le
\\[2pt]
& \le C\cdot |\zeta|^{2\mu -2 +(\alpha -1)\nu}\norm{J}_{\calc^{1,\alpha}(B)}
|\tilde w_2-\tilde w_1|^{\alpha}\norm{\nabla\tilde u_0}_{L^{\infty}(\Delta)}^2 + 
\\[2pt]
& + C\cdot \left(\norm{J}_{\calc^1(B)} |\tilde w_2 - \tilde w_1| + 
|\tilde w_1|\norm{J}_{\calc^{1,\alpha}}|\tilde w_2 - \tilde w_1|^{\alpha}\right)
\norm{\tilde u_0}_{L^{1,p}(\Delta)} \le 
\\[2pt]
& \le C\cdot \norm{J}_{\calc^{1,\alpha}(B)}\norm{\tilde u_0}_{L^{2,p}(\Delta)}
\left(|\tilde w_2 - \tilde w_1| +  |\tilde w_2 - \tilde w_1|^{\alpha}\right),
\end{split}
\]
provided $\norm{\tilde u_0}_{L^{2,p}(\Delta)}\le 1$ and $2\mu -2 +(\alpha -1)\nu
\ge 0$.

\smallskip Notice the following: take $\frac{1}{q} = \frac{\alpha}{p}$ and $\frac{1}{r} = 
\frac{1-\alpha}{p}$ to have $\frac{1}{p} = \frac{1}{q} + \frac{1}{r}$ and, using the 
H\"older inequality, write 
\begin{equation}
\eqqno(r-p-q)
\Big\Vert\, |w_1 - w_2|^{\alpha}\Big\Vert_{L^p(\Delta)} \le  
\Big\Vert\, |w_1 - w_2|^{\alpha}\Big\Vert_{L^q(\Delta)} \cdot
\norm{\adyn}_{L^r(\Delta)} = C\, \Big\Vert\, w_1 - w_2
\Big\Vert_{L^p(\Delta)}^{\alpha}. 
\end{equation}
 As the result we obtain 
\begin{equation}
\eqqno(est-a2)
\norm{A_2}_{L^p(\Delta)} \le C\cdot \norm{J}_{\calc^{1,\alpha}(B)}\, 
\norm{\tilde u_0}_{L^{2,p}(\Delta)} \big(\norm{\tilde w_2 - \tilde w_1}_{L^p(\Delta)}
+ \norm{\tilde w_2 - \tilde w_1}_{L^p(\Delta)}^{\alpha}\big).
\, 
\end{equation}

\smallskip Third,
\[
\begin{split}
& \big\vert A_3\big\vert \le \left(\Big\vert \nabla_zJ (\tilde u_2) - 
\nabla J_z(\tilde u_1)\Big\vert \Big\vert\frac{\db \tilde w_2}{\db \xi}
\Big\vert + \norm{J}_{\calc^1(B)}\Big\vert \frac{\db \tilde w_2}{\db \xi}
- \frac{\db \tilde w_1}{\db \xi}\Big\vert\right)\norm{\tilde u_0}_{L^{1,p}(\Delta)} \le 
\\[2pt]
& \le C\cdot \norm{J}_{\calc^{1,\alpha}(B)} \left(|\tilde w_2-\tilde w_1|^{\alpha} + \Big\vert \frac{\db \tilde w_2}{\db \xi}- \frac{\db \tilde w_1}{\db \xi}\Big\vert\right)\norm{\tilde u_0}_{L^{1,p}(\Delta)}.
\end{split}
\]
This gives us 
\begin{equation}
\eqqno(est-a3)
\norm{A_3}_{L^p(\Delta)} \le C\cdot \norm{J}_{\calc^{1,\alpha}(B)}\, 
\norm{\tilde u_0}_{L^{1,p}(\Delta)} \big(\norm{\tilde w_2 - 
\tilde w_1}_{L^{1,p}(\Delta)}
+ \norm{\tilde w_2 - \tilde w_1}_{L^p(\Delta)}^{\alpha}\big).
\, 
\end{equation}

Fourth, 
\[
\begin{split}
&\big\vert A_4\big\vert \le  |\zeta |^{-\nu}\Big\vert \nabla_{\bar z}J
(\tilde u_2) -  \nabla_z J(\tilde u_1)\Big\vert \cdot \Big\vert
\frac{\db\tilde u_0}{\db \xi}\Big\vert |\db_{\eta}\tilde u_0 | + 
\\[2pt]
& + |\zeta|^{-\nu}
\Big\vert\nabla_{\bar z}J(\tilde u_2)\bar \zeta^{\nu}\frac{\db \overline{\tilde w_2}}
{\db \xi} - \nabla_{\bar z}J(\tilde u_1)\bar\zeta^{\nu}\frac{\overline{\tilde w_1}}
{\db \xi} \Big\vert |\db_{\eta}\tilde u_0 | \le
\\[2pt]
& \le |\zeta|^{2\mu - 2 + (\alpha -1)\nu}\norm{J}_{\calc^{1,\alpha}(B)}|\tilde w_2 - 
\tilde w_1|^{\alpha}\norm{\nabla\tilde u_0}_{L^{\infty}(\Delta)}^2 + 
\\[2pt]
& + \left(\Big\vert \nabla_{\bar z}J(\tilde u_2) -  \nabla_z J(\tilde u_1)\Big\vert 
\Big\vert \frac{\db \overline{\tilde w_2}}{\db \xi}\Big\vert
+ \Big\vert \nabla_z J(\tilde u_1)\Big\vert \Big\vert \frac{\db \overline{\tilde w_2}}
{\db \xi} - \frac{\db \overline{\tilde w_1}}{\db \xi}\Big\vert\right)
|\db_{\eta}\tilde u_0 | \le 
\\[2pt]
&  \le C\cdot \norm{\tilde u_0}_{L^{2,p}(\Delta)}\norm{J}_{\calc^{1,\alpha}(B)}
\left(|\tilde w_2 - \tilde w_1|^{\alpha} + |\nabla \tilde w_2 - \nabla\tilde w_1|
\right).
\end{split}
\]
Therefore we obtain 
\begin{equation}
\eqqno(est-a4)
\norm{A_4}_{L^p(\Delta)} \le C\cdot \norm{J}_{\calc^{1,\alpha}(B)}\, 
\norm{\tilde u_0}_{L^{2,p}(\Delta)} \big(\norm{\tilde w_2 - 
\tilde w_1}_{L^{1,p}(\Delta)}
+ \norm{\tilde w_2 - \tilde w_1}_{L^p(\Delta)}^{\alpha}\big),
\, 
\end{equation}
again, provided $\norm{\tilde u_0}_{L^{2,p}(\Delta)}\le 1$ and $2\mu -2 +(\alpha -1)\nu
\ge 0$.

\smallskip Finally,
\[
\begin{split}
\big\vert A_5\big\vert \le Lip_B(J)|\tilde w_2 - \tilde w_1||\db^2_{\zeta  \eta}
\tilde u_0|,
\end{split}
\]
and therefore 
\begin{equation}
\eqqno(est-a5)
\norm{A_5}_{L^p(\Delta)} \le C\cdot Lip_B(J) \norm{\tilde u_0}_{L^{2,p}(\Delta)}.
\norm{\tilde w_2 - \tilde w_1}_{L^p(\Delta)}.
\end{equation}
 
Estimation of $\frac{\db I_1}{\db \bar \zeta}$ is analogous (even simpler),
and therefore all together we get 
\begin{equation}
\eqqno(est-i1)
\norm{I_1}_{L^{1,p}(\Delta)} \le C\cdot \norm{J}_{\calc^{1,\alpha}(B)}\, 
\norm{\tilde u_0}_{L^{2,p}(\Delta)} \big(\norm{\tilde w_2 - 
\tilde w_1}_{L^{1,p}(\Delta)}
+ \norm{\tilde w_2 - \tilde w_1}_{L^p(\Delta)}^{\alpha}\big)
\end{equation}

\smallskip\noindent
{\slsf Estimation of the term $I_2 = \zeta^{-\nu} J(\tilde u_0) \;\db_\eta 
\big( \zeta^\nu\,( \tilde w_1 -\tilde w_2) \big)$.} Write 
\[
\begin{split}
& \frac{\db I_2}{\db \xi} = -\nu \zeta^{-\nu -1}J(\tilde u_0) \db_{\eta}
(\zeta^{\nu}(\tilde w_1 - \tilde w_2)) + \zeta^{-\nu}\nabla_zJ
\frac{\db \tilde u_0}{\db \xi}\db_{\eta}(\zeta^{\nu}(\tilde w_1 - \tilde w_2))
\\[2pt]
&+ \zeta^{-\nu}\nabla_{\bar z}J\frac{\db \overline{\tilde u_0}}{\db \xi}\db_{\eta}
(\zeta^{\nu}(\tilde w_1 - \tilde w_2)) + \zeta^{-\nu}J(\tilde u_0)
\nu \zeta^{\nu -1}\db_{\eta}(\tilde w_1 - \tilde w_2) + 
\\[2pt]
& + \zeta^{-\nu}J(\tilde u_0)\nu(\nu -1)i\zeta^{\nu -2}(\tilde w_1 -\tilde w_2) + 
\zeta^{-\nu}J(\tilde u_0) \zeta^{\nu}\db^2_{\zeta\eta}(\tilde w_1 - \tilde w_2) + 
\\[2pt]
& + \zeta^{-\nu}J(\tilde u_0)\nu i\zeta^{\nu -1}(\tilde w_1 -\tilde w_2)
= B_1 + B_2 + B_3 + B_4 + B_5 + B_6 + B_7.
\end{split}
\]

In the case $\nu =0$ we have $B_1=0$. Otherwise  write
\[
\begin{split}
& B_1 =  -\nu \zeta^{-\nu -1}J(\tilde u_0) \nu i \zeta^{\nu -1}
(\tilde w_1 - \tilde w_2) - \nu \zeta^{-\nu -1} J(\tilde u_0) 
\zeta^{\nu }\db_{\eta} (\tilde w_1 - \tilde w_2)
\end{split}
\]
Since $\mu\ge 2$ we obtain the needed estimate from Lemma \ref{lip-cont}
using \eqqref(w2-w1).

\smallskip Second, 
\[
\begin{split}
& \big\vert B_2 \big\vert \le \big\vert \norm{J}_{\calc^1(B)}
\norm{\tilde u_0}_{L^{2,p}(\Delta)}
\left(\Big\vert \frac{\db \tilde w_1}{\db \xi} - \frac{\db \tilde w_2}{\db \xi}
\Big\vert + |\tilde w_1 - \tilde w_2|\right),
\end{split}
\]
because \eqqref(nrm-frm3) and $\tilde u_0(\zeta) = \vect_0\zeta^{\mu} + \zeta^{2\mu -2}\zeta v(\zeta)$ with  $\mu \ge 2$.

\smallskip The third term can be treated exactly as  the second one. For the fourth
and the fifth we obviously by Lemma \ref{lip-cont} and condition $\mu \ge 2$ obtain
the estimate by $Lip_B(J)\norm{\tilde u_0}_{L^{2,p}(\Delta)}$. Therefore we conclude with
the estimate
\begin{equation}
\eqqno(est-i2)
\norm{I_2}_{L^{1,p}(\Delta)} \le C\norm{J}_{\calc^1(B)} 
\norm{\tilde w_2 - \tilde w_1}_{L^{1,p}(\Delta)}.
\end{equation}


\medskip\noindent{\slsf  Estimation of the term $I_3=\zeta^{-\nu} \big(J(\tilde u_2) - 
J(\tilde u_1) \big)\, \db_\eta ( \zeta^\nu\, \tilde w_1)$.}  
Write 
\[
\begin{split}
& \frac{\db I_3}{\db \xi} = -\nu \zeta^{-\nu-1}\big(J(\tilde u_2) - 
J(\tilde u_1) \big)\, \left(\nu i\zeta^{\nu -1}\tilde w_1 + \zeta^{\nu} 
\db_\eta \tilde w_1\right) \,  + 
\\[2pt]
& + \zeta^{-\nu}\Big[\nabla_zJ(\tilde u_2)\cdot\left(\frac{\db\tilde u_0}{\db \xi}
+\nu\zeta^{\nu -1}\tilde w_2\right)- \nabla_zJ(\tilde u_1)\cdot\left(\frac{\db
\tilde u_0}{\db \xi}+\nu\zeta^{\nu -1}\tilde w_1\right)\Big]\db_\eta ( \zeta^\nu\, 
\tilde w_1) \, + 
\\[2pt]
& + \zeta^{-\nu}\Big[\nabla_zJ(\tilde u_2)\zeta^{\nu}\frac{\db\tilde w_2}{\db \xi}
- \nabla_zJ(\tilde u_1)\zeta^{\nu}\frac{\db\tilde w_1}{\db \xi}\Big]
\db_\eta ( \zeta^\nu\, \tilde w_1) \, + 
\\[2pt]
& + \zeta^{-\nu}\Big[\nabla_{\bar z}J(\tilde u_2)\cdot\left(\frac{\db\bar u_0}{\db \xi}
+\nu\bar\zeta^{\nu -1}\bar w_2\right)- \nabla_{\bar z}J(\tilde u_1)\cdot\left(\frac{\db
\bar u_0}{\db \xi}+\nu\bar\zeta^{\nu -1}\bar w_1\right)\Big]\db_\eta ( \zeta^\nu\, 
\tilde w_1) \, + 
\\[2pt]
& + \zeta^{-\nu}\Big[\nabla_{\bar z}J(\tilde u_2)\bar\zeta^{\nu}\frac{\db\bar w_2}{\db \xi}
- \nabla_{\bar z}J(\tilde u_1)\bar\zeta^{\nu}\frac{\db\bar w_1}{\db \xi}\Big]
\db_\eta ( \zeta^\nu\, \tilde w_1) \, + 
\\[2pt]
& = C_1 + C_2 +C_3+C_4+C_5.
\end{split}
\]
In the case $\nu =0$ we have that $C_1=0$. Otherwise we have 
\[
 \big\vert C_1\big\vert \le C\cdot Lip_B(J)|\zeta|^{-2}|\zeta|^{\nu}|\tilde w_1 - \tilde w_2|
 \le C\cdot Lip_B(J)|\zeta|^{\nu -1} \norm{\tilde w_2-\tilde w_1}_{\calc^1(\Delta)},
\]
due to \eqqref(w2-w1). This gives us the estimate when $\nu \ge 1$. The second term can be 
estimated as follows 
\[
\big\vert C_2\big\vert \le |\zeta|^{-1}\norm{J}_{\calc^{1,\alpha}(B)}
|\zeta|^{\nu\alpha}|\tilde w_2 - \tilde w_1|^{\alpha}
\norm{\tilde u_0}_{\calc^1(\Delta)}\norm{\tilde w_1}_{\calc^1(\Delta)}, 
\]
by H\"older $\alpha$-regularity of $\nabla J$. This gives us the needed estimate since 
we have that $|\nabla u_0(\zeta)|\ = O(|\zeta|^{\mu -1})$ and $\mu \ge 2$. 
Indeed, $|\frac{\db u_0}{\db \zeta}|\le \mu |\zeta|^{\mu -1} + 
O(|\zeta|^{\mu -1 + \alpha})$ by \eqqref(nrm-frm3).
 
\smallskip Estimation of $C_3$ is obvious and $C_4$ resp. $C_5$ is the same as of
$C_2$ resp. $C_3$.

\medskip\noindent{\slsf Estimation of the term $I_4= \zeta^{-\nu} \, J(\tilde u_2) \;\db_\eta \big( \zeta^\nu\, ( \tilde w_2 -\tilde w_1) \big)$.} Case $\nu = 0 $ is obvious. Therefore
assume that $\nu \ge 1$ and write
\[
\begin{split}
& \frac{\db I_4}{\db \xi} = -\nu \zeta^{-\nu -1}J(\tilde u_2) \db_{\eta}
(\zeta^{\nu}(\tilde w_1 - \tilde w_2)) + \zeta^{-\nu}\nabla_zJ
\frac{\db \tilde u_2}{\db \xi}\db_{\eta}(\zeta^{\nu}(\tilde w_1 - \tilde w_2)) \, +
\\[2pt]
& + \zeta^{-\nu}\nabla_{\bar z}J\frac{\db \overline{\tilde u_2}}{\db \xi}\db_{\eta}
(\zeta^{\nu}(\tilde w_1 - \tilde w_2)) + \zeta^{-\nu}J(\tilde u_2)
\nu i\zeta^{\nu -1}\db_{\eta}(\tilde w_1 - \tilde w_2) + 
\\[2pt]
& + \zeta^{-\nu}J(\tilde u_2)\nu(\nu -1)i\zeta^{\nu -2}(\tilde w_1 -\tilde w_2) + 
\zeta^{-\nu}J(\tilde u_2) \zeta^{\nu}\db^2_{\zeta\eta}(\tilde w_1 - \tilde w_2) + 
\\[2pt]
& + \zeta^{-\nu}J(\tilde u_2)\nu i\zeta^{\nu -1}(\tilde w_1 -\tilde w_2)
= D_1 + D_2 + D_3 + D_4 + D_5 + D_6 + D_7.
\end{split}
\]
Estimations of all term except $B_5$ are analogous to the previous estimtions. What concerns $B_5$ just remark that it is zero when 
$\nu 1,2$, and therefore it needed to be estimated for $\nu\ge 2$
only. And this is obvious. Estimations of the derivatives with respect
to $\eta$ are analogous. Step is proved.

\bigskip 
Now we can repeat step by step the proof of Lemma \ref{disk-pert-r}, but estimating
$L^{2,p}$-norms instead of $L^{1,p}$ ones. As in the quoted proof  we 
intend to solve the same initial value problem
\begin{equation}
\eqqno(syst1-1)
\begin{cases}
D_{J^{(\nu)}_{u_0},R^{(\nu)}}\tilde w = F^{(\nu)}(\zeta,\tilde w),\cr
\tilde w(0) = \wect_0,
\end{cases}
\end{equation}


\smallskip\noindent
by setting
\begin{equation} 
\eqqno(def-w-1-p)
\tilde w_1(\zeta) = \wect_0 - \wt T_{J^{(\nu)}_{u_0},R^{(\nu)}}^0 \big(D_{J^{(\nu)}_{u_0},R^{(\nu)}}
\wect_0\big), 
\end{equation}
where $\wect_0$ is considered as a constant function,  and then recursively for 
$n=1,2,3,\ldots$
\begin{equation}
\eqqno(newton-1-p)
\tilde w_{n+1} = \wt T_{J^{(\nu)}_{u_0},R^{(\nu)}}^0\big[F^{(\nu)}(\zeta,
\tilde w_n)\big] + \tilde w_1.
\end{equation}
Here $\wt T_{J^{(\nu)}_{u_0},R^{(\nu)}}^0$ is the right inverse to 
$D_{J^{(\nu)}_{u_0}, R^{(\nu)}}$ constructed in Proposition \ref{morrey2}. 
Notice that the functions $\tilde w_n(\zeta)$ satisfy the same relations as 
\eqqref(syst2), \ie 
\begin{equation}
\eqqno(syst-w_n-p)
\begin{cases}
D_{J^{(\nu)}_{u_0}, R^{(\nu)}}\tilde w_1 = 0, \cr
\tilde w_1(0) = \wect_0,
\end{cases}
\qquad
\begin{cases}
D_{J^{(\nu)}_{u_0}, R^{(\nu)}}\tilde w_n = F^{(\nu)}(\zeta, \tilde w_{n-1}), \cr
\tilde w_n(0) = \wect_0, \qquad \text{for }n\ge 2.
\end{cases}
\end{equation}
Notice that on $L^p(\Delta)$ operators $T_{CG}$  and $\tilde T_{CG}$ coincide and 
therefore $\tilde w_n$ converge to the solution of our system in $L^{1,p}$. Moreover,
they do satisfy the reality condition. Therefore the only thing we need to insure 
is the $L^{2,p}$-bound. Since $\tilde w_n$ are bounded 
in $L^{1,p}(\Delta , \cc^n)$ we remark that they do satisfy the a priori bound 
$\norm{\tilde w_n}_{L^{\infty}}\le \frac{1}{2}$ established in the proof of 
Lemma \ref{disk-pert-r}. 

\medskip\noindent{\slsf Step 7.} {\it There exists a constant $C$ independent 
of $n$ such that 
\begin{equation}
\eqqno(est-w-n)
\norm{\tilde w_n}_{L^{2,p}(\Delta)}\leq  C\,|\wect_0| + C\delta |\wect_0|^{\alpha},
\end{equation}
where $\delta = \norm{J}_{\calc^{1,\alpha}(\Delta)}$.
}

\smallskip The reasoning's here are analogous to that of Step 3 in the 
proof of Lemma \ref{disk-pert-r}. First of all from \eqqref(def-w-1-p) we see that
\[
\norm{\tilde w_1}_{L^{2,p}(\Delta)} \le a_p|\wect_0| + 
\norm{R^{(\nu)}\wect_0}_{L^{1,p}(\Delta)} \le (a_p + C\delta)|\wect_0|.
\]

Furthermore, from \eqqref(newton-1-p) we see that 
\[
\begin{split}
& \norm{\tilde w_{n+1}}_{L^{2,p}(\Delta)} \le C \norm{F^{(\nu)}(\zeta,\tilde w_n)}_{L^{1,p}
(\Delta)} + C\, |\wect_0| \le  
\\[2pt]
& \le C\, |\wect_0|  + C\delta \left(\norm{\tilde w_n}_{L^{2,p}(\Delta)} + 
\norm{\tilde w_n}_{L^p(\Delta)}^{\alpha}\right) \le 
\\[2pt]
& \le C\, |\wect_0|  + C\delta \norm{\tilde w_n}_{L^{2,p}(\Delta)} + 
C\delta |\wect_0|^{\alpha} \le 
\\[2pt]
& \le C\, |\wect_0| + C^2\delta\, |\wect_0| + (C\delta)^2\, \norm{\tilde w_{n-1}}_{L^{2,p}
(\Delta)}  + \left(C\delta + (C\delta)^2 \right) |\wect_0|^{\alpha}.
\end{split}
\]
By induction we prove that 
\[
\norm{\tilde w_{n+1}}_{L^{2,p}(\Delta)} \le C\, \frac{|\wect_0|}{1-C\delta} +
(C\delta)^n\norm{\tilde w_1}_{L^{2,p}(\Delta)} + C\delta\frac{|\wect_0|^{\alpha}}{1-C\delta}
\le C\,|\wect_0| + C\delta |\wect_0|^{\alpha},
\]
as stated. This gives us the a priori bound on $\norm{\tilde w_n}_{L^{2,p}(\Delta)}$, which was used for the previous estimates.

\smallskip\noindent{\slsf Step 8. $L^{2,p}$-convergence of 
approximations.}
Analogously to the estimate \eqqref(est-w-n) we have the following
\begin{equation}
\eqqno(l2p-w1-w_2)
\norm{\tilde w_{n+1}- \tilde w_n}_{L^{2,p}(\Delta)}\leq
C\eps\big(\norm{\tilde w_n -\tilde w_{n-1}}_{L^{2,p}(\Delta)} + 
\norm{\tilde w_n - \tilde w_{n-1}}^{\alpha}_{L^p(\Delta)}\big).
\end{equation}
Indeed, using  \eqqref(est4) we can write 
\[
\begin{split}
&\norm{\tilde w_{n+1}-\tilde w_n}_{L^{2, p}(\Delta)}\leq C\cdot
\norm{F^{(\nu)}(\zeta,\tilde w_n) - F^{(\nu)}(\zeta, \tilde 
w_{n-1})}_{L^{1,p}(\Delta)}\leq
\\[2pt]
& \le C\cdot \norm{J}_{\calc^{1,\alpha}(\Delta)}\; \norm{\tilde u_0}_{L^{2,p} \Delta)}
\left(\norm{\tilde w_n -\tilde w_{n-1}}_{L^{2,p}(\Delta)} + \norm{\tilde w_n - 
\tilde w_{n-1}}_{L^p(\Delta)}^{\alpha}\right) \le
\end{split}
\]
\begin{equation}
\eqqno(l2p-est)
\leq  C_1\cdot\eps \cdot \left(\norm{\tilde w_{n}- \tilde w_{n-1}}_{L^{2,p}
(\Delta)} +  \norm{\tilde w_n - \tilde w_{n-1}}_{L^p(\Delta)}^{\alpha}\right),
\end{equation}
where $\eps >0$ as small as we wish.  

\smallskip Notice that due to \eqqref(6.20) we have that 
\[
\norm{\tilde w_{n+1}-\tilde w_n}_{L^p(\Delta)}
\le C\, Lip_B(J)\norm{\tilde w_n - \tilde w_{n-1}}_{L^p(\Delta)}.
\]
Rescaling we can insure that $C\cdot Lip_B(J) <\eps$ to obtain from here that 
\[
\norm{\tilde w_{n+1}-\tilde w_n}_{L^p(\Delta)} \le \norm{\tilde w_2 - \tilde w_1}_{L^p(\Delta)}\cdot \eps^{n-1},
\]
and therefore 
\[
\norm{\tilde w_{n+1}- \tilde w_{n}}_{L^p(\Delta)}^{\alpha} \le
\norm{\tilde w_2 - \tilde w_1}_{L^p(\Delta)}^{\alpha}\cdot (\eps^{\alpha})^{n-1}.
\]

\medskip\noindent{\slsf Claim 1. } {\it  One has the following estimate:
\begin{equation}
\eqqno(eta-tau1)
\begin{split}
& \norm{\tilde w_{n+1} -\tilde w_n}_{L^{2,p}(\Delta)} \le (C_1\eps)^{n-1}
\norm{\tilde w_2 -\tilde w_1}_{L^{2,p}(\Delta)} + 
\\[2pt]
& + \norm{\tilde w_2 - \tilde w_1}_{L^p(\Delta)}^{\alpha}\, 
\sum_{k=0}^{n-2}\eps^{k\alpha}(C_1\eps)^{n-k-1},
\end{split}
\end{equation}
where $C_1$ is the constant from \eqqref(l2p-est).
} Indeed, for $n=2$ we have by \eqqref(l2p-est)
\[
\norm{\tilde w_3 -\tilde w_2}_{L^{2,p}(\Delta)} \le C_1\eps \norm{\tilde w_2 - 
\tilde w_1}_{L^{2,p}(\Delta)} + C_1\eps \norm{\tilde w_2 - \tilde w_1}_{L^p(\Delta)}^{\alpha},
\]
which is \eqqref(eta-tau1) in this case. After that we write
\[
\begin{split}
&\norm{\tilde w_{n+2} -\tilde w_{n+1}}_{L^{2,p}(\Delta)} \le 
C_1\cdot\eps \cdot \left(\norm{\tilde w_{n+1}- \tilde w_{n}}_{L^{2,p}
(\Delta)} +  \norm{\tilde w_{n+1} - \tilde w_{n}}_{L^p(\Delta)}^{\alpha}\right)\leq 
\\[2pt]
& \leq (C_1\eps)^{n}\norm{\tilde w_2 -\tilde w_1}_{L^{2,p}(\Delta)} + 
 C_1\eps \norm{\tilde w_2 - \tilde w_1}_{L^p(\Delta)}^{\alpha}\, (\eps^{\alpha})^{n-1}
+
\\[2pt]
& + C_1\eps \norm{\tilde w_2 - \tilde w_1}_{L^p(\Delta)}^{\alpha}\, \sum_{k=1}
^{n-2}\eps^{k\alpha}(C_1\eps)^{n-k-1} = (C_1\eps)^{n}\norm{\tilde w_2 -\tilde w_1}_{L^{2,p}(\Delta)} +
\\[2pt]
& + \norm{\tilde w_2 - \tilde w_1}_{L^p(\Delta)}^{\alpha} 
\sum_{k=1} ^{n-1}\eps^{k\alpha}(C_1\eps)^{n-k},
\end{split}
\]
as claimed.

\medskip\noindent{\slsf Claim 2.} {\it  Under the assumptions of the theorem the series
\[
\tilde w_N = \sum_{n=0}^{N}\left(\tilde w_{n+1} - \tilde w_{n}\right ) +\wect_0
\]
converge in $L^{2,p}$-norm to a solution $\tilde w$ of \eqqref(syst1-1). 
}

\smallskip Indeed,
\[
\sum_{n=2}^{\infty}\norm{\tilde w_{n+1} -\tilde w_n}_{L^{2,p}(\Delta)}
\le \frac{\norm{\tilde w_2 -\tilde w_1}_{L^{2,p}(\Delta)}}{1-C_1\eps} + 
\norm{\tilde w_2 -\tilde w_1}_{L^{p}(\Delta)}^{\alpha}\sum_{n=2}^{\infty}
\sum_{k=0}^{n-1}\eps^{k\alpha}(C_1\eps)^{n-k-1}.
\]
And the latter series converge. From \eqqref(est-w-n) we conclude that 
\begin{equation}
\eqqno(tilde-2p)
\norm{\tilde w}_{L^{2,p}(\Delta)} \le C_p\cdot |\wect_0| + C_p\eps|\wect_0|^{\alpha},
\end{equation}
which implies \eqqref(pert5). Lemma is proved.

\smallskip\qed

\begin{rema} \rm 
From \eqqref(tilde-2p) we obtain that for a given $p>2$ and  $\beta \deff 1-\frac{2}{p}$
one has 
\begin{equation}
\eqqno(tilde-beta)
\norm{\tilde w}_{\calc^{1,\beta}(\Delta)} \le C_p\cdot |\wect_0| + 
C_p\eps |\wect_0|^{\alpha}.
\end{equation}
In particular, for  $p= \frac{2}{1 - \alpha}$ we get 
\begin{equation}
\eqqno(tilde-alpha)
\norm{\tilde w}_{\calc^{1,\alpha}(\Delta)} \le C_p\cdot |\wect_0| +
C_p\eps |\wect_0|^{\alpha},
\end{equation}
and the same for $w = \tilde w|_{\Delta^+}$. Here as $\eps$ we can take 
$\norm{J}_{\calc^{1,\alpha}(B)}$, provided the latter is small enough.
\end{rema}

\newprg[COMP.no-cusp]{Smoothing cusps by perturbation}

We prove the following proposition. 

\begin{prop}
\label{no-cusps}
Let $J$ be an almost complex structure on $\rr^{2n}$, $J(0)$ $= J\st$ of class 
$\calc^{1,\alpha}$ and $u_0:(\Delta^+,\db_0\Delta^+,0) \to (\rr^{2n}, \rr^2, 0)$ 
a $J$-holomorphic map. Then the $J$-complex curve 
\begin{equation}
u(\zeta) = u_0(\zeta) + \zeta w, 
\end{equation}
constructed in Lemma \ref{cusp-pert-h} for $\nu =1$  has no cusps in a 
neighbourhood of zero, provided  the vector $\wect_0 = w(0)\in \rr^n$ is 
transverse to to the tangent vector $\vect_0\in \rr^n$ and is small enough. 
\end{prop}
\proof Proof is essentially that of Lemma 6.1 in \cite{IS4} applied to the extension 
$\tilde u_0$ of $u_0$ by reflection as the section of the bundle $E$. Function $w$ 
here is actually $w=\tilde w|_{\Delta^+_r}$, where $\tilde w$ is a solution of 
\eqqref(syst1-r) with $\nu =1$ and initial data $\tilde w(0)=\wect_0$. 

\smallskip 
Write $u_0$ in the normal form \eqqref(nrm-frm1) 
\begin{equation}
\eqqno(nrm-frm-u0)
u_0(\zeta) = \zeta^{\mu}P(\zeta) + \zeta^{2\mu -1}v(\zeta),
\end{equation}
where  $v\in L^{1,p}_{loc}$ with $v(0)=0$ and $P$ is a polynomial of degree
$\leq \mu -1$, $P(0) = \vect_0=e_1$. The differential of $\tilde u_0$ 
satisfies
\begin{equation}
\eqqno(nrm-frm-dif)
d\tilde u_0 (\zeta)  = \mu \vect_0\zeta^{\mu - 1} + O(|\zeta|^{\mu-1+\alpha}),
\end{equation}
and this implies the same relation for $u=\tilde u|_{\Delta^+}$. Rewrite this as 
\begin{equation}
\eqqno(du0)
\nabla \tilde u_0(\zeta) = \mu \zeta^{\mu -1}e_1 + R(\zeta),
\end{equation}
with the pointwise estimate
\begin{equation*}
|R(\zeta)| \le 
C_0\norm{\nabla \tilde u_0}_{L^p(\Delta)} \cdot  |\zeta|^{\mu -1 +\alpha} \le 
\sigma |\zeta|^{\mu -1 +\alpha}
\end{equation*}
for some $\sigma >0$ in the disc $\Delta(\half)$, this follows from \eqqref(est-v). 
Wlog we can assume that $\wect_0 =ae_2$, where $|a|$  can be taken as small as we need. 
Write 
\[
\nabla \tilde u = \nabla \tilde u_0 + \nabla(\zeta \tilde w)
\]
and use the estimate $\norm{\nabla w}_{L^{1,p}(\Delta)}\le C|\wect_0|$, which implies 
$|\nabla (\zeta w) - \wect_0| \le C_1|\zeta|^{\alpha}|\wect_0|$.

\smallskip We see that the $e_2$-component of $\nabla \tilde u$ will not vanish provided 
\begin{equation}
\eqqno(e2-estim)
C_1\cdot |\zeta|^{\alpha} \le \frac{1}{3} 
\qquad\text{and}\qquad
\sigma\cdot |\zeta|^{\mu-1+ \alpha} 
\le \frac{|a|}{3}.
\end{equation}
The $e_1$-component will not vanish provided 
\begin{equation}
\eqqno(e1-estim)
\sigma\cdot |\zeta|^{\alpha} \le \frac{\mu}{3} 
\qquad\text{and}\qquad
|a|\cdot C\cdot |\zeta|^{\alpha} \le \frac{1}{3} \mu |\zeta|^{\mu-1}.
\end{equation}
Since the constants $\sigma$ and $C$ are independent of $a$ we may assume 
that the first conditions in \eqqref(e1-estim) and \eqqref(e2-estim) are 
satisfied, namely
\begin{equation}
\eqqno(zeta<Auw)
|\zeta| \le \min \left( \left| \frac{1}{3C} \right|^{1/\alpha}, 
    \left| \frac{\mu}{3\sigma} \right|^{1/\alpha} \right).
\end{equation}
Therefore it remains to show that for any sufficiently small $\zeta$ (i.e., satisfying \eqqref(zeta<Auw))  \emph{at least one} of the following two conditions: 
\begin{equation}
    |\zeta|^{\mu-1+ \alpha} \le  \frac{|a|}{3\sigma} 
\qquad\text{or}\qquad
\frac{3|a|C}{\mu} \le  |\zeta|^{\mu-1- \alpha}
\end{equation}
is satisfied provided $|a|\not= 0$ is sufficiently small.
A sufficient condition for this is the following.
\begin{equation}
\left( \frac{3|a|C}{\mu} \right)^{\frac{1}{\mu-1 - \alpha}} \le
\left( \frac{|a|}{3\sigma} \right)^{\frac{1}{\mu-1 + \alpha}},
\end{equation}
which is equivalent to
\begin{equation}
\eqqno(a-condit)
|a|^{\frac{2\alpha}{(\mu -1)^2-\alpha^2}} = a^{\frac{1}{\mu-1 - \alpha} -\frac{1}
{\mu-1 + \alpha}}\le  \Big( \frac{\mu}{3C} \Big)^{\frac{1}{\mu-1 - \alpha}}
\Big( \frac{1}{3\sigma} \Big)^{\frac{1}{\mu-1 + \alpha}}.
\end{equation}
The latter holds for $|a|$ sufficiently small since $\mu \ge 2$ and $0<\alpha <1$.
Proposition and therefore  Theorem \ref{b-pert-thm} are proved.

\smallskip\qed

\ifx\undefined\bysame
\newcommand{\bysame}{\leavevmode\hbox to3em{\hrulefill}\,}
\fi

\def\entry#1#2#3#4\par{\bibitem[#1]{#1}
{\textsc{#2 }}{\sl{#3} }#4\par

\end{document}